# Interior point methods for an algebraic system involving complementarity equations for geomechanical fractures

A dissertation submitted by
HOANG Trung Hau

Master Thesis Report

based on a work conducted at

INSTITUT NATIONAL DES SCIENCES APPLIQUÉES RENNES

under the supervision of:

Professor Mounir Haddou: University professor, INSA-IRMAR Rennes
Ibtihel BEN GHARBIA: Research engineer, IFPEN
TRAN Quang Huy: Research engineer, HDR, IFPEN

June 21, 2019

# Acknowledgement

In this page, I want to express my gratitude to my supervisors **Professor Mounir Haddou, Research Engineer Ibtihel BEN GHARBIA, and Research Engineer TRAN Quang Huy**. Thank you for accepting me as your internship. I am very grateful to **Professor Mounir Haddou** who welcomed me warmly in Rennes. Thank you for the first step to help me get acquainted with the INSA work environment and the life in Rennes.

I would like to thank the enthusiastic and vast knowledge of **Professor Mounir Haddou, Research Engineer Ibtihel BEN GHARBIA, and Research Engineer TRAN Quang Huy**. I have learned a lot from this Internship, It's really opening my eyes bigger in Optimization in particular and in science in general, Thank you for give me really good questions, give me the question at the bottom of things, give me hard problems, comment, pointed out my mistake during I solve the problems and answer my questions for me to understand that I have many shortcomings and I realized I needed to try harder, It also helps me to grow up so much. Especially, I feel warm when my question is answered clearly and easy to understand, and help me to see with an expert perspective. So that I feel myself stronger, and and I feel more confident about myself.

I would like to thank **Ms. Xhensila Lachambre**& **Ms. Patricia Soufflet**& **Ms.** Research Engineer& **Ms. Minh Phuong** for helping me with the paperwork when I was in France and Vietnam. I would like to thank **Henri Lebesgue Center, INSA, and IRMAR and It's partners** which allowed me to spend my Internship in France and to discover this beautiful country both scientifically and culturally.

I would like to thank my parents, and all my family and my colleagues and my friends for their constant support, encourage,... It's my great motivation to help me complete this dissertation.

In the end, I would like to say it is my great honor of me to be the student of **Professor Mounir Haddou, Research Engineer Ibtihel BEN GHARBIA, and Research Engineer TRAN Quang Huy**.

*At Rennes, France, June 21, 2019*

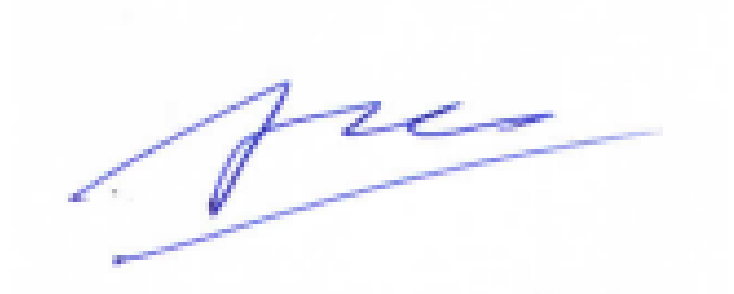

*HOANG Trung Hau*

# Contents

# Chapter 1

# Introduction

## 1.1 Origin of the problem

- Many applications like subseismic fault modeling, fractured reservoir modeling and interpretation/validation of fault connectivity involve the solution to an elliptic boundary value problem in a background medium perturbed by the presence of cracks that take the form of one or many pieces of surface (with boundary).
- When the background medium can be considered as homogeneous, boundary integral equations appear as a method of choice for the numerical solution to fractures problems.
- With such an approach, the problem is reformulated as a fully non-local equation posed at the surface of cracks. This is the strategy adopted by IFP Energies Nouvelles (IFPEN) for the evolution of the deformation and perturbed stress field associated with the solution of an elastostatic problem around a network composed of multiple cracks surfaces
- Discretization of boundary integral resulting in the so-called Boundary Element Method (BEM) [4] leads to densely populated matrices due to the full non-locality of the operators under consideration. After the discretization process, geologists are faced with a system of equations that turns out difficult to solve numerically.

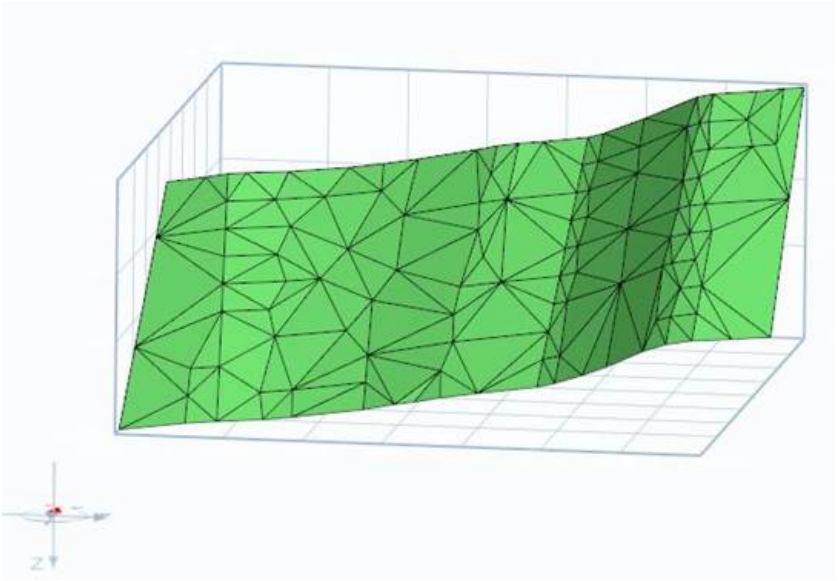

Figure 1.1: An example of cracks network with N = 500

## 1.2 Statement of the problem

- To describe this system, let u $\in \mathbb{R}^{3N}$ be the vector of unknowns (displacements), where $N \geqslant 1$ is the number of cells in the surface. This vector **u** consists of $N$ three-dimensional vectors $\mathbf{v}_l = (u_{3I-2}, u_{3I-1}, u_{3I}) \in \mathbb{R}^3$, for $1 \leqslant I \leqslant N$, each representing the displacement of the *I-th* cell.
- Let
  - $\mathbf{A} \in \mathbb{R}^{3N\times 3N}$ (generated by BEM), encapsulating the interactions between the cells.
  - $\mathbf{b} \in \mathbb{R}^{3N}$ a (known) vector corresponding to some external force.

we want to find $\mathbf{u} \in \mathbb{R}^{3N}$, solution of

$$\begin{cases} (Au-b)_1 & = 0 \\ (Au-b)_2 & = 0 \\ min\,\{u_3; (Au-b)_3\} & = 0 \\ (Au-b)_4 & = 0 \\ (Au-b)_5 & = 0 \\ min\,\{u_6; (Au-b)_6\} & = 0 \\ \dots & \\ \dots & \\ \dots & \\ (Au-b)_{3N-2} & = 0 \\ (Au-b)_{3N-1} & = 0 \\ min\,\{u_{3N}; (Au-b)_{3N}\} & = 0 \end{cases} \tag{1.1}$$

where $(Au-b)_i$ denotes the *i-th* component of $Au-b$. In other words,

$$(Au-b)_i = \sum_{j=1}^{3N} A_{ij}u_j - b_i \quad \forall 1 \leqslant i \leqslant 3N \tag{1.2}$$

**Remark 1.1.** *The system (1.1) can be seen as a modification of the linear system **Au=b**, in which every third equation has been replaced by a complementarity equation. Physically speaking, the latter expresses a "non-interpenetration" condition to be imposed on the normal component of each displacement vector* $(u_{3I-2}, u_{3I-1}, u_{3I}) \in \mathbb{R}^3$ *corresponding to the I-th cell.*

Loosely speaking, the previous set of equations can be regarded as a "constrained" version of the linear system $\mathbf{Au} = \mathbf{b}$ (the mechanical behavior of the fractures is modeled as a contact problem).

Many empirical algorithms have been proposed [2],[3] by geologists to solve this system of equations. Most of them are iterative methods based on some algebraic reformulation of the system at issue. Unfortunately, none of them is guaranteed to converge in theory (in particular when faults (fractures) intersect each other forming a geometrically highly irregular structure). In practice, none of them appears to be either robust or efficient (i.e., to run with reasonable computational time).

The objective of this project is to investigate another approach, referred to as interior point methods, for which convergence can be ensured (even if faults are too close). Interior point methods have proved their efficiency in a wide variety of domains, most notably for linear programming. Here, even though we do not have any optimization problem, we can adapt ideas from interior point methods for the numerical resolution of the system considered.

## 1.3 Standard method

Currently, in the **IFPEN** geomechanical software, the following Jacobi-like (It looks like the Jacobi method for the linear system) iterative method is implemented.

**Algorithm 1** Standard Method

1: Starting from an initial guess $\mathbf{u}^{(\mathbf{0})}$
2: Update

$$\begin{aligned} A_{3I-2,3I-2}u_{3I-2}^{(k+1)} &= b_{3I-2} - \sum_{j=1, j\neq 3I-2}^{3N} A_{3I-2,j}u_j^{(k)} \\ A_{3I-1,3I-1}u_{3I-1}^{(k+1)} &= b_{3I-1} - \sum_{j=1, j\neq 3I-1}^{3N} A_{3I-2,j}u_j^{(k)} \\ A_{3I,3I}u_{3I}^{(k+1)} &= \left( b_{3I-2} - \sum_{j=1, j\neq 3I}^{3N} A_{3I-2,j}u_j^{(k)} \right)^+ \end{aligned} \tag{1.3}$$

3: If $F(\mathbf{u}^{(\mathbf{k})}) = \left\| u^{(k+1)} - u^{(k)} \right\| = 0$; stop; else back to step 2

for $1 \leqslant I \leqslant N$. In (1.3) the symbol

$$a^+ = max\{a, 0\}$$

denoted the positive part of any real number a.

**Theorem 1.1.** *Assume that*

$$A_{i,i} > 0 \quad \forall 1 \leqslant i \leqslant 3N$$

*Then the iteration (3) are well defined. If the iterates $u^{(k)}$ converge to a limit $u$, then $u$ is a solution of* (1.1)*.*

*Proof.* Firstly, we talk about the iteration (3) are well defined . Because the left hand side is $A_{i,i}u_i$ where $A_{i,i} > 0$ then the iteration can go from step $k \to k+1$ without any trouble for all I s.t $1 \leqslant I \leqslant 3N$

Secondly, we talk about the remained part. we can suppose that

$$u^{(k)} = \begin{pmatrix} u_1^{(k)} \\ u_2^{(k)} \\ u_3^{(k)} \\ u_4^{(k)} \\ \cdots \\ u_{3N}^{(k)} \end{pmatrix} \to \begin{pmatrix} u_1 \\ u_2 \\ u_3 \\ u_4 \\ \cdots \\ u_{3N} \end{pmatrix} = u$$

so that

$$\begin{pmatrix} u_1^{(k)} \to u_1 \\ u_2^{(k)} \to u_2 \\ u_3^{(k)} \to u_3 \\ u_4^{(k)} \to u_4 \\ ... \\ u_{3N}^{(k)} \to u_{3N} \end{pmatrix}$$

Consider the index 3I-2 we have

$$\lim_{k\to\infty} u_{3I-2}^{(k)} = u_{3I-2}$$

Or

$$\lim_{k\to\infty} b_{3I-2} - \sum_{j=1, j\neq 3I-2}^{3N} A_{3I-2,j} u_j^{(k-1)} = A_{3I-2,3I-2} u_{3I-2}$$

Or

$$b_{3I-2} - \sum_{j=1, j\neq 3I-2}^{3N} A_{3I-2,j} u_j = A_{3I-2,3I-2} u_{3I-2}$$

Rearranging we obtain

$$\sum_{j=1}^{3N} A_{3I-2,j} u_j - b_{3I-2} = 0$$

Or

$$(Au - b)_{3I-2} = 0$$

The same with the index 3I-1

Now we go to the equation

$$min\,\{u_{3I}, (Au - b)_{3I}\} = 0 \quad \forall I = 1, ..., N \tag{1.4}$$

we have that

$$\lim_{k\to\infty} u_{3I}^{(k)} = u_{3I}$$

Or

$$\lim_{k\to\infty} (b_{3I-2} - \sum_{j=1, j\neq 3I}^{3N} A_{3I,j} u_j^{(k-1)})^+ = A_{3I,3I} u_{3I}$$

By the continuity of the function $x \to x^+$ we get

$$(b_{3I-2} - \sum_{j=1, j\neq 3I}^{3N} A_{3I,j} u_j)^+ = A_{3I,3I} u_{3I} \tag{1.5}$$

Now if

$$(b_{3I-2} - \sum_{j=1, j\neq 3I}^{3N} A_{3I,j} u_j)^+ > 0$$

then (1.5) become :

$$(b_{3I-2} - \sum_{j=1, j\neq 3I}^{3N} A_{3I,j}u_j) = A_{3I,3I}u_{3I}$$

Rearranging we get

$$\sum_{j=1}^{3N} A_{3I,j}u_j - b_{3I-2} = 0$$

Or

$$(Au - b)_{3I} = 0$$

By the definition of $u_{3I}^{(k)} \geqslant 0$ we conclude that $u_{3I} \geqslant 0$ so

$$min\{u_{3I}, (Au - b)_{3I}\} = 0 \quad \forall I = 1, ..., N$$

Now if

$$(b_{3I-2} - \sum_{j=1, j\neq 3I}^{3N} A_{3I,j}u_j)^{+} \leqslant 0$$

then

$$A_{3I,3I}u_{3I} = 0$$

By $A_{3I,3I} > 0$ then

$$u_{3I} = 0$$

but

$$Au - b)_{3I} = \sum_{j=1}^{3N} A_{3I,j}u_j - b_{3I-2} \geqslant 0$$

then

$$min\{u_{3I}, (Au - b)_{3I}\} = 0 \quad \forall I = 1, ..., N$$

and the result follows. □

**Remark 1.2.** *The way we choose the initial point is very important when we choose a good initial point we can make the first error is as small as possible. So from this, we can save time, save iteration and memory also (save memory for the case N is large our algorithm still work). One idea is as follow.*

*(i) Let D is a diagonal matrix of A. We consider D as an approximation of A*

*(ii) then we have the system to solve is*

$$\begin{cases} (Du - b)_1 & = 0 \\ (Du - b)_2 & = 0 \\ min\{u_3, (Du - b)_3\} & = 0 \\ (Du - b)_4 & = 0 \\ (Du - b)_5 & = 0 \\ min\{u_6, (Du - b)_6\} & = 0 \\ \cdots \\ \cdots \\ \cdots \\ (Du - b)_{3N-2} & = 0 \\ (Du - b)_{3N-1} & = 0 \\ min\{u_{3N}, (Du - b)_{3N}\} & = 0 \end{cases} \tag{1.6}$$

*(iii) the third equation is equivalent to*

$$0 \leqslant u_{3I} \perp (Du - b)_{3I} \geqslant 0 \tag{1.7}$$

*for* $1 \leqslant I \leqslant N$ *, so we need an approximation of* $b$ *(denoted* $\tilde{b}$ *) to ensure that* $(Du - b)_{3I} \geqslant 0$ *. we consider an approximation*

$$\tilde{b} = \begin{bmatrix} b_1 \\ b_2 \\ b_3^+ \\ \cdots \\ \cdots \\ b_{3N-2} \\ b_{3N-1} \\ b_{3N}^+ \end{bmatrix} \tag{1.8}$$

*So now if we solve the system* $Du = \tilde{b}$ *so* $u = D \backslash \tilde{b}$ *.* $u$ *is satisfies the system* (1.6) *. With the ways choosing* $u$ *above .*

### 1.3.1 Numerical results of Standard Method

- 3 algorithms here (Standard method, IPM, non-parametric IPM) are implemented on a standard laptop (2.5 GHz, 2 GoM) in Matlab R 2018 an update 2
- using command \ in Matlab to compute Newton direction
- $It_{max}$ is set to 40

**The parameter Number of restart meaning the number of initial point, e.g. Consider the picture below**

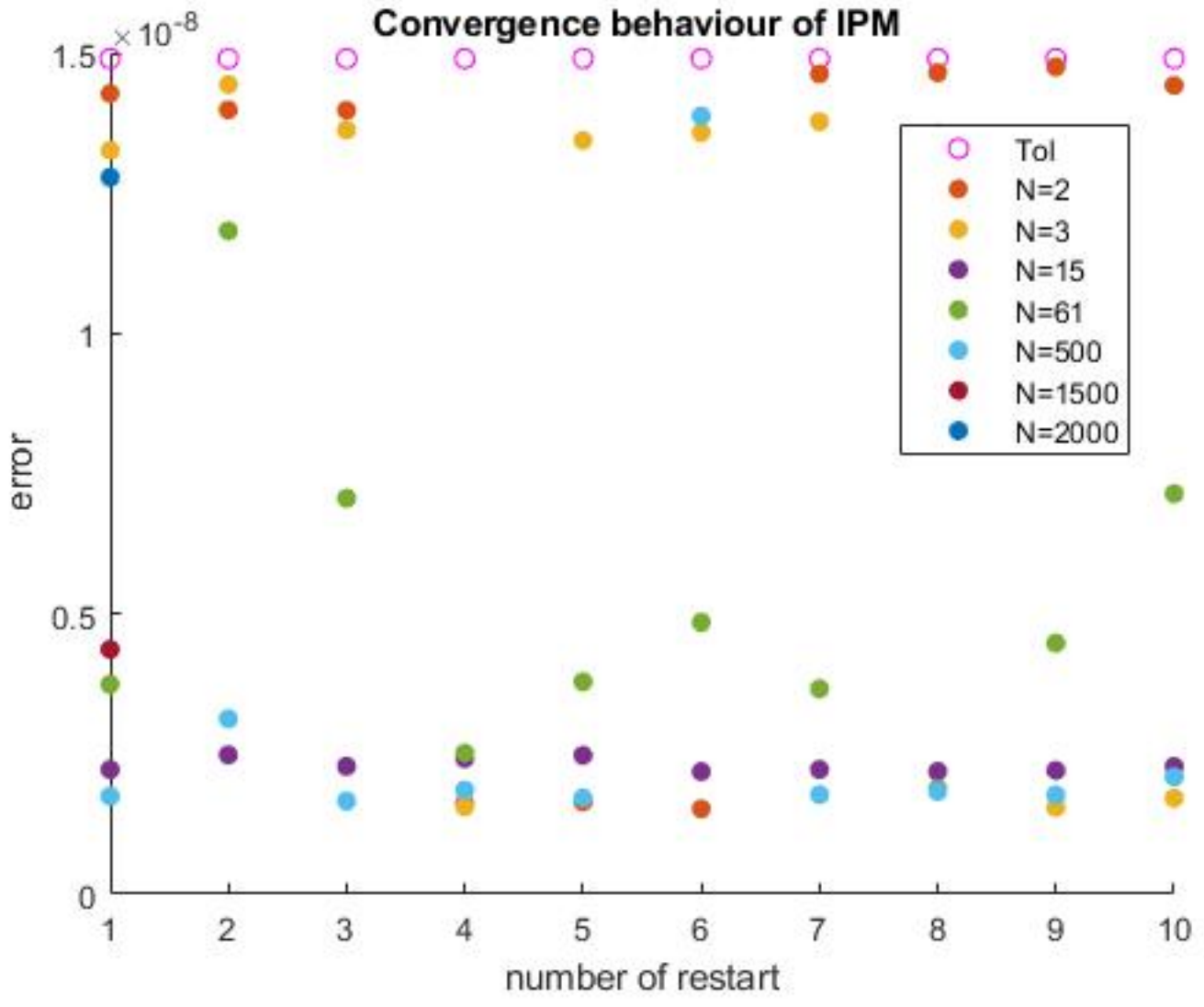

Figure 1.2: The convergence of IPM for various option N

**Look at the axis x, we can see that the number of restart = 1,2,3,...,10. It means the first initial point, the second initial point, the third initial point,...,**

**the tenth initial point. Look at the end of the axis x, the number of restart = 10, It means we use 10 difference initial point. Number of restarts also has this meaning throughout this thesis**.

Figure below show in cases N = 61,500,1500,2000 . **"Standard Method"** is diverge.

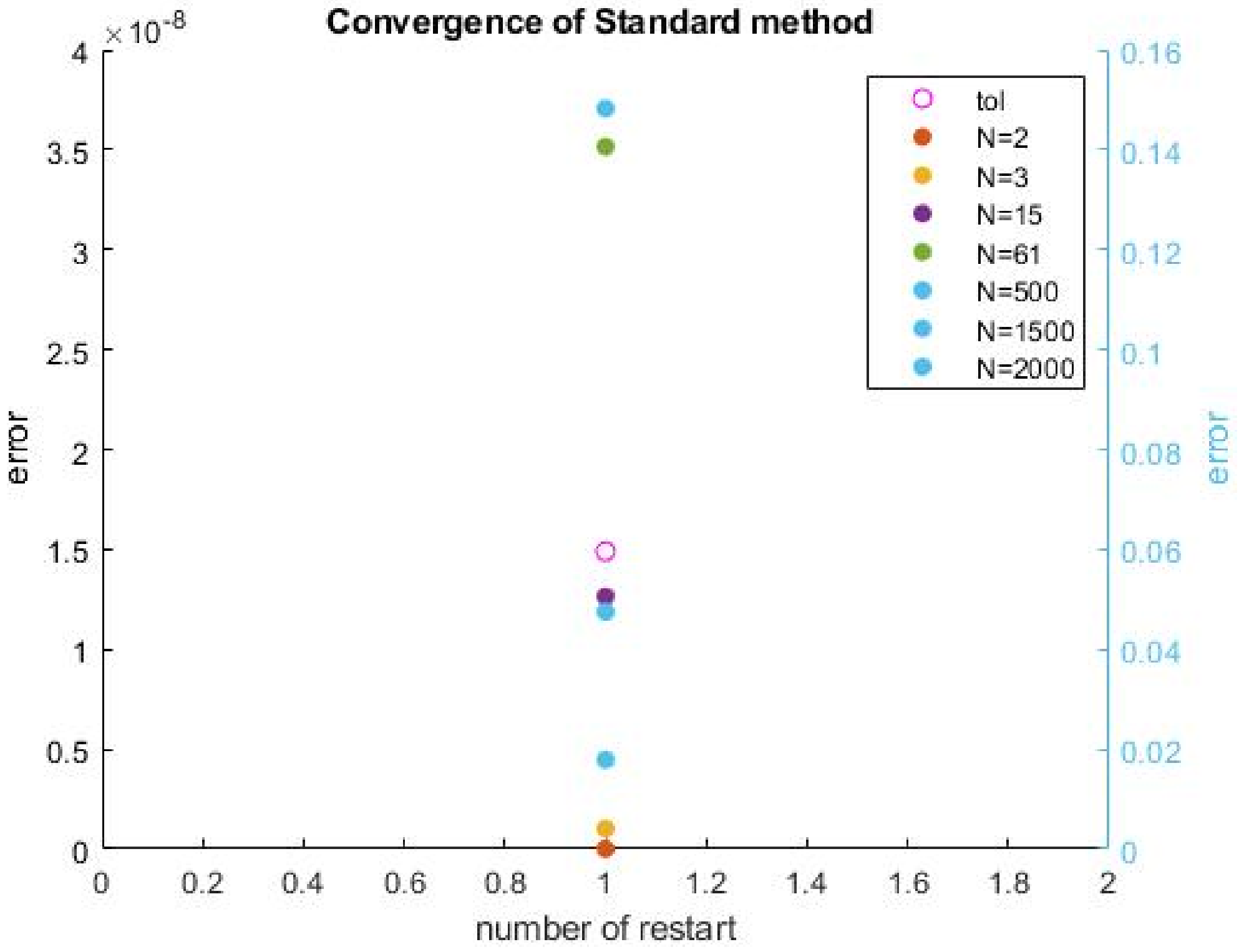

Figure 1.3: The convergence of Standard method for the various option N

**Remark 1.3.** *3 blue dots on the picture* (1.3) *(N=500, N=1500, N=2000), follow the axis of error on the right, the others colored dots follow the axis of error on the left.*

Secondly, we propose to how many iterations we pay for convergence with **"Standard Method"**

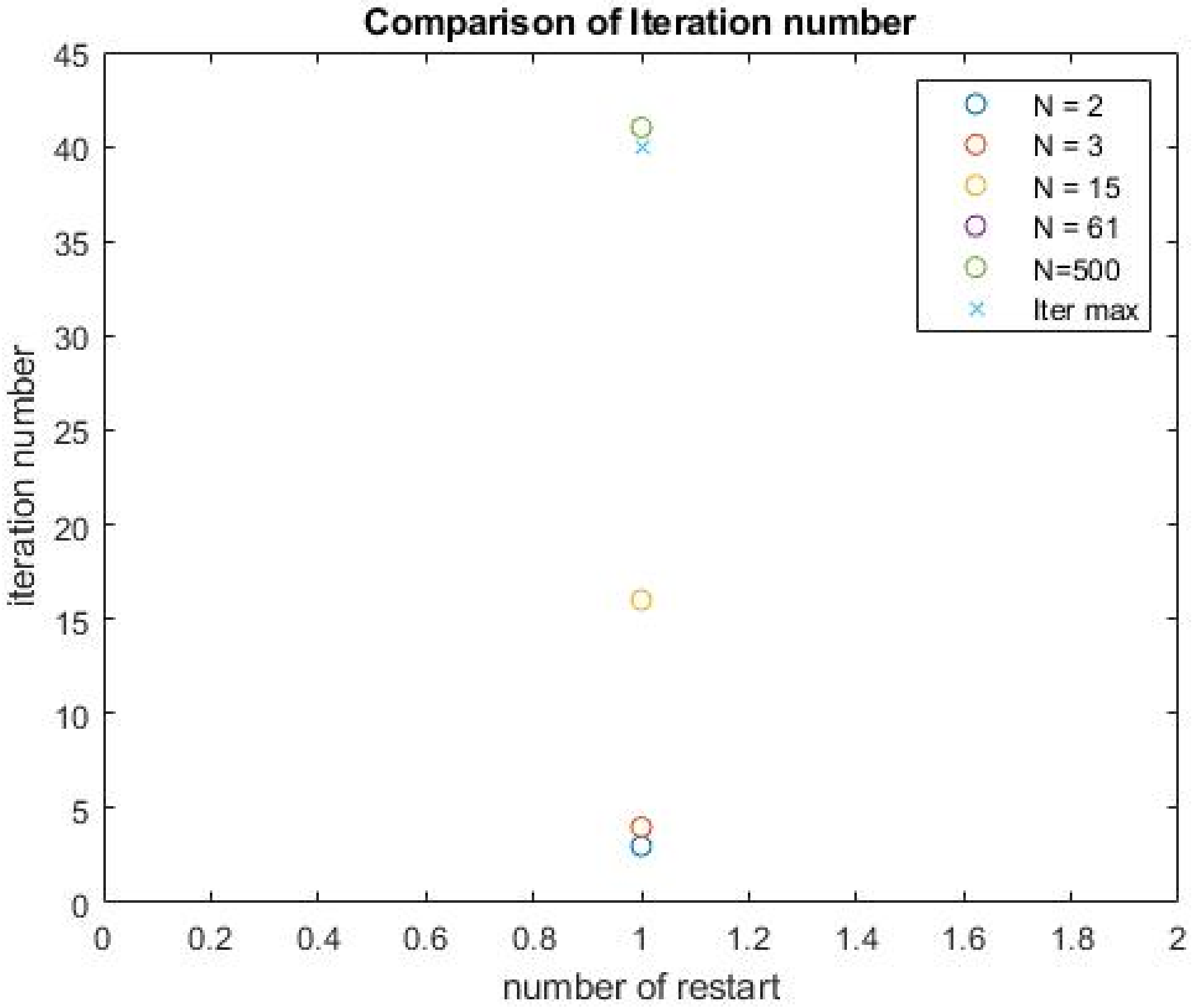

Figure 1.4: Iteration of Standard methods for the various option N

**Remark 1.4.** *The trouble with Standard method is that it does not always converge. In other words, Standard method is lack of robustness, slowness... We wish to propose a more robust one.*

## 1.4 Interior Point Method

The idea is to consider a sequence of approximate problems $\mathrm{P}_\mu$ defined by

$$(\mathrm{P}_\mu)\begin{cases}(Au-b)_{3I-2}=\mathbf{0}\\(Au-b)_{3I-1}=\mathbf{0}\\u_{3I}\bullet(Au-b)_{3I}=\mu\mathbf{e}\end{cases}\tag{1.9}$$

for $1\leqslant I\leqslant N$ and to drive $\mu>0$ to 0 in some "smart" way.

**Remark 1.5.** *The reason why we need to consider a "Perturbed system" because with the original ones we can't solve it. Illustration as below*

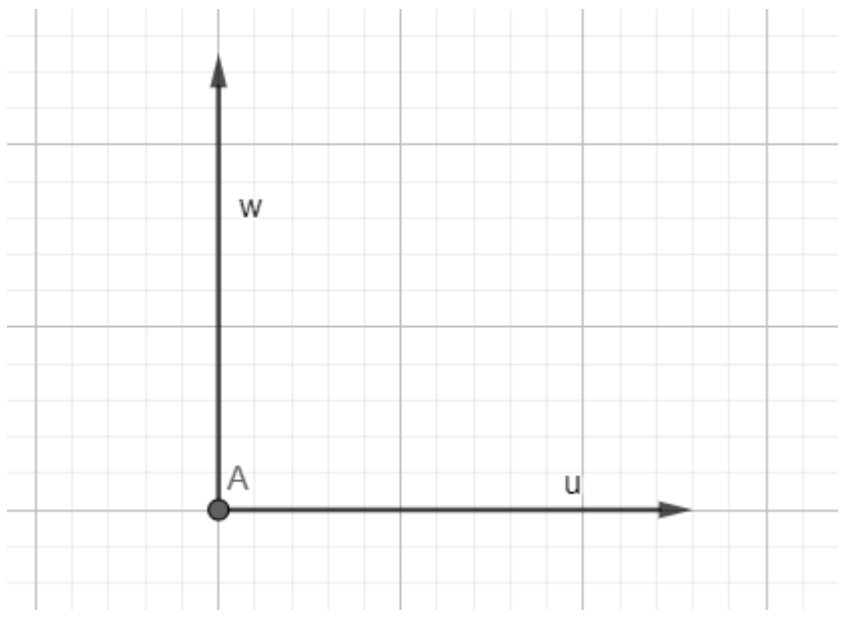

Figure 1.5: Picture of the original complementary problem

*So if we are looking at the picture, then 2 axis is perpendicular, we can't slide or move. Now if we look at the "Perturbed system".*

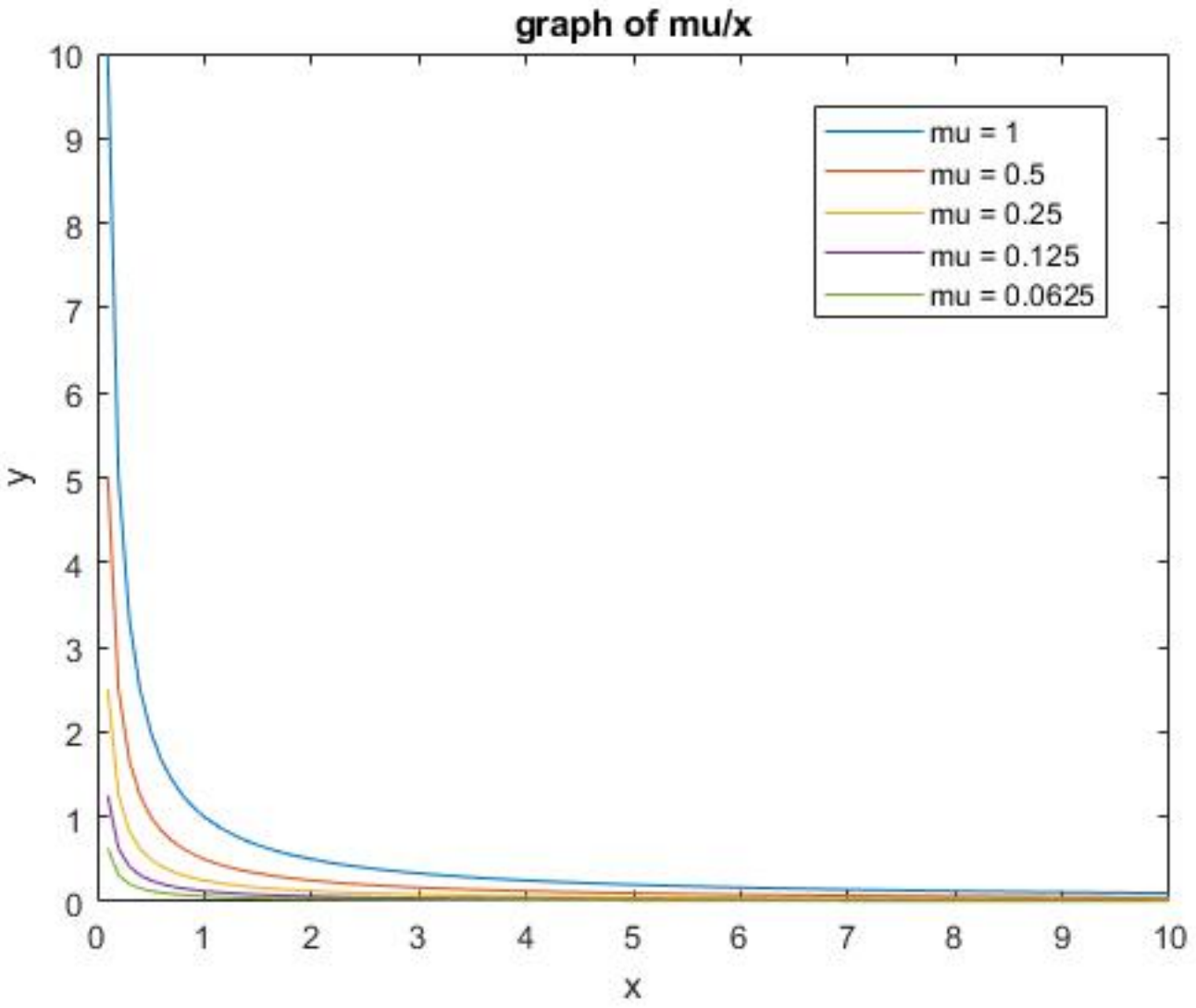

Figure 1.6: Picture of the relaxation complementary problem

*which make the problems easier, we can slide on it, which maybe we can easier to obtain the solution. If* $\mu\downarrow 0$ *then we get a good approximation of original problems*

To make things easier, let us introduce the slack variables.

$$w = Au - b \in \mathbb{R}^{3N} \tag{1.10}$$

system (1.9) then reads

$$\begin{cases} (Au - b - w)_{3I-2} = \mathbf{0} \\ (Au - b - w)_{3I-1} = \mathbf{0} \\ (Au - b - w)_{3I} = \mathbf{0} \\ w_{3I-2} = \mathbf{0} \\ w_{3I-1} = \mathbf{0} \\ u_{3I} \bullet w_I = \mu \mathbf{e} \end{cases} \tag{1.11}$$

with 6N unknowns and 6N equations for a fixed $\mu$ for $1 \leqslant I \leqslant N$.

**Remark 1.6.** *the system above consist of 6N unknown which makes the size of the Jacobi matrix is large. Another approach in the next chapter will consider only 4N unknown. So we can reduce the size of the Jacobi matrix which makes the time of convergence is faster and we can save the memory on the computer also the iteration number. Some test is runnable on the case 4N but 6N which make an improvement to run when N is large.*

Thus , a first IPM algorithm could be

**Algorithm 2** IPM algorithm

1: Starting from an initial guess $\mu^{(0)} > 0, (\mathbf{u}^{(0)}/(u_{3I}^{(0)}, w_I^{(0)}) > 0)$
2: If $F(\mathbf{u}^{(k)}) = 0$ or $DF(\mathbf{u}^{(k)})$ is not invertible: stop; else
3: Calculate $d^k$ by solving
$$F(\mathbf{u}^k) + DF(\mathbf{u}^k)d^k = 0$$
4: Truncate $(u_{3I}^{(k+1)}, w_I^{(k+1)})$ to ensure $u_{3I}^{(k+1)} \geqslant 0$ and $w_I^{(k+1)} \geqslant 0$
5: Update $\mathbf{u}^{k+1} = \mathbf{u}^k + \alpha_{max} d^k$, where $\alpha_{max}$ is found by step 4
6: Update $\mu^{(k+1)} = min \left\{0.8\mu^{(k)}, [\mu^{(k)}]^2\right\}$

**Remark 1.7.** *The idea to truncate is we want to find d at position is multiply by 3, for example, 3,6,9,12,... If d at this position is negative then we have to solve the system.*

$$x_- + \alpha d_- = 0 \tag{1.12}$$

*to ensure that $x_-$ is the position correspond to our position d is non-negative. So the $\alpha$ we get is*

$$\alpha = -\frac{x_-}{d_-} \tag{1.13}$$

*that means if we move with step size $> \alpha$ then we get $x_-$ is negative and we false. Now we have a large number of the equations like this, I mean the number of alpha so the alpha we take is*

$$\alpha_{max} = min(-x(ind)./d(ind)) \tag{1.14}$$

*where ind is the vector as follow*

$$ind = 3:3:6N \tag{1.15}$$

*that means the vector of position multiply by 3 in our vector u and w. But we need to make a safety parameter (let call it $\theta$) to ensure that we don't touch the boundary. Because if we take the boundary, our Jacobi matrix is near to singular (this is the bad case and we can't compute the Newton direction and then the iteration will stop). That is we choose the $\theta$ is positive and $< 1$ . Now the new $\alpha_{max}$ should be*

$$\alpha_{max} = \theta * \alpha_{max} \tag{1.16}$$

*In practice, we need to make a comparison between the choice of $\theta$ to get the good $\theta$ to make our algorithm converge with less iteration . Also, some case our $\alpha_{max}$ is $> 1$ then we need to make a correction with $\alpha_{max}$*

$$\alpha_{max} = min(\alpha_{max}, 1) \tag{1.17}$$

***(the reason why we need $\alpha_{max} \leqslant 1$ because if $\alpha_{max} > 1$ then we lose all theorical results about convergence of Newton methods)*** *to ensure we don't go to the negative zone . Now we can use this step size to update our newton direction .*

**Remark 1.8.** *we prove that*

$$min\left(1, \frac{x_1}{y_1}, \frac{x_2}{y_2}, ..., \frac{x_n}{y_n}\right) max\left(1, \frac{y_1}{x_1}, \frac{y_2}{x_2}, ..., \frac{y_n}{x_n}\right) = 1 \tag{1.18}$$

*we can see 1 = $\frac{x_0}{y_0}$ and suppose that $\exists\ i \in \{1, ..., n\}$ s.t*

$$min\left(1, \frac{x_1}{y_1}, \frac{x_2}{y_2}, ..., \frac{x_n}{y_n}\right) = \frac{x_i}{y_i}$$

*then*

$$\left(\frac{x_j}{y_j} \geqslant \frac{x_i}{y_i} \forall j \in \{1, ..., n\}\right) \rightarrow \left(\frac{y_j}{x_j} \leqslant \frac{y_i}{x_i} \forall j \in \{1, ..., n\}\right)$$

*so*

$$max\left(1, \frac{x_1}{y_1}, \frac{x_2}{y_2}, ..., \frac{x_n}{y_n}\right) = \frac{y_i}{x_i}$$

*so*

$$min\left(1, \frac{x_1}{y_1}, \frac{x_2}{y_2}, ..., \frac{x_n}{y_n}\right) max\left(1, \frac{y_1}{x_1}, \frac{y_2}{x_2}, ..., \frac{y_n}{x_n}\right) = \frac{x_i}{y_i}\frac{y_i}{x_i} = 1$$

*and the result follows. We prove that because we will use the step size*

$$\alpha_{max} = min(\alpha_{max}, 1) = \frac{1}{max(1, \frac{1}{\alpha_{max}})} \tag{1.19}$$

*In theory, it is the same, but in practice, the latter ones will more stable because when $d_i$ is too small then $\frac{-x_i}{d_i} >> 0$ is too large and we make some numerical error when divide something is too small. To avoid this problem we consider the max with type $-\frac{d_i}{x_i}$ when we scale $\alpha_{max} = 0.9 * \alpha_{max}$ that means we can guarantee that $x_i$ is not too near 0, it's small but not too small so that we can divide to avoid numeric problems.*

**Ex. 1** — 1. Write down the Jacobian matrix and the Newton iteration for (1.11)

**Answer (Ex. 1)** — If we rearrange $(Au-b-w)_{3I-2}, (Au-b-w)_{3I-1}, (Au-b-w)_{3I}$ on top, I mean the way we arrange the system (1.11) as follow:

$$\begin{cases} (Au-b-w)_1 & = 0 \\ (Au-b-w)_2 & = 0 \\ (Au-b-w)_3 & = 0 \\ (Au-b-w)_4 & = 0 \\ (Au-b-w)_5 & = 0 \\ (Au-b-w)_6 & = 0 \\ \quad ... & \\ \quad ... & \\ (Au-b-w)_{3N-2} & = 0 \\ (Au-b-w)_{3N-1} & = 0 \\ (Au-b-w)_{3N} & = 0 \\ w_1 & = 0 \\ w_2 & = 0 \\ u_3 w_3 & = \mu \\ w_4 & = 0 \\ w_5 & = 0 \\ u_6 w_6 & = \mu \\ ... & \\ ... & \\ w_{3N-2} & = 0 \\ w_{3N-1} & = 0 \\ u_{3N} w_{3N} & = \mu \end{cases} \tag{1.20}$$

then we get the Jacobian as follow :

$$DF(X) = \begin{bmatrix} A & -I \\ C & D \end{bmatrix} \in \mathbb{R}^{2m,2m} \tag{1.21}$$

where N is the number of elements , m = 3N , $A \in \mathbb{R}^{m,m}$ , $-I \in \mathbb{R}^{m,m}$, $-I$ is the identity matrix .

$$C = \begin{bmatrix} 0 & 0 & 0 & ... & 0 & 0 \\ 0 & 0 & 0 & ... & 0 & 0 \\ 0 & 0 & w_3 & 0 & ... & 0 \\ ... & & & & & \\ ... & & & & & \\ ... & & & & & \\ 0 & 0 & ... & 0 & 0 & 0 \\ 0 & 0 & ... & 0 & 0 & 0 \\ 0 & 0 & ... & 0 & 0 & w_{3N} \end{bmatrix} \in \mathbb{R}^{m,m}, D = \begin{bmatrix} 1 & 0 & 0 & ... & 0 & 0 \\ 0 & 1 & 0 & ... & 0 & 0 \\ 0 & 0 & u_3 & 0 & ... & 0 \\ ... & & & & & \\ ... & & & & & \\ ... & & & & & \\ 0 & 0 & ... & 1 & 0 & 0 \\ 0 & 0 & ... & 0 & 1 & 0 \\ 0 & 0 & ... & 0 & 0 & u_{3N} \end{bmatrix} \in \mathbb{R}^{m,m} \tag{1.22}$$

we can use \ in Matlab to compute Newton direction .
Newton iteration , firstly we choose $u_{3I}^{(0)}, w_{3I}^{(0)}, \mu^{(0)} > 0$

.we have the system to solve is F(x) = 0

.compute the newton direction : $d_N(x) = -DF(x)^{-1}F(x)$

.compute $\alpha_{max}$

.update $x^{(k+1)} = x^{(k)} + \alpha_{max} d_N(x)$

### 1.4.1 Mathematical background

Let E and F is 2 norm space with corresponding norm $\|.\|_E$,$\|.\|_F$, U is an open set in E, and e is a vector in E and x in U.

**Proposition 1.1.** *There exists a real number $\alpha$ s.t the open interval $(-\alpha, \alpha)$ included in this set*

$$I_{x,h} = \{t : t \in \mathbb{R}, x + te \in U\}$$

*set* $U_{x,h} = \{y : t \in \mathbb{R}, y = x + te \in U\}$

**Definition 1.1.** *Let f be a map from U to F, e be a vector in E and $x \in U$, we say that f has a directional derivative at x if and only if there is a bounded linear mapping T from E into G s.t*

$$T(h) = \lim_{t \to 0} \frac{f(x + th) - f(x)}{t} \quad \forall h \in E$$

*in this case, we call T the directional derivative at x of f and denote it by $Df(x)$ If Df(x) exists for any x in U, we say f is directional differentiable on U.*

**Definition 1.2.** *Let f be a map from U to F , $x \in U$ . we say that*

*. f is Gâteaux differentiable at x if f has directional derivative at x and Df(x) $\in$ L(E,F)*

*. f is Gâteaux differentiable on U if f has Gâteaux differentiable for any $x \in$ U*

**Theorem 1.2.** *Let f be a map be Gâteaux differentiable from an open set U in a norm space E into a norm space F. Suppose that the map x$\to$ Df(x) is a continuous mapping from U into L(E,F). Then f is Frechet differentiable on U*

**Definition 1.3.** *Let f be a map from U into F . we say that*

1. *f Gâteaux continuously differentiable if f Gâteaux differentiable on U and mapping $x \to Df(x)$ continuous from U into L(E,F)*

2. *f Frechet continuously differentiable if f Gâteaux differentiable on U and mapping $x \to Df(x)$ continuous from U into L(E,F). In this case we say that f is of class $C^1$*

### 1.4.2 Globalization with line search

In this section, we restrict X is a Hilbert space. A function F of class $C^1 : X \to X$ generates potential of least squares.

$$\Phi(x) := \frac{1}{2} \|F(x)\|^2 \tag{1.23}$$

If there exists a zero $\bar{x}$ of F, then inf$\Phi = 0$ and $\bar{x}$ is a solution of minimization problem of $\Phi$

**Lemma 1.1.** *function $\Phi$ is in class $C^1$ and $D\Phi(x)d = (F(x), DF(x)d)$. In particular, if $x \in X$ s.t $DF(x)$ is invertible, then*

$$D\Phi(x)d_N(x) = -2\Phi(x) \tag{1.24}$$

*Proof.* We begin with the directional derivative of $\Phi$

$$\begin{aligned}
D\Phi(x)(d) &= \lim_{t\to 0} \frac{\Phi(x+td)-\Phi(x)}{t} \\
&= \frac{1}{2}\lim_{t\to 0} \frac{(F(x+td),F(x+td))-(F(x),F(x))}{t} \\
&= \frac{1}{2}\lim_{t\to 0} \frac{(F(x+td),F(x+td))-(F(x+td),F(x))+(F(x+td),F(x))-(F(x),F(x))}{t} \\
&= \frac{1}{2}\lim_{t\to 0} \frac{(F(x+td),F(x+td)-F(x))}{t} + \lim_{t\to 0} \frac{(F(x),F(x+td)-F(x))}{t} \\
&= \left(F(x), \lim_{t\to 0} \frac{(F(x),F(x+td)-F(x))}{t}\right) \\
&= (F(x),DF(x)(d)) \quad \exists, \forall x\in X, d\in X
\end{aligned}$$

Now we have that

$$D\Phi(x)d_N(x) = (F(x),DF(x)(d)) = \left(F(x),DF(x)-DF(x)^{-1}F(x)\right) = -(F(x),Fx) = -2\Phi(x)$$

and the result follow □

**We define $d\in X$ be a descent direction of $\Phi$ at x if $D\Phi(x)d<0$**

**Remark 1.9.** *Lemma (1.1) implies that Newton direction (if it well defined) is a descent direction of potential of least squares. We suggest a replacement of form $\rho d_N(x)$, where $\rho\in(0,1)$ of the form $\beta^j$ (because we have to truncate s.t $u_{3I}, w_{3I} \geqslant 0$ if we move with step size 1 with Newton direction we can get the negative zone), where $j$ is a small natural number such that.*

$$\Phi(x+\rho d) \leqslant \Phi(x) + m\rho D\Phi(x)(d) \tag{1.25}$$

*where $m\in\left(0,\frac{1}{2}\right)$, and $\beta\in(0,1)$. That is the type of line search for our algorithm.*

*But there is trouble with Armijo line search we can't guarantee that $u_{3I}, w_{3I} \geqslant 0$ because on Armijo condition we just have (we also provide another form of Armijo line search in which we change the merit function for effort control positivity of $u_{3I}, w_{3I}, \forall 1\leqslant I\leqslant N$ which is described later.)*

$$\Phi(x^k+\beta^{j_k}d_k) \leqslant \Phi(x^k) + m\beta^{j_k}D\Phi(x^k)d^k \tag{1.26}$$

*also if we move with $t_0 = 1$ then we can get into the negative zone, and we can't control the value of $\Phi(x) = \frac{1}{2}\left\|F(x)\right\|^2$ where F(x) = [Au-b-w,w] while we move with the step size $t_0$ =1 and then scale the step size to $t = \beta t_0$. Here is the picture for illustration.*

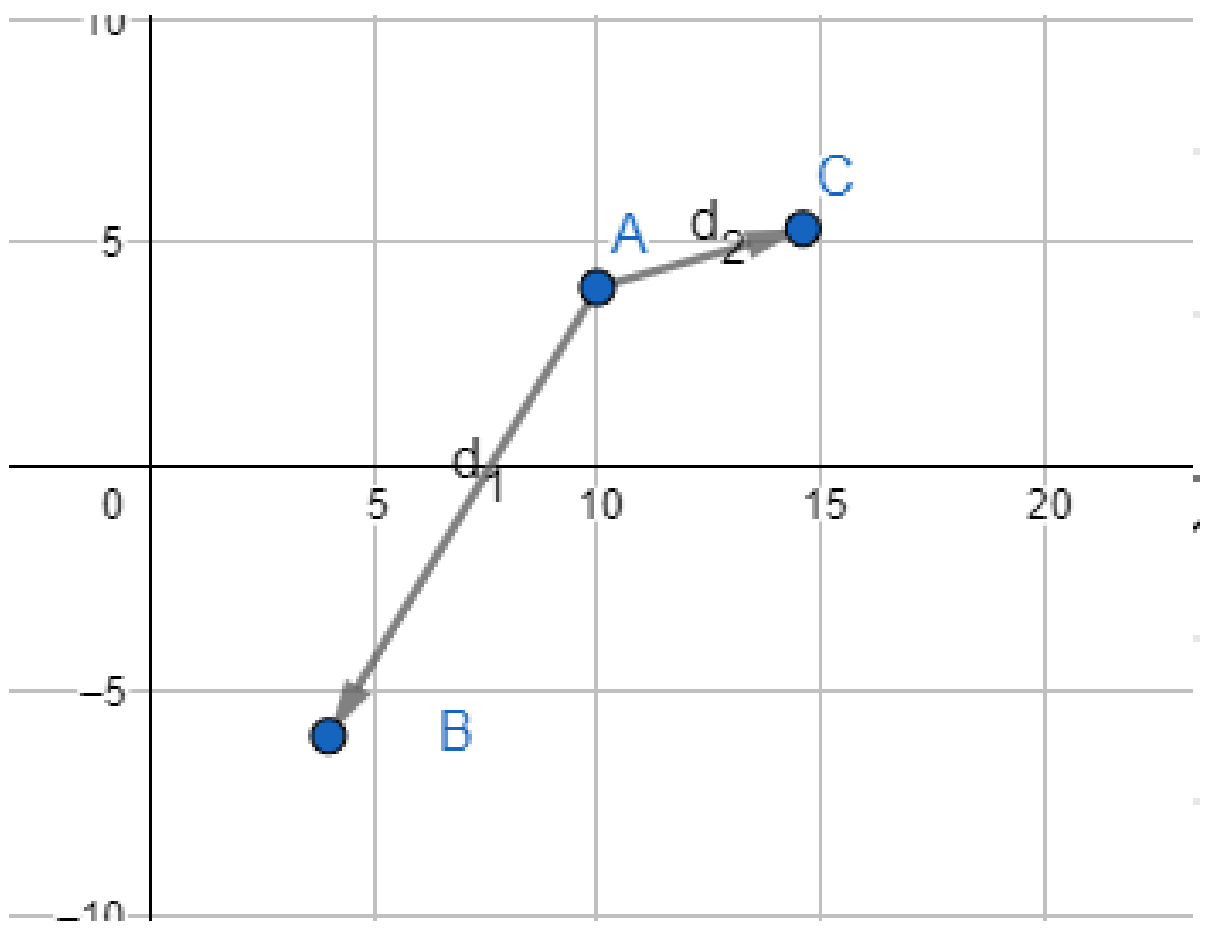

Figure 1.7: For 2 newton direction , One is leads to the negative zone we can't control the positivity of $u_{3I}, w_{3I}, \forall 1 \leqslant I \leqslant N$ also we can't control the value of $\Phi(x)$ , one is lead to the positive zone but we can't control the value of $\Phi(x)$

In practice, when we run the Armijo line search, we usually met the condition max iteration number but we don't meet the condition (even when Newton direction go to positive zone (zone directed by $d_2$)), negative zone (zone directed by $d_1$)

$$\Phi(x^k + \beta^{j_k} d_k) \leqslant \Phi(x^k) + m\beta^{j_k} D\Phi(x^k) d^k \tag{1.27}$$

This show that Armijo line search should be starting from another $t_0$ should guarantee that $u_{3I}, w_{3I} \geqslant 0$, an idea is to start with the $\alpha_{max}$ we analyze before, and we can guarantee that the non-negative value of $u_{3I}, w_{3I}$ also start with that we can have the convergence of our algorithm which confirms by numerical results as follow.

Thus , a second IPM algorithm could be

**Algorithm 3** IPM parameter with Armijo line search

1: Starting from an initial guess $\mu^{(0)} > 0, (\mathbf{u}^{(0)}/(u_{3I}^{(0)}, w_I^{(0)}) > 0)$, $m \in \left(0, \frac{1}{2}\right)$ , $\beta \in (0,1)$, k := 0
2: If $F(\mathbf{u}^{(k)}) = 0$ or $DF(\mathbf{u}^{(k)})$ is not invertible: stop; else
3: Calculate $d^k$ by solving
$$F(\mathbf{u}^k) + DF(\mathbf{u}^k) d^k = 0$$
4: Truncate $(u_{3I}^{(k+1)}, w_I^{(k+1)})$ to ensure $u_{3I}^{(k+1)} \geqslant 0$ and $w_I^{(k+1)} \geqslant 0$
5: Armijo type line search : choose $\rho^k \in (0,1)$ of the form $\beta^{j_k}$ , where $j_k$ is a small natural number s.t
$$\Phi(\mathbf{u}^{(k)} + \beta^{j_k} d^k) \leqslant \Phi(\mathbf{u}^{(k)}) + m\beta^{j_k} D\Phi(\mathbf{u}^{(k)}) d^k \tag{1.28}$$
where
$$\Phi(\mathbf{u}^{(k)}) = \frac{1}{2} \left\| F(\mathbf{u}^{(k)}) \right\|^2$$
6: Update $\mathbf{u}^{k+1} = \mathbf{u}^k + \alpha_{max} d^k$, where $\alpha_{max}$ is found by step 4
7: Update $\mu^{(k+1)} = min \left\{0.8\mu^{(k)}, [\mu^{(k)}]^2\right\}$

### 1.4.3 Numerical Result of Interior-Point method with parameter

Figure below show that our methods is always convergence with any initial point.

**Remark 1.10.** *The word "always" meaning that we run our algorithm with 10 difference initial point and it converge. The word "always" is also throughout this thesis in this meaning*

*Actually, in practice we don't need it, we just need 1 good initial point. When initial point is fixed, we change N ( the size of the problem ) and we claim that we can solve with N is as large as possible . In other words, we claim that we can solve a problems as large as possible.*

The way i calculating the error as below, I create an error vector

$$error = \begin{pmatrix} (Au-b)_1 \\ (Au-b)_2 \\ min\{u_3,(Au-b)_3\} \\ (Au-b)_4 \\ (Au-b)_5 \\ min\{u_6,(Au-b)_6\} \\ ... \\ ... \\ ... \\ (Au-b)_{3N-2} \\ (Au-b)_{3N-1} \\ min\{u_{3N},(Au-b)_{3N}\} \end{pmatrix} = \begin{pmatrix} w_1 \\ w_2 \\ min\{u_3,w_3\} \\ w_4 \\ w_5 \\ min\{u_6,w_6\} \\ ... \\ ... \\ ... \\ w_{3N-2} \\ w_{3N-1} \\ min\{u_{3N},w_{3N}\} \end{pmatrix} \tag{1.29}$$

Then i take $\|error\|_\infty$. This figure below to show that the IPM methods is always converge I use the safety parameter ($\theta$) to be 0.9

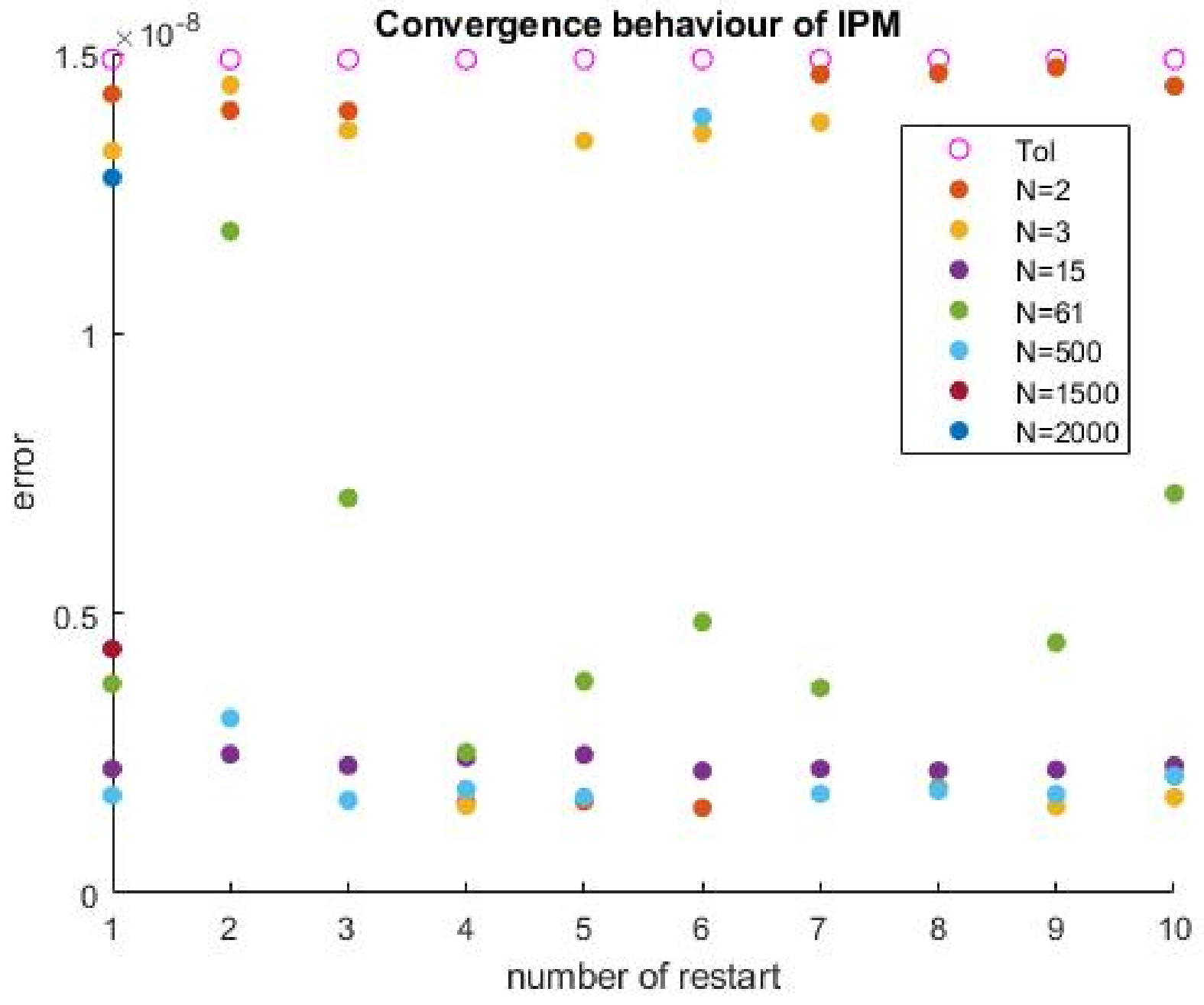

Figure 1.8: The convergence of IPM for various option N

Secondly, we propose a comparison with difference safety parameter $\theta$

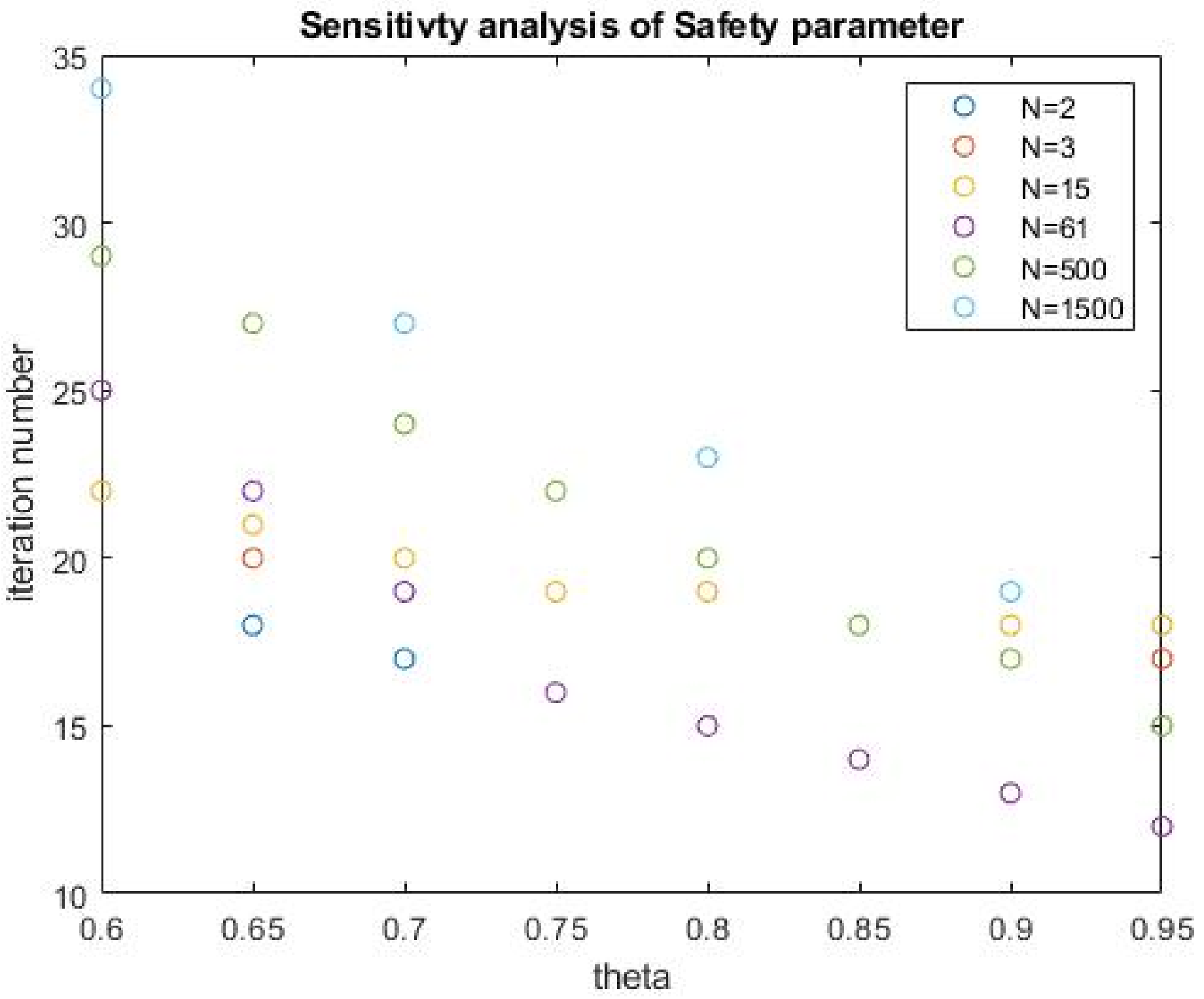

Figure 1.9: Iteration of IPM for various option N with difference safety parameter

Now based on the result we can see that the safety parameter near 1 takes less iteration than the far from 1. I think it's good if we take a safety parameter near 1.

### 1.4.4 Numerical Result of Interior Point method with parameter and Armijo line search.

In this section we used Armijo line search with $t_0 = \alpha_{max}$ to guarantee that $u_{3I}, w_{3I} \geqslant 0, \forall 1 \leqslant I \leqslant N$. Safety parameter is set to be 0.9

Firstly we show that Interior-Point method with parameter and Armijo line search convergence rate is 100%

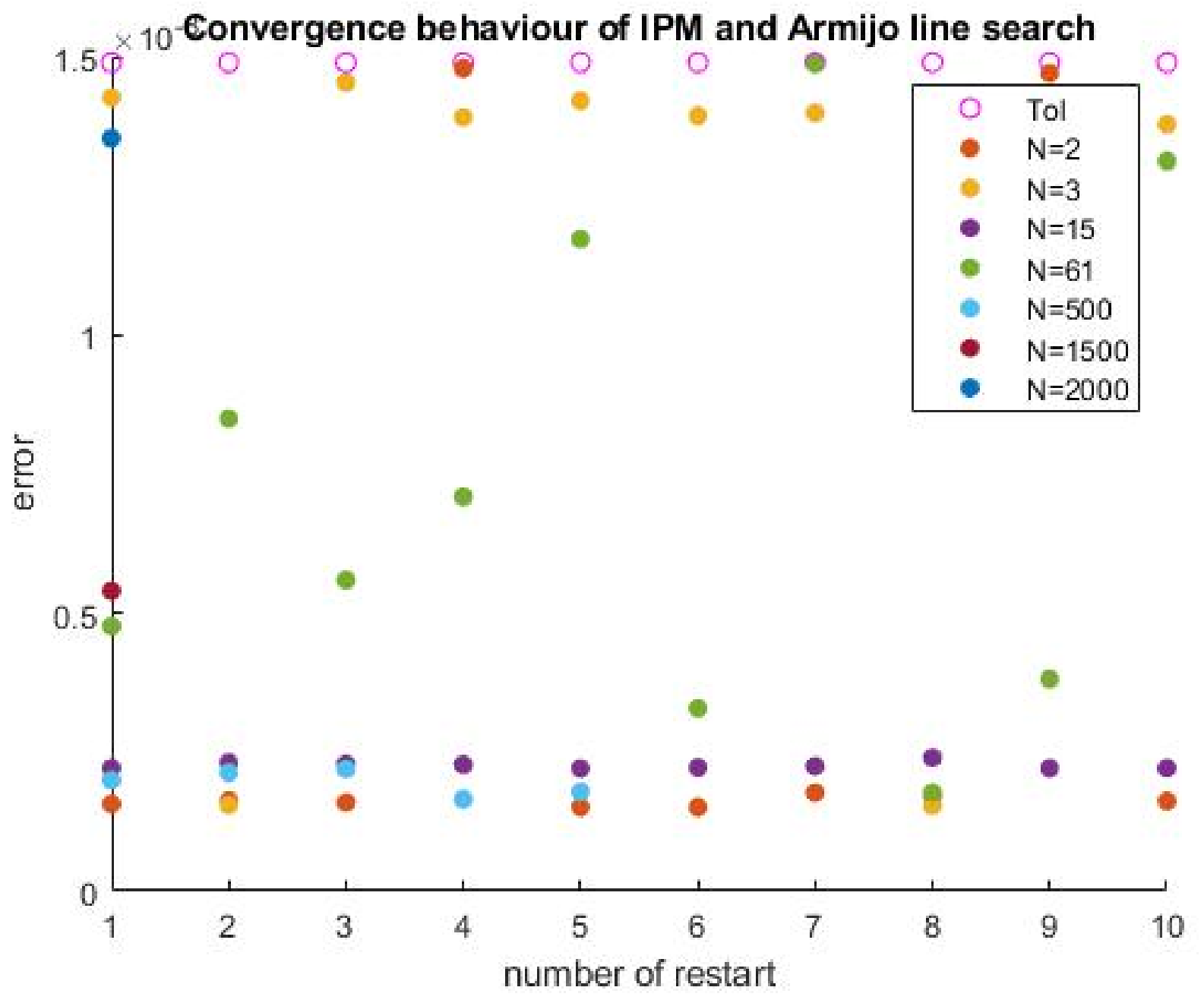

Figure 1.10: The convergence of IPM with Armijo line search for various option N

Secondly, we propose a comparison of iteration with various safety parameter $\theta$

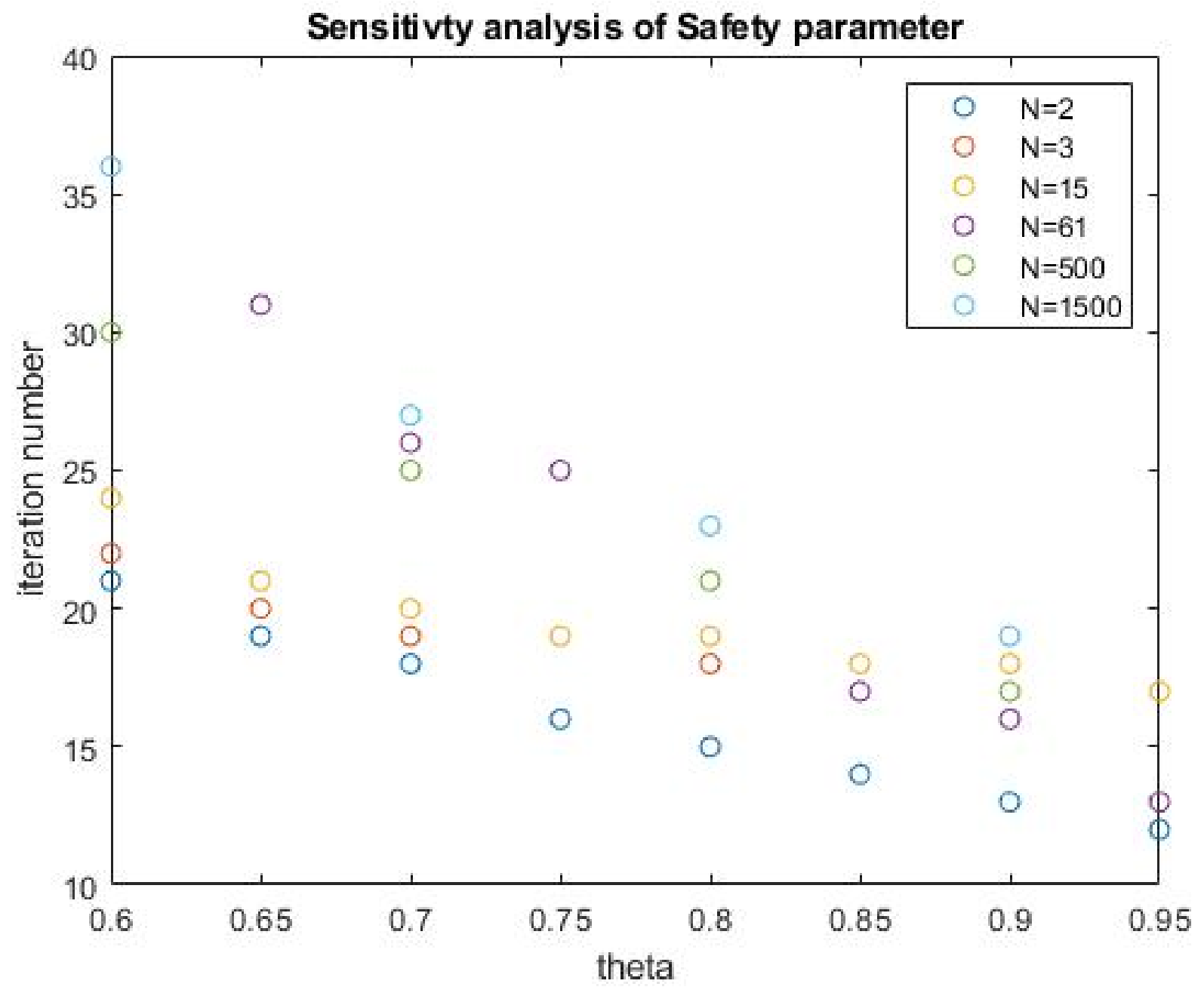

Figure 1.11: Iteration of IPM Armijo for various option N with difference safety parameter

The same with the previous, It is good to take $\theta$ to be near 1 to save the iterations.

**Remark 1.11.** *One observation that it's not difference too much between **"IPM "** and **"IPM with Armijo line search"**. Define $z = |\alpha_{max} - \rho|$, with $\rho$ is the step size find in Armijo function with $t_0 = \alpha_{max}$. Run our function with various N we see that .*

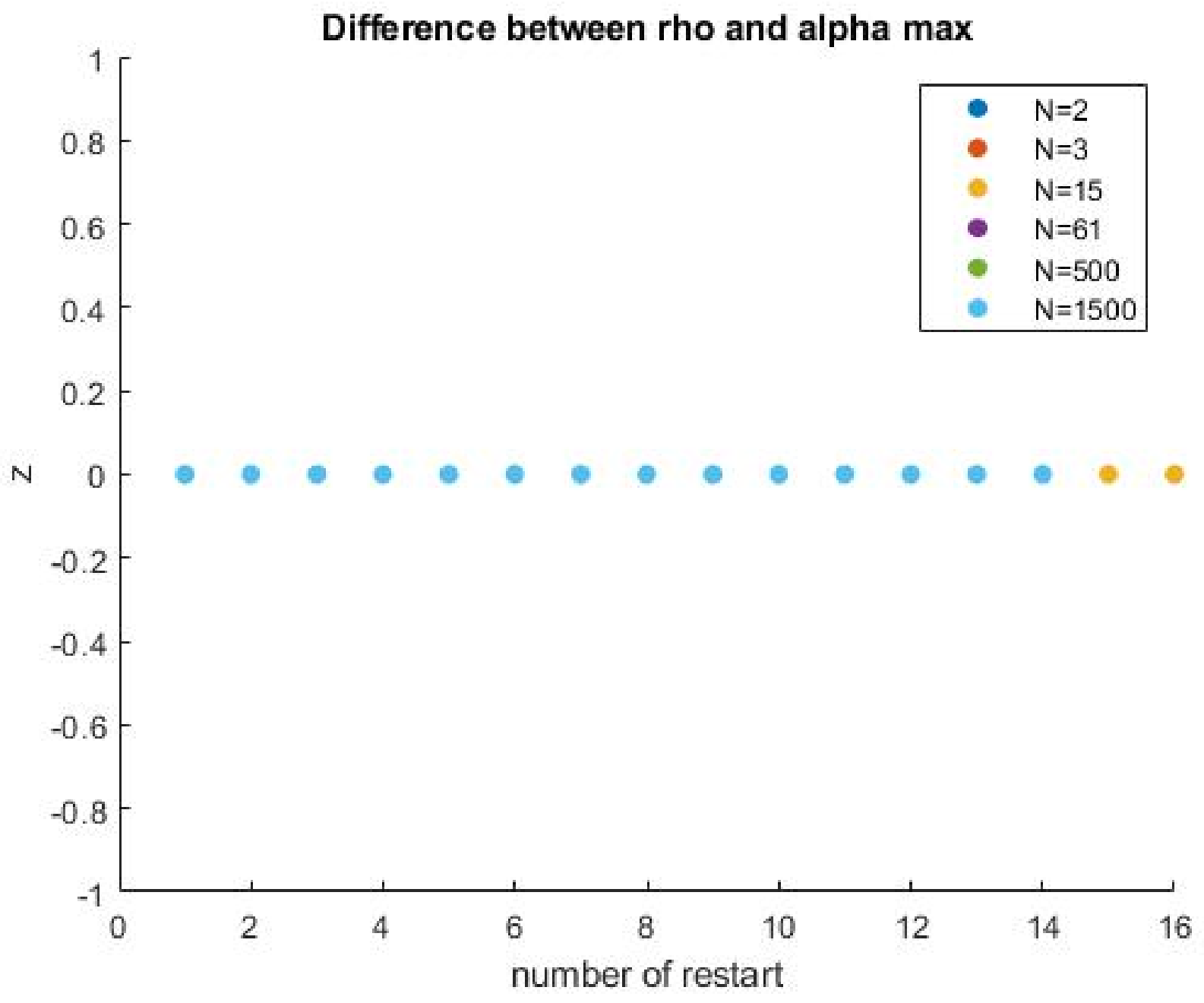

Figure 1.12: Difference between $\rho$ & $\alpha_{max}$ for various option N

*That means with a large number of restart. Armijo start with $t_0 = \alpha_{max}$ and return $t_0$. It's no meaning that Armijo doesn't always work. Sometimes it works sometimes not we don't know it.*

**Remark 1.12.** *Ones natural idea is to decrease the value of $\mu$ faster than $\mu = min\left\{0.8\mu, \mu^2\right\}$ **in order to save computational time.** However, we have to pay attention to the convergence. We considered 3 strategies*

- $\mu_1 = min\left\{0.8\mu, \mu^2\right\}$
- $\mu_2 = min\left\{0.1\mu, \mu^2\right\}$
- $\mu_3 = min\left\{0.01\mu, \mu^3\right\}$

| **Cases** | Lowest nb_iter | Highest nb_iter |
|---|---|---|
| $N = 2$ | (12,10,15) | (13,20,19) |
| $N = 3$ | (17,13,9) | (17,16,12) |
| $N = 15$ | (17,12,12) | (17,15,26) |
| $N = 61$ | (12,20,22) | (30,35,37) |
| $N = 500$ | (12,13,12) | (13,13,13) |
| $N = 1500$ | (15,11,10) | (15,29,12) |

*we can see that some times it is not effective when we decrease the value of $\mu$ in a fast way in case N = 2, N= 61, but sometimes we need to decrease the value of $\mu$ in a fast way in case N = 3, N = 15, N = 1500. It does not vary too much in the case N =500. That means we have to choose the value of $\mu$ in some "smart" way but not in general, even with the "smart" ways we choose, With $\mu_1 = min\left\{0.8\mu, \mu^2\right\}$ sometimes it's better sometimes not. It also serves a motivation to the latter ones which mean we don't need to think about the "smart" way we have to choose for parameter $\mu$ for each position. It becomes our variables (an unknown of the system).*

## 1.5 New algorithm for complementary condition problems

A more sophisticated algorithm corresponds to what is called **"Non-parametric IPM"** in Son's work (come from Ph.D. thesis of Son from **IFPEN**). The idea is to regard $\mu$ as an unknown and to supply the system with an additional equation enforcing $\mu = 0$ in a judicious way.

$$\begin{cases} (Au - b - w)_{3I-2} = \mathbf{0} \\ (Au - b - w)_{3I-1} = \mathbf{0} \\ (Au - b - w)_{3I} = \mathbf{0} \\ w_{3I-2} = \mathbf{0} \\ w_{3I-1} = \mathbf{0} \\ u_{3I} \bullet w_I = \mu \mathbf{e} \end{cases} \tag{1.30}$$

for $1 \leqslant I \leqslant N$, and a last equation:

$$\frac{1}{2}\sum_{I=1}^{N} min\left\{u_{3I};0\right\}^2 + \frac{1}{2}\sum_{I=1}^{N} min\left\{w_{3I};0\right\}^2 + \mu^2 + \epsilon\mu = 0 \tag{1.31}$$

for a small value of $\epsilon > 0$ (unchanged during the iterations). System (5.4) - (1.31) consists of 6N+1 unknowns and 6N+1 equation. An ordinary Newton method is then applied to solve it.

**Ex. 2** — 1. Check with Professor Mounir Haddou about the additional equation (1.31)

2. Write down the Jacobian matrix and the Newton iteration for (5.4) - (1.31)

**Answer (Ex. 2)** — We need an equation to ensure the non-negative of slack variables. We consider function $f(u) = \frac{1}{2}(min\left\{u;0\right\})^2$ . It is easy to see that $f(u) = 0$ when $u \geqslant 0$. Then we can add following equation

$$\frac{1}{2}\sum_{I=1}^{N} min\left\{u_{3I};0\right\}^2 + \frac{1}{2}\sum_{I=1}^{N} min\left\{w_{3I};0\right\}^2 + \mu^2 = 0 \tag{1.32}$$

from this, we can imply that $\mu = 0$ and all $u_{3I}, w_{3I}$ are non-negative.

$$F(\chi, \mu) = \begin{bmatrix} \Im(\chi) - \mu e \\ \frac{1}{2}\sum_{I=1}^{N} min\left\{u_{3I};0\right\}^2 + \frac{1}{2}\sum_{I=1}^{N} min\left\{w_{3I};0\right\}^2 + \mu^2 \end{bmatrix} \in \mathbb{R}^{2m+1,1} \tag{1.33}$$

Let $D\Im(\chi)$ be the Jacobian matrix of $\Im$, since $w_{3I} \geqslant 0$ and $u_{3I} \geqslant 0, \forall \quad 1 \leqslant I \leqslant N$, then we have Jacobian matrix of F

$$DF(\chi, \mu) = \begin{bmatrix} D\Im(\chi) & -e \\ 0 & 2\mu \end{bmatrix} \in \mathbb{R}^{2m+1,2m+1} \tag{1.34}$$

where N is the number of elements, m = 3N, e = $(0, ..., 0, 0, 0, -1, 0, 0, -1, ..., 0, 0, -1)^T \in \mathbb{R}^{2m,1}$ with m component above is zeros, component m+3 = -1 ,m+6 = -1,...,m+m = -1. $D\Im(\chi) \in \mathbb{R}^{2m,2m}$ is as in exercise before. $\mathbb{O} = (0, 0, 0, ..., 0) \in \mathbb{R}^{1,2m}$. If $F(\chi, \mu) = 0$ then $\mu = 0$. It leads to matrix $DF(\chi, \mu)$ is singular since $\det DF(\chi, \mu) = 0$. Here is the numerical confirm our theory results when we run the code sometimes Matlab inform that

```
mu =

   5.4707e-11


mu =

   2.7353e-11

Warning: Matrix is close to singular or
badly scaled. Results may be inaccurate.
RCOND =  2.190305e-16.
> In IPM_nonpara (line 54)
  In main (line 36)
```

Figure 1.13: Matrix is close to singular

**Remark 1.13.** *In theory, we don't compute the $\alpha_{max}$ we just the Armijo line search so we don't guarantee $\mu > 0$ as always. So maybe sometimes $\mu \leqslant 0$ can happen. So it leads to we have to put $\epsilon$ small. the reason why we need $\epsilon > 0$ is small but not too small because if it too small The Jacobi matrix near to bad cases ( singular). But it is not too large if we look at the graph of the function $\mu^2 + \epsilon\mu$.*

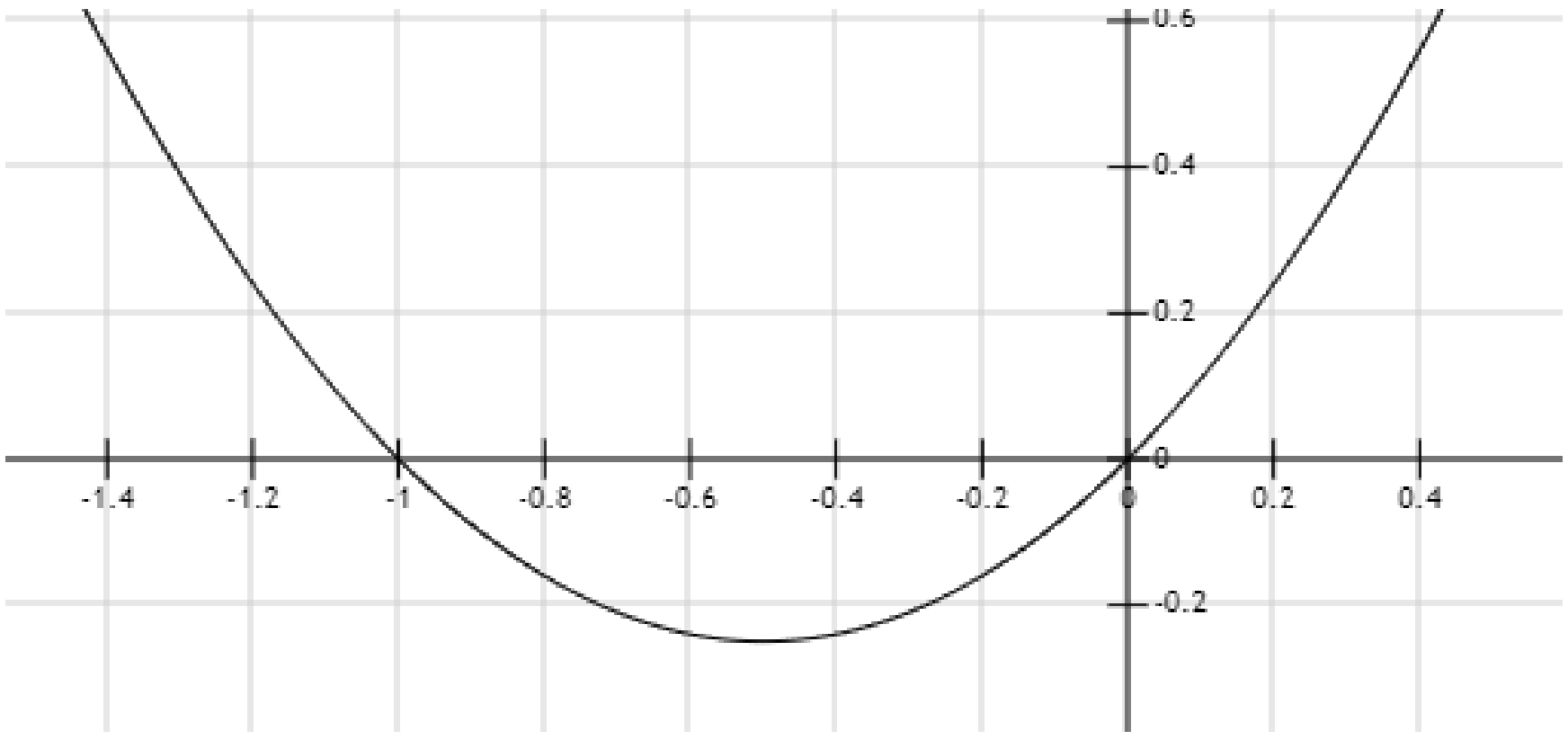

Figure 1.14: Graph of $x^2 + x$

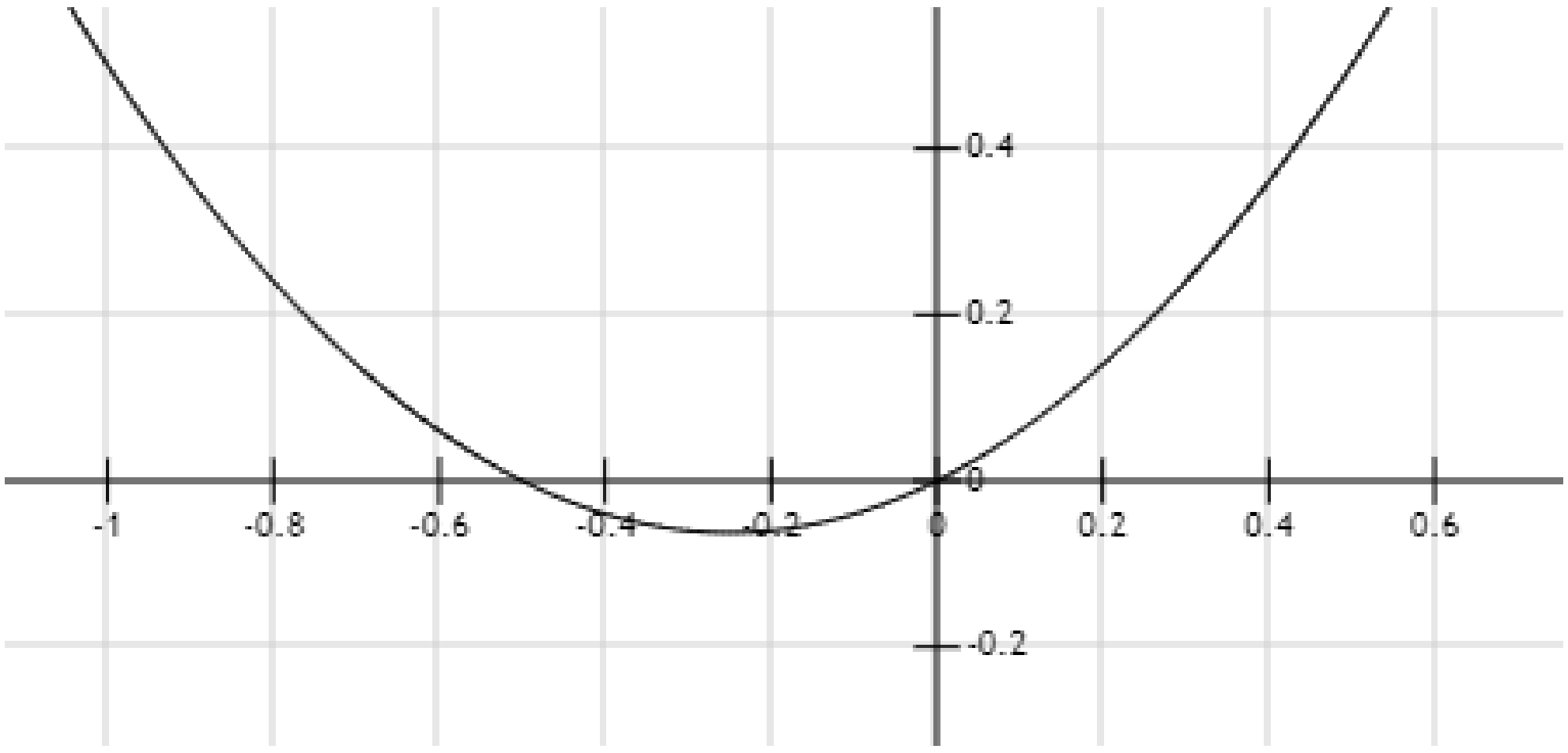

Figure 1.15: Graph of $x^2 + \frac{1}{2}x$

*we can see that the negative zone is narrowing down when we decrease the value of $\epsilon$ and that makes sense because if we are in negative zone then our function*

$$\frac{1}{2}\sum_{I=1}^{N} min\left\{u_{3I};0\right\}^2 + \frac{1}{2}\sum_{I=1}^{N} min\left\{w_{3I};0\right\}^2 = \tau > 0 \tag{1.35}$$

*$\tau$ sometimes is small when we do numerical test but we can accept that. So if we take $\epsilon > 0$ is too large $\mu$ can be large enough make $\tau$ become larger in comparison with the previous ones and we can't accept that (because it's too large). It leads to we don't have the positivity of $u_{3I}, w_{3I}$.*

**Remark 1.14.** *In practice, we choose the ways to convergence is choose an initial guess $[u_{3I}^{(0)}, w_{3I}^{(0)}] > 0\ \forall i \in \{1, ..., N\}$, we have the thirds equation is*

$$u_{3I}w_{3I} = \mu, \forall I \in \{1, ..., N\} \tag{1.36}$$

*so that we have $(u, w) = N\mu$ so we choose*

$$\mu^{(0)} = \frac{1}{N}(u, w) = \frac{1}{N}\sum_{I=1}^{N} u_{3I}w_{3I}$$

*and we have the positivity of $\mu^{(0)}$. We also calculate $\alpha_{max}$ to guarantee $u_{3I}, w_{3I} \geqslant 0$.*
Another form of this is

$$\frac{1}{2}\sum_{I=1}^{N} min\left\{u_{3I};0\right\}^2 + \frac{1}{2}\sum_{I=1}^{N} min\left\{w_{3I};0\right\}^2 + \epsilon\mu = 0 \tag{1.37}$$

but in theory when we look at the graph of function $\epsilon\mu$

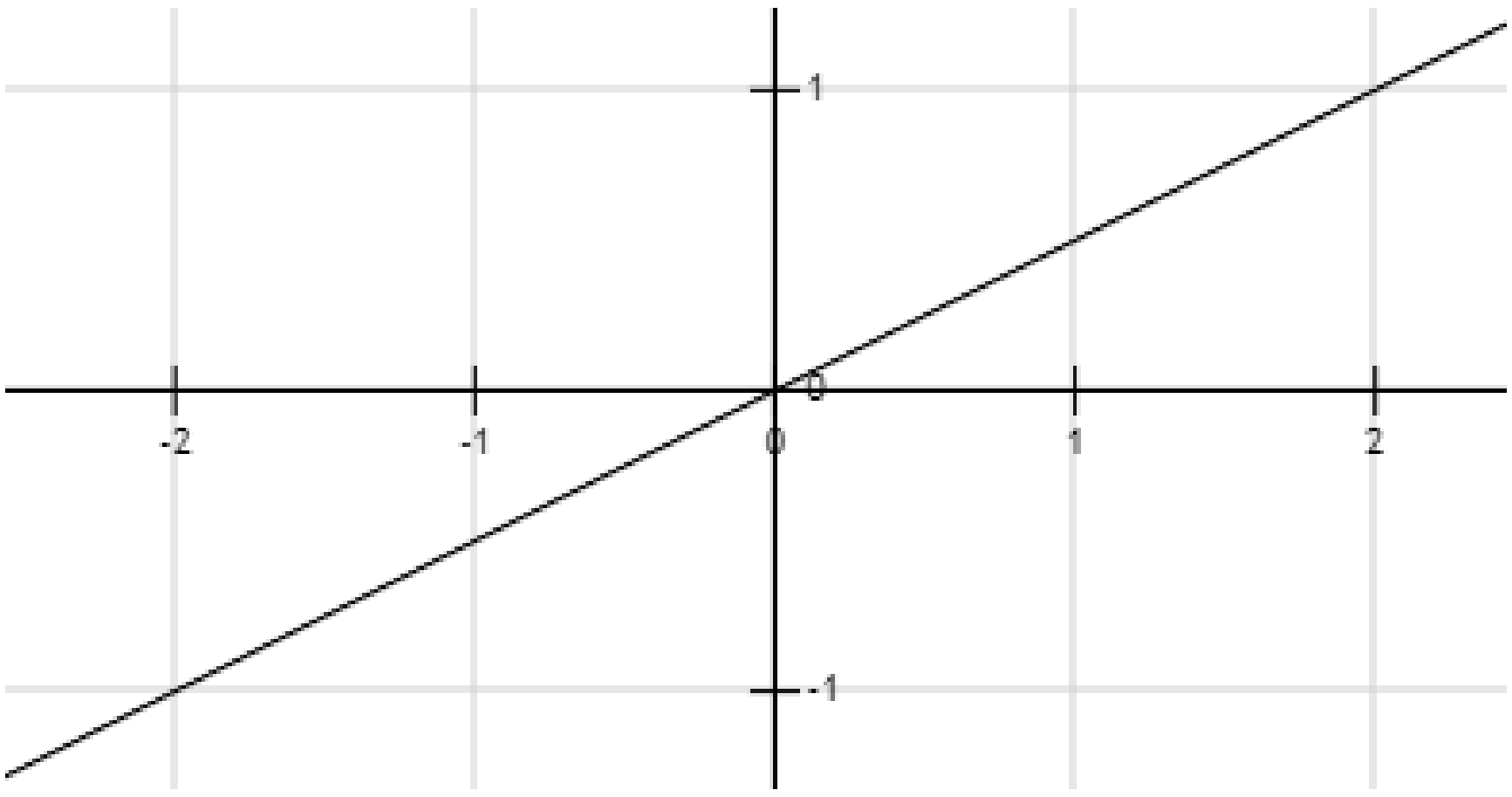

Figure 1.16: Graph of $\frac{1}{2}x$

we can see that the negative zone is too large and that's why we don't use it(reason is the same as above). One more reason that's we use linear approximation to linear scheme with respect to u. That is $d\mu^{(k)} = -\frac{\epsilon\mu^{(k)}}{\epsilon} = -\mu^{(k)}$ then the Newton direction $\mu^{(k)} = \mu^{(k)} + d\mu^{(k)} = \mu^{(k)} - \mu^{(k)} = 0$ and the iteration stop at 1 step(that means after that we back to the origin problems which things we can't solve it as the beginning so we have to make the perturbed system ). Or it become constant because $\alpha_{max} = 0$. Let's say we move with step size = 1 then $\exists i \in \{3I, 6I, ...\}$ s.t $x_i + \alpha_{max} d_i = 0$ then $x_i$ at the next step is zeros, so that $\alpha_{max}$ after that is 0 so $\mu = \mu + \alpha_{max} d_n = \mu$ so we don't move anymore and our algorithm stop. Numerical result is follow

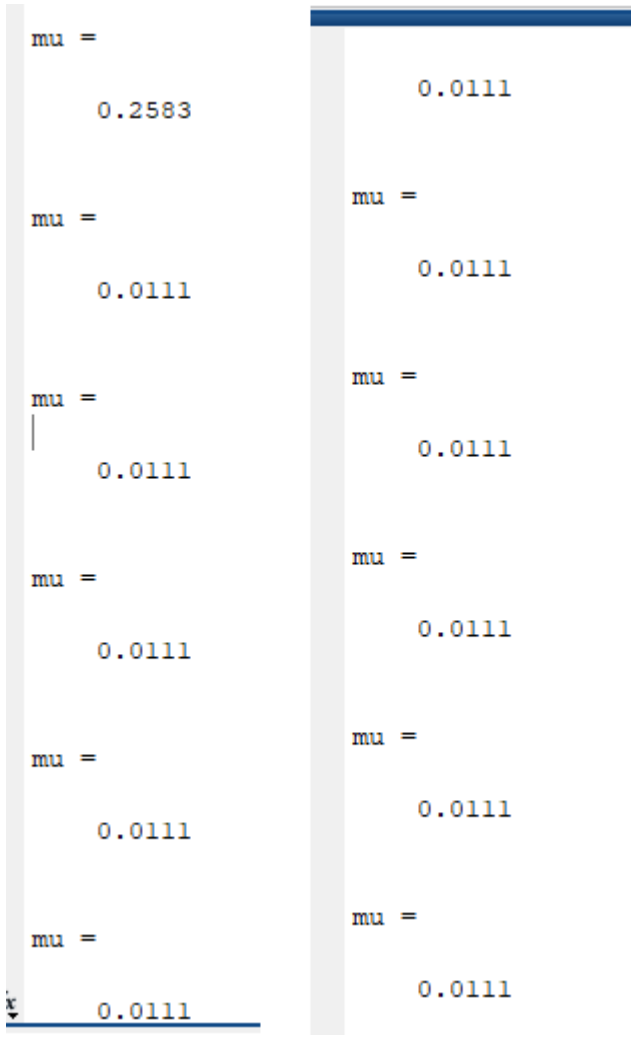

Figure 1.17: how $\mu$ change during the iteration .

by visualization we can say that $\mu$ just change 1 time at the beginning then it doesn't move anymore, or we can calculate the error after step 1 (the error is based on how does $\mu$ change at each step) and we get the results is 0 (means that after 1st step, $\mu$ don't change anymore)
To avoid this problem, we will add a small enough positive parameter $\epsilon$ to equation (1.32) and get

$$\frac{1}{2}\sum_{I=1}^{N} min\left\{u_{3I};0\right\}^2 + \frac{1}{2}\sum_{I=1}^{N} min\left\{w_{3I};0\right\}^2 + \mu^2 + \epsilon\mu = 0 \tag{1.38}$$

set X = $(\chi, \mu) \in \mathbb{R}^{2m+1,1}$ then we get

$$DF(X) = DF(\chi, \mu) = \begin{bmatrix} D\Im(\chi) & -e \\ \mathbb{O} & 2\mu + \epsilon \end{bmatrix} \in \mathbb{R}^{2m+1,2m+1} \tag{1.39}$$

where e $= (0, ..., 0, 0, 0, -1, 0, 0, -1, ..., 0, 0, -1)^T \in \mathbb{R}^{2m}$ with m component above is zeros, component m+3 = -1 ,m+6 = -1,...,m+m = -1. $D\Im(\chi) \in \mathbb{R}^{2m}$ is as in exercise before. $\mathbb{O} = (0, 0, 0, ..., 0) \in \mathbb{R}^{2m}$ Now at the solution $(\chi, \mu)$, we have $detDF(\chi, \mu) = \epsilon detD\Im(\chi)$. It means that Jacobian matrix DF is singular iff Jacobian matrix $D\Im$ is singular **(also meaning that we shouldn't take $\epsilon$ to be too small - it can make some numeric error make our matrix here is close to singular )**

**Remark 1.15.** *Here I make the sensitivity analysis of the problems when $\epsilon \downarrow 0$.*

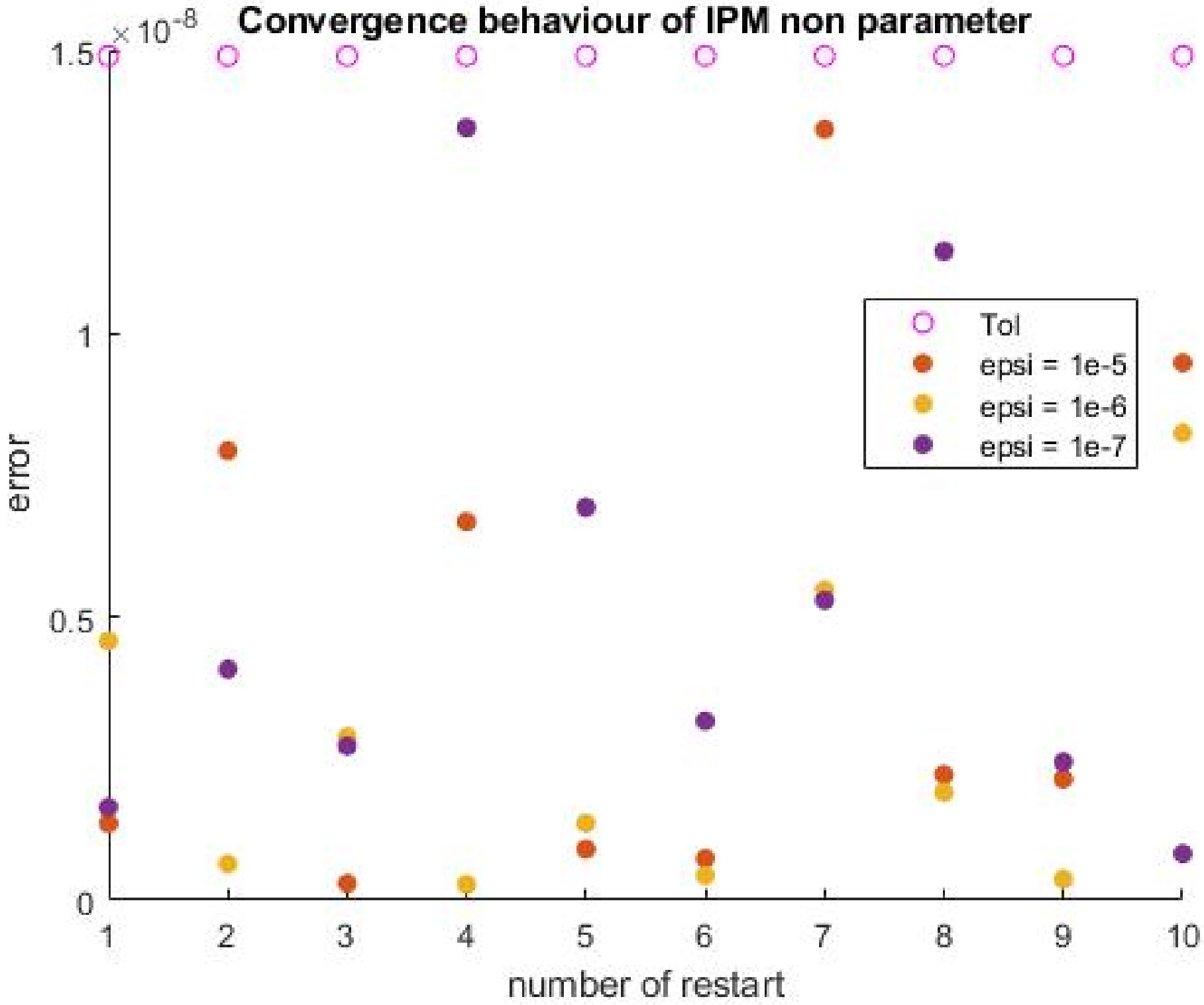

Figure 1.18: The convergence of non-parametric IPM when N = 3

*Figure above show that when $\epsilon >> 1e-7$, our algorithm can run but when $\epsilon = 1e-7$. Matlab inform that*

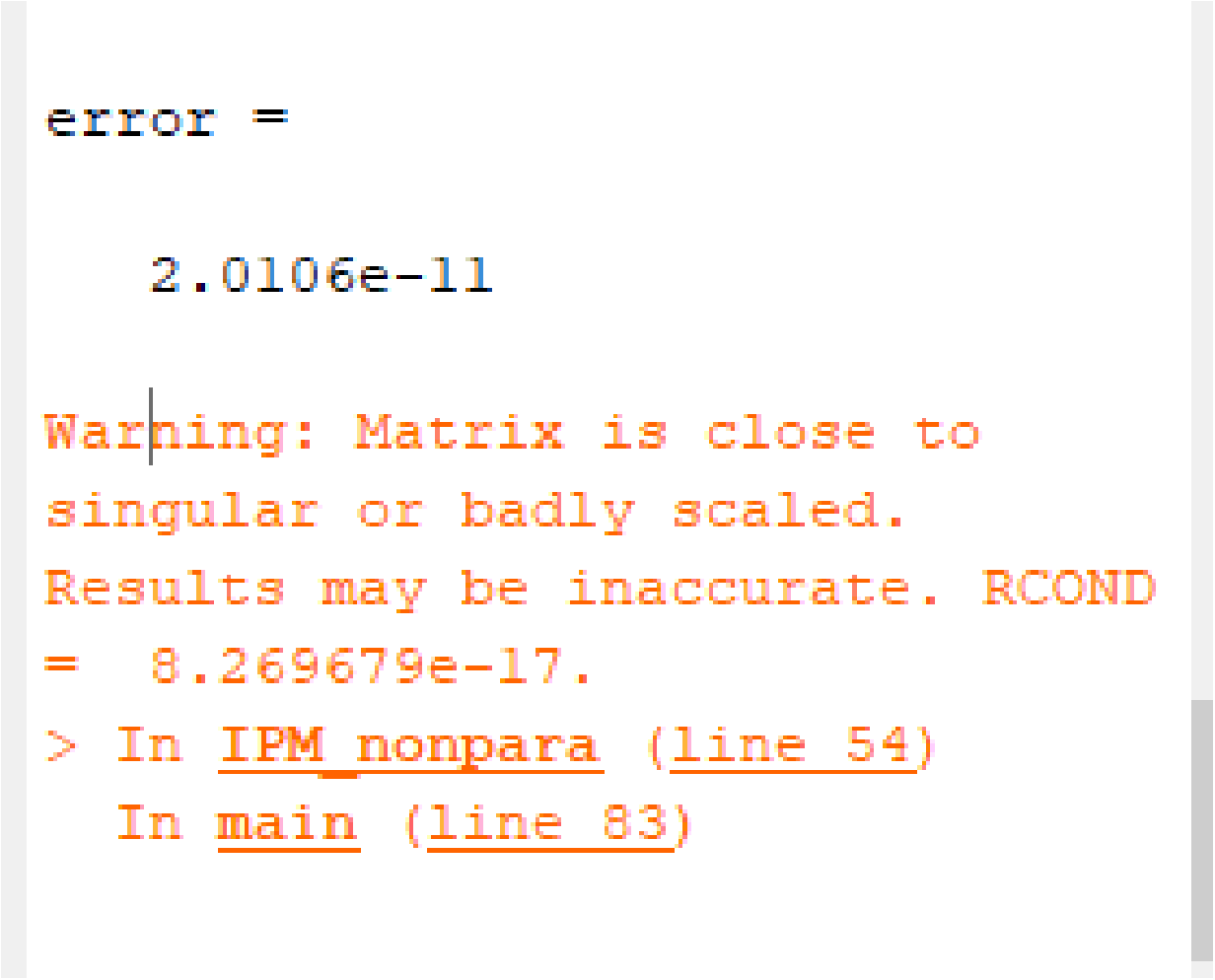

```
error =

   2.0106e-11

Warning: Matrix is close to
singular or badly scaled.
Results may be inaccurate. RCOND
=  8.269679e-17.
> In IPM_nonpara (line 54)
  In main (line 83)
```

Figure 1.19: Matrix is close to singular

*So we don't need to take $\epsilon$ become too small. It has no further profits for us.*

**The non-parametric IPM algorithm is as follow**

---

**Algorithm 4** "Non-parametric IPM"

---

1: Starting from an initial guess $\mu^{(0)} > 0, (\mathbf{u}^{(0)}/(u_{3I}^{(0)}, w_I^{(0)}) > 0)$

2: If $F(\mathbf{u}^{(k)}) = 0$ or $DF(\mathbf{u}^{(k)})$ is not invertible: stop; else

3: Calculate $d^k$ by solving

$$F(\mathbf{u}^k) + DF(\mathbf{u}^k)d^k = 0$$

If all $u_{3I}$ and slack variables are non-negative, then we have

$$d^k = \begin{bmatrix} d\mathfrak{X}^k \\ d\mu^k \end{bmatrix} = -\begin{bmatrix} D\mathfrak{F}(\mathfrak{X}) & -e \\ 0 & 2\mu+\epsilon \end{bmatrix}^{-1} \begin{bmatrix} F\left(\mathfrak{X}^k\right) - \mu^k e \\ \left(\mu^k\right)^2 + \epsilon\mu^k \end{bmatrix}$$

and

$$d\mu^k = -\frac{\left(\mu^k\right)^2 + \epsilon\mu^k}{2\mu^k + \epsilon}$$

4: Truncate $(u_{3I}^{(k+1)}, w_I^{(k+1)})$ to ensure $u_{3I}^{(k+1)} \geqslant 0$ and $w_I^{(k+1)} \geqslant 0$

5: Update $\mathbf{u}^{k+1} = \mathbf{u}^k + \alpha_{max} d^k$, where $\alpha_{max}$ is found by step 4

---

When we consider "non-parametric IPM with Armijo line search", our strategy is the same with **"IPM with Armijo line search"**. We also calculate $\alpha_{max}$ to guarantee $u_{3I}, w_{3I} \geqslant 0$ and then we use the Armijo line search to compute step size of Newton direction.

**The non-parametric IPM with Armijo line search algorithm is as follow**

**Algorithm 5** "Non-parametric IPM with Armijo line search"

1: Starting from an initial guess $\mu^{(0)} > 0$, $(\mathbf{u}^{(0)}/(u_{3I}^{(0)}, w_I^{(0)}) > 0)$
2: If $F(\mathbf{u}^{(k)}) = 0$ or $DF(\mathbf{u}^{(k)})$ is not invertible: stop; else
3: Calculate $d^k$ by solving

$$F(\mathbf{u}^k) + DF(\mathbf{u}^k)d^k = 0$$

If all $u_{3I}$ and $w_{3I}$ are non-negative, then we have

$$d^k = \begin{bmatrix} d\mathfrak{X}^k \\ d\mu^k \end{bmatrix} = -\begin{bmatrix} D\mathfrak{F}(\mathfrak{X}) & -e \\ 0 & 2\mu + \epsilon \end{bmatrix}^{-1} \begin{bmatrix} F\left(\mathfrak{X}^k\right) - \mu^k e \\ \left(\mu^k\right)^2 + \epsilon\mu^k \end{bmatrix}$$

and

$$d\mu^k = -\frac{\left(\mu^k\right)^2 + \epsilon\mu^k}{2\mu^k + \epsilon}$$

4: Truncate $(u_{3I}^{(k+1)}, w_I^{(k+1)})$ to ensure $u_{3I}^{(k+1)} \geqslant 0$ and $w_I^{(k+1)} \geqslant 0$
5: Armijo type line search:
Choose $\rho^k \in (0, 1)$ of the form $\beta^{j_k}$, where $j_k$ is a small natural number such that

$$\Phi\left(\mathbf{u}^{(k)} + \beta^{j_k k} d^k\right) \leq \Phi\left(\mathbf{u}^{(k)}\right) + m\beta^{j_k} D\Phi\left(\mathbf{u}^{(k)}\right) d^k$$

where

$$\Phi(\mathbf{u}^{(k)}) = \frac{1}{2}\left\|F(\mathbf{u}^{(k)})\right\|^2$$

6: Update $\mathbf{u}^{k+1} = \mathbf{u}^k + \alpha_{max} d^k$, where $\alpha_{max}$ is found by step 5

**Remark 1.16.** *We need a truncation because we want during the iteration we always have* $u_{3I} \geqslant 0$ *and* $w_{3I} \geqslant 0$

**Remark 1.17.** *With this step size on* $\mu$ *(*$d\mu$*), when we start with initial* $\mu^{(0)} > 0$*, we still remain* $\mu^{(k)} > 0, \forall k$ *during the iteration. Indeed, suppose* $\mu^{(k)} > 0$ *then*

$$\begin{aligned} \mu^k + d^k &= \mu^k - \frac{\left(\mu^k\right)^2 + \epsilon\mu^k}{2\mu^k + \epsilon} \\ &= \frac{\mu^k(2\mu^k + \epsilon) - \left(\mu^k\right)^2 - \epsilon\mu^k}{(2\mu^k + \epsilon)} \\ &= \frac{2\left(\mu^k\right)^2 + \mu^k\epsilon - \left(\mu^k\right)^2 - \mu^k\epsilon}{(2\mu^k + \epsilon)} \\ &= \frac{\left(\mu^k\right)^2}{2\mu^k + \epsilon} < \frac{\left(\mu^k\right)^2}{\epsilon} \end{aligned}$$

*since* $d_k < 0$ *and* $1 \geqslant \alpha_{max} > 0$ *so* $d_k \leqslant \alpha_{max} d_k < 0$ *we get*

$$\mu^k > \mu^{k+1} = \mu^k + \alpha_{max} d_k \geqslant \mu^k + d_k > 0$$

*That is we drive* $\mu \downarrow 0$*, In theory,* $\mu$ *can be negative but in practice,* $\mu$ *is always positive,* $\mu$ *is small enough and we can accept this error on numerical. With this "converge" ways, we can* *avoid trouble to choose met $\mu < 0$. Also from this if we are near to the solution then we get the quadratic convergence.*

***Here I check in practice $\epsilon$ shound't be too large*** *, when $\epsilon$ is large, then we met the trouble with the convergence of our algorithm as in the case N =61, we met the maximum iteration condition (It need more iteration to get the convergence). So it's not good when $\epsilon$ too large in theory also in practice. The figure below shows what I just wrote.*

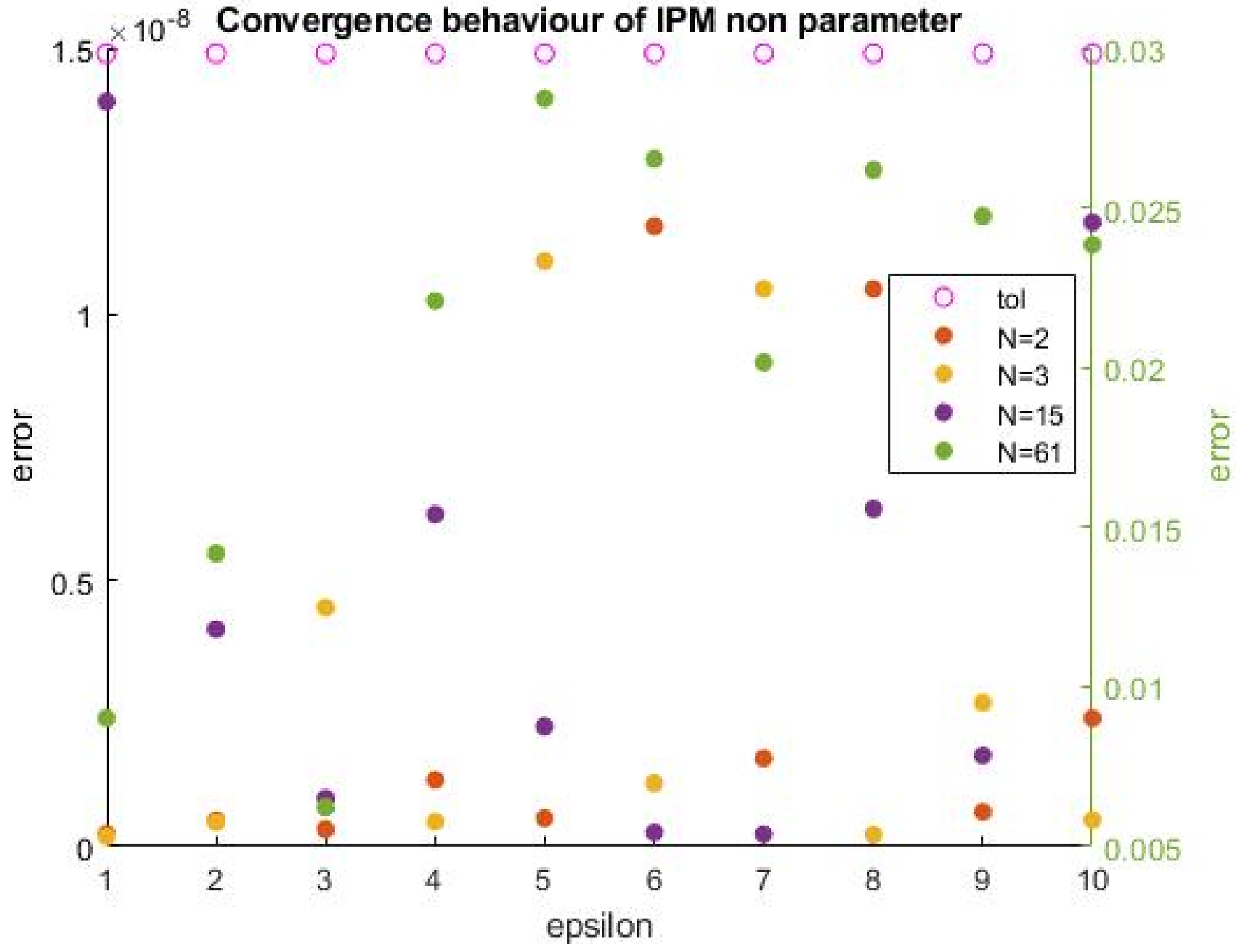

Figure 1.20: The convergence of non-parametric IPM for various option N

**Remark 1.18.** *The green dots on the picture* (1.20) *(N=61), follow the axis of error on the right, the others colored dots follow the axis of error on the left.*

### 1.5.1 Numerical Result of non-parametric Interior Point method

Before going to convergence, we make the sensitivity analysis of $\epsilon \downarrow 0$ the have a good "$\epsilon$" for non-parametric IPM.

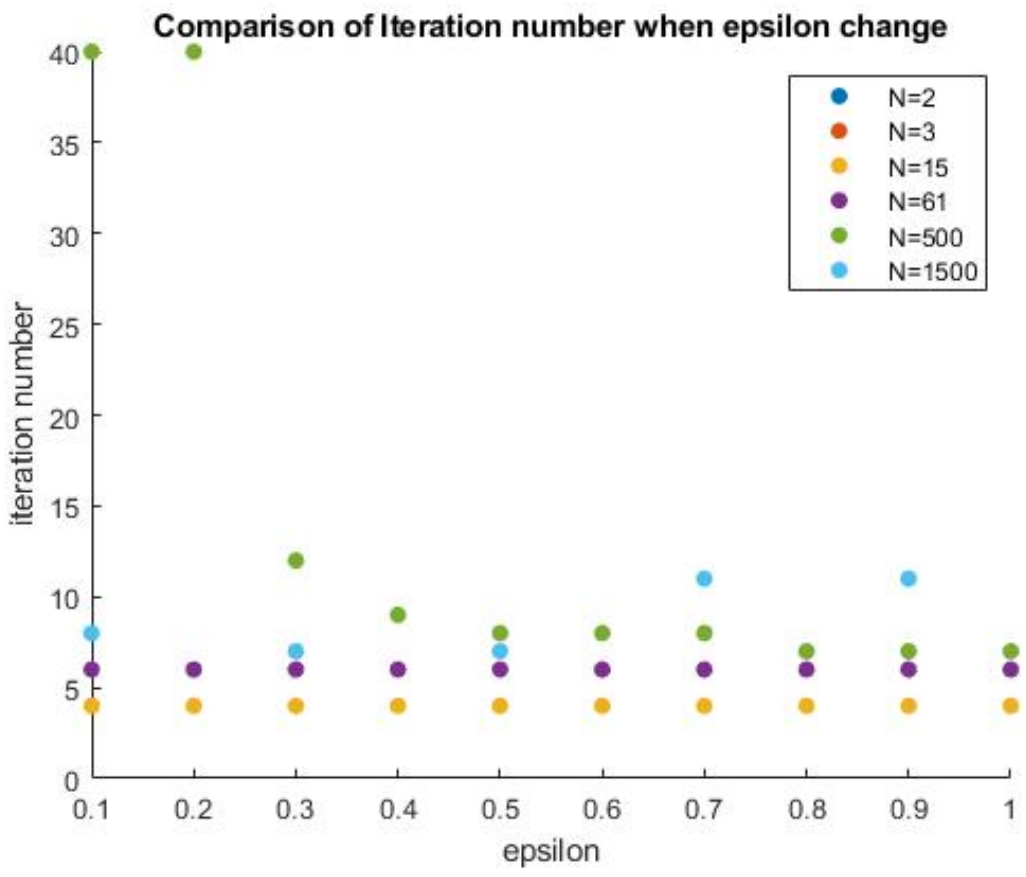

(a) the first choice of the initial point

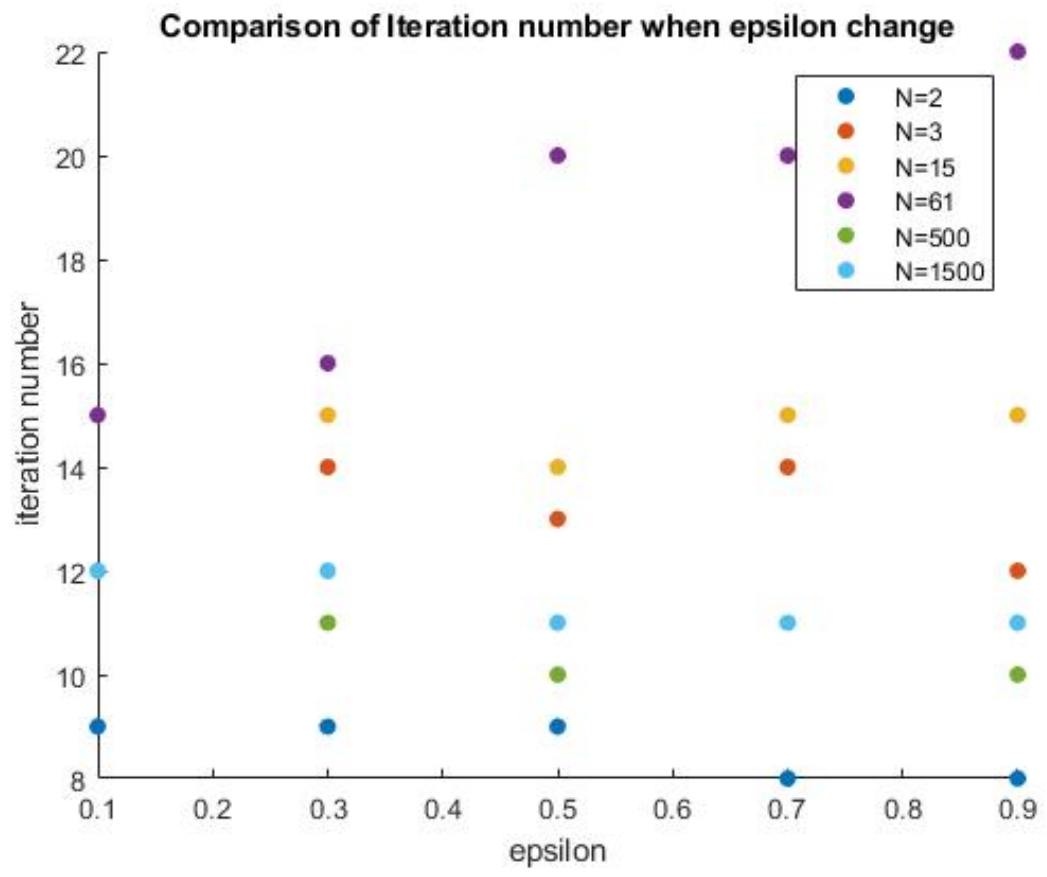

(b) the second choice of the initial point

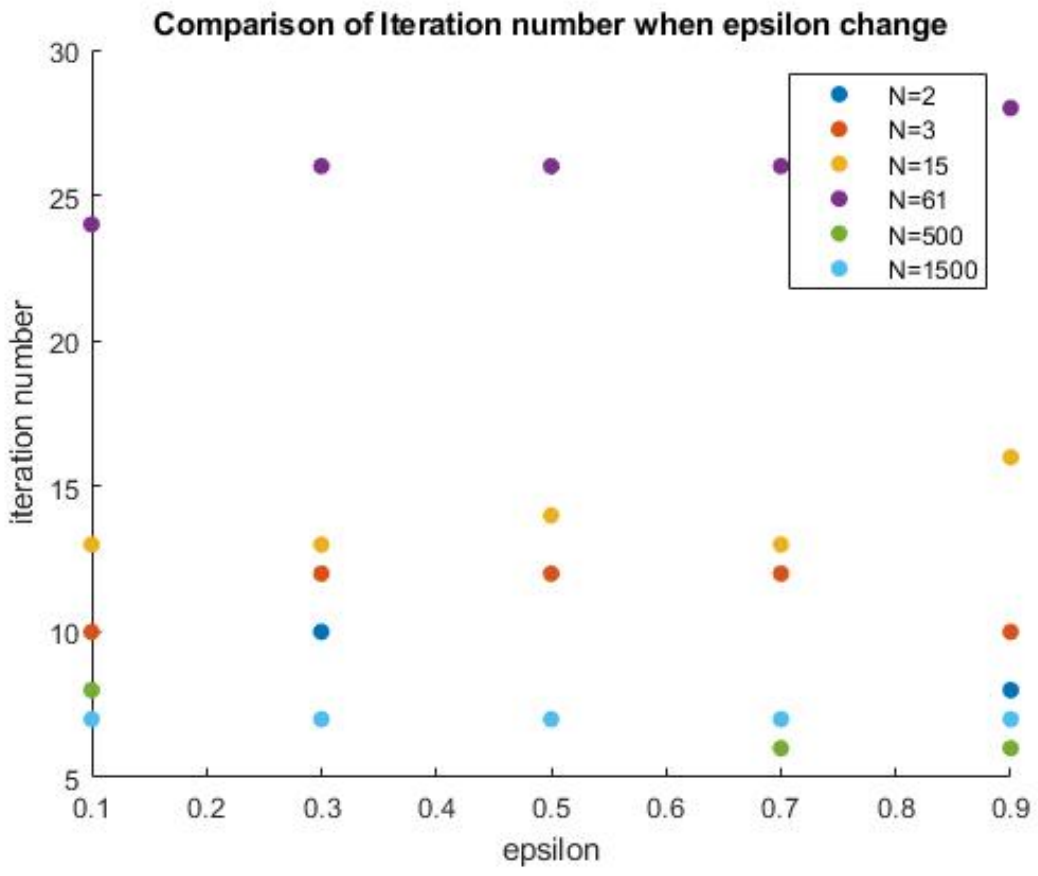

(c) For the third choice of the initial point

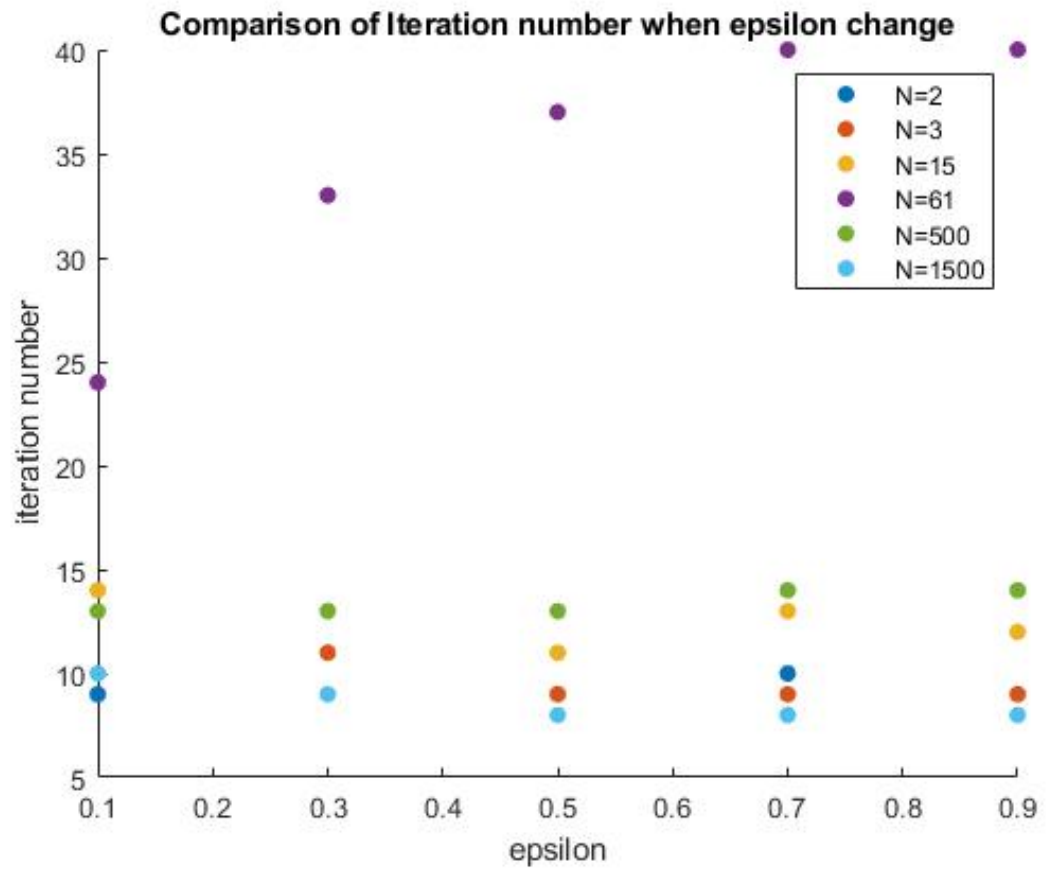

(d) Random of the initial point

We can see that it depend on cases, we choose a suitable initial point and $\epsilon$, for example with N = 61 we can choose the first choice of the initial point and epsilon near to 0, but when $\epsilon = \frac{1}{2}$ or 1 we can see that it's the value with taking the smallest number of iteration of a large number of cases except for the case N = 61. In general, we can choose that when we run various cases to predict the good iteration results we obtain.

The figure below shows that our methods are always convergence with any initial point.

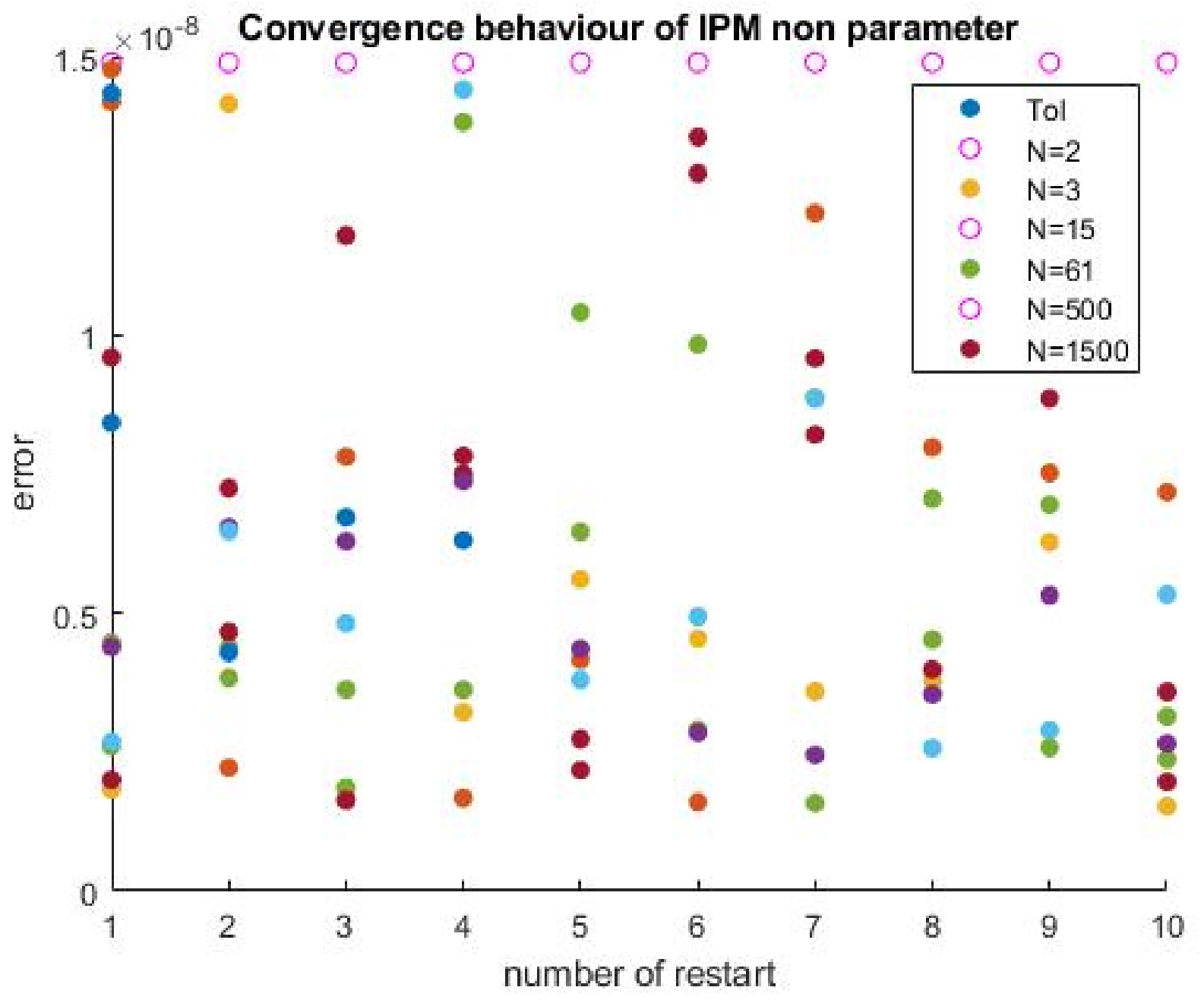

Figure 1.22: The convergence of non-parametric IPM for various option N

But the time to run it is quite high. Takes 161.4375s to run 200 times restart to run option case 2.

Secondly, we propose a comparison with difference safety parameter $\theta$

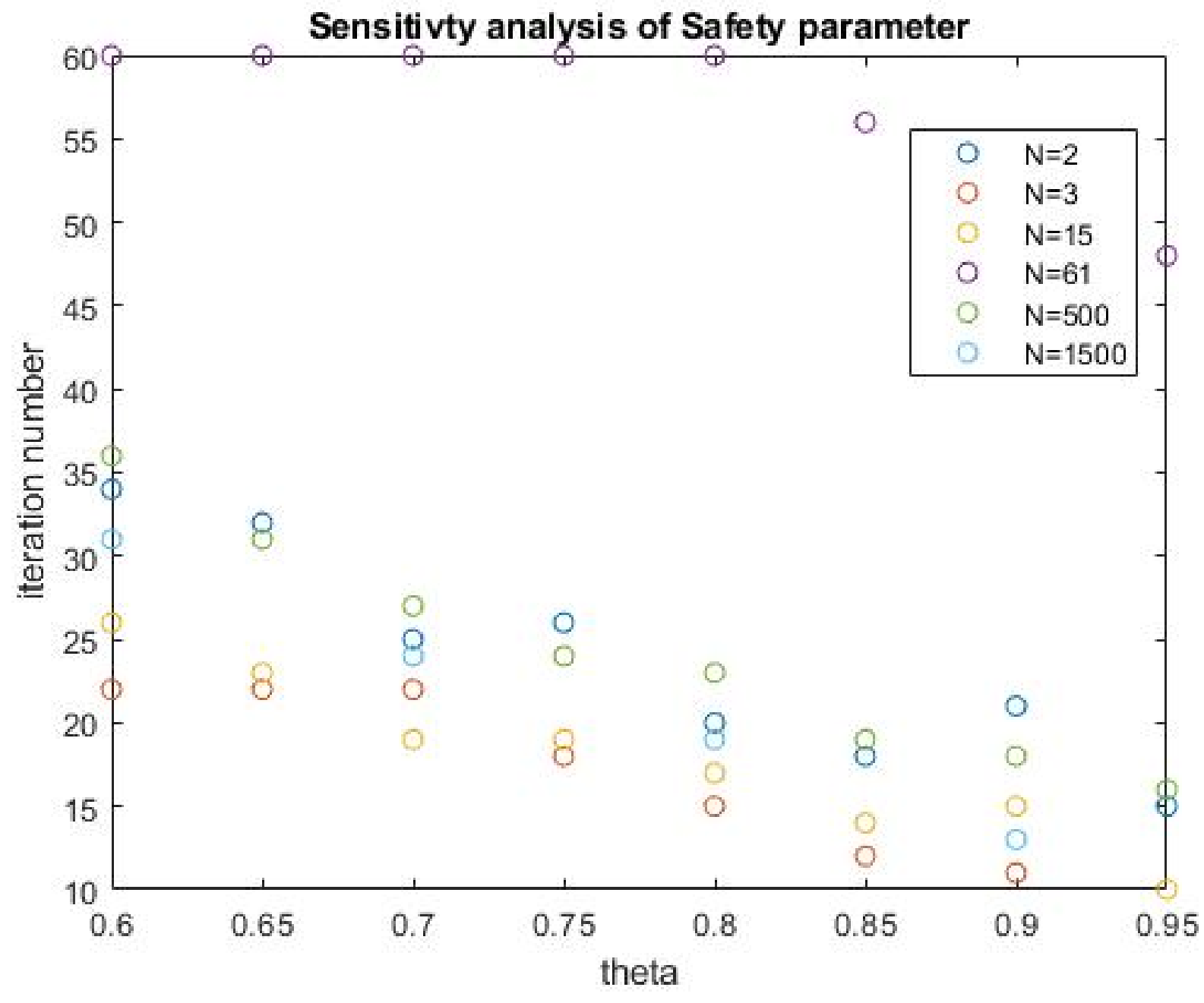

Figure 1.23: Iteration of non-parametric IPM for various option N

Now based on the result we can see that the safety parameter near 1 takes less iteration than the far from 1. I think it's good if we take a safety parameter near 1.

### 1.5.2 Numerical Result of non-parametric Interior Point method with Armijo line search

The figure below shows that our methods are always convergence with any initial point.

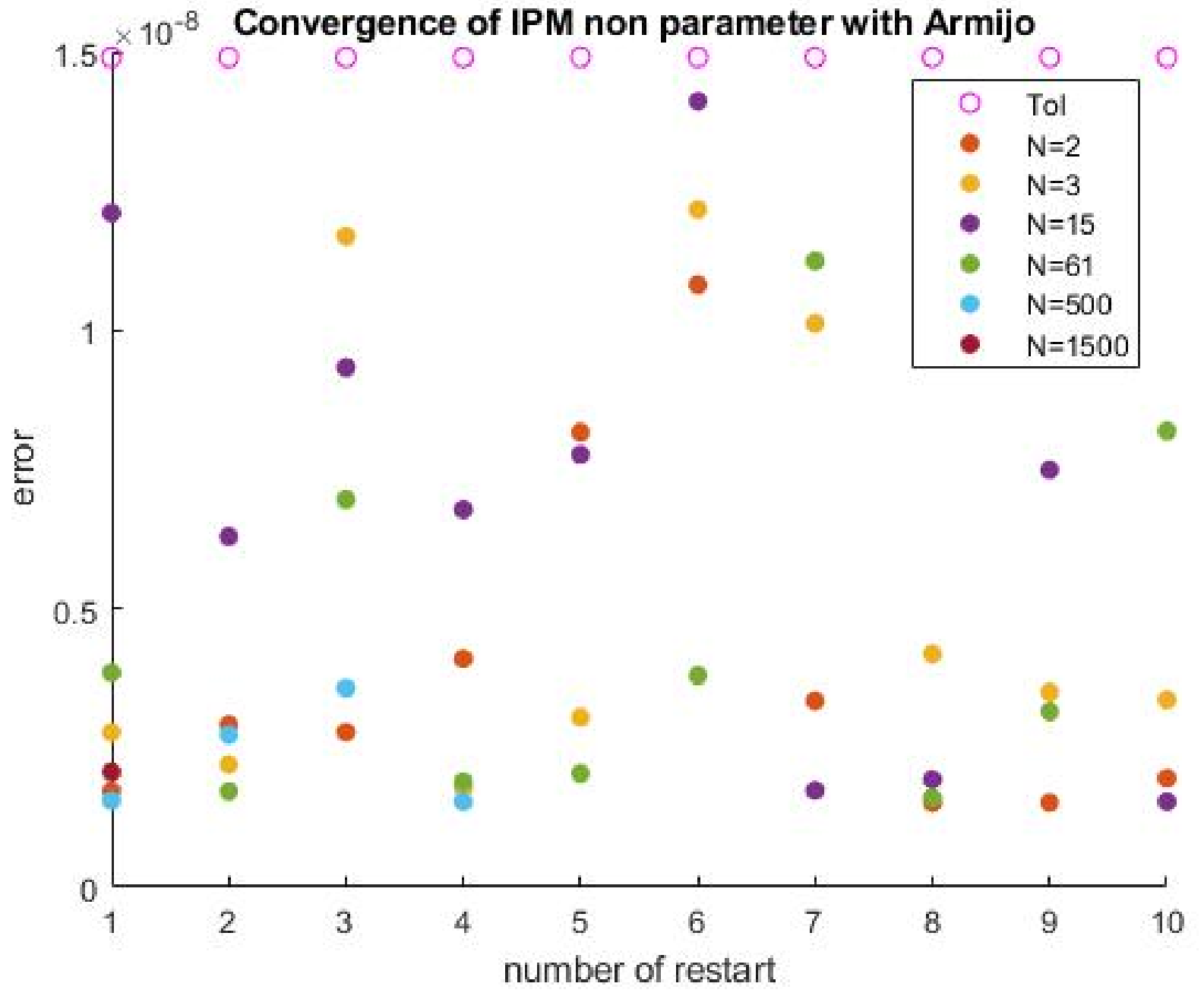

Figure 1.24: The convergence of non-parametric IPM for various option N

But the time to run it is quite high. Takes 161.1875s to run 200 times restart to run option case 2. To compare with IPM parameter it takes 245.0469s to run 200 times restart. If we use **"non-parametric IPM "** we can save 83.8594s in compare with **"IPM with parameter"**

Secondly, we propose a comparison with difference safety parameter $\theta$

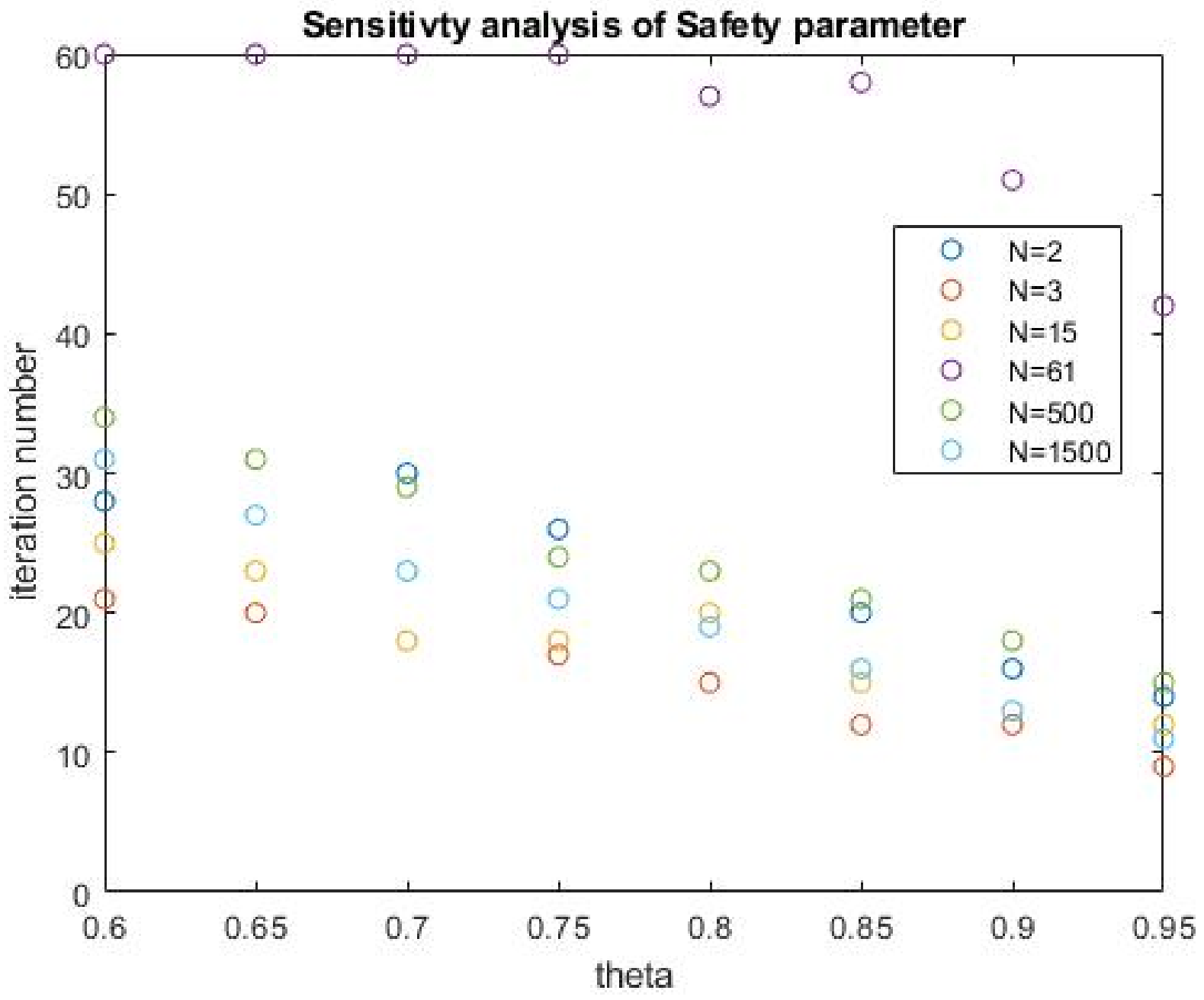

Figure 1.25: Iteration of non-parametric IPM for various option N

Now based on the result we can see that the safety parameter near 1 takes less iteration than the far from 1. I think it's good if we take a safety parameter near 1.

**Remark 1.19.** *One observation that it's not difference too much between **"IPM non parametric"** and **"IPM non parametric with Armijo line search"**. Define $z = |\alpha_{max} - \rho|$, with rho is the step size find in Armijo function with $t_0 = \alpha_{max}$. Run our function with various N we can see that*

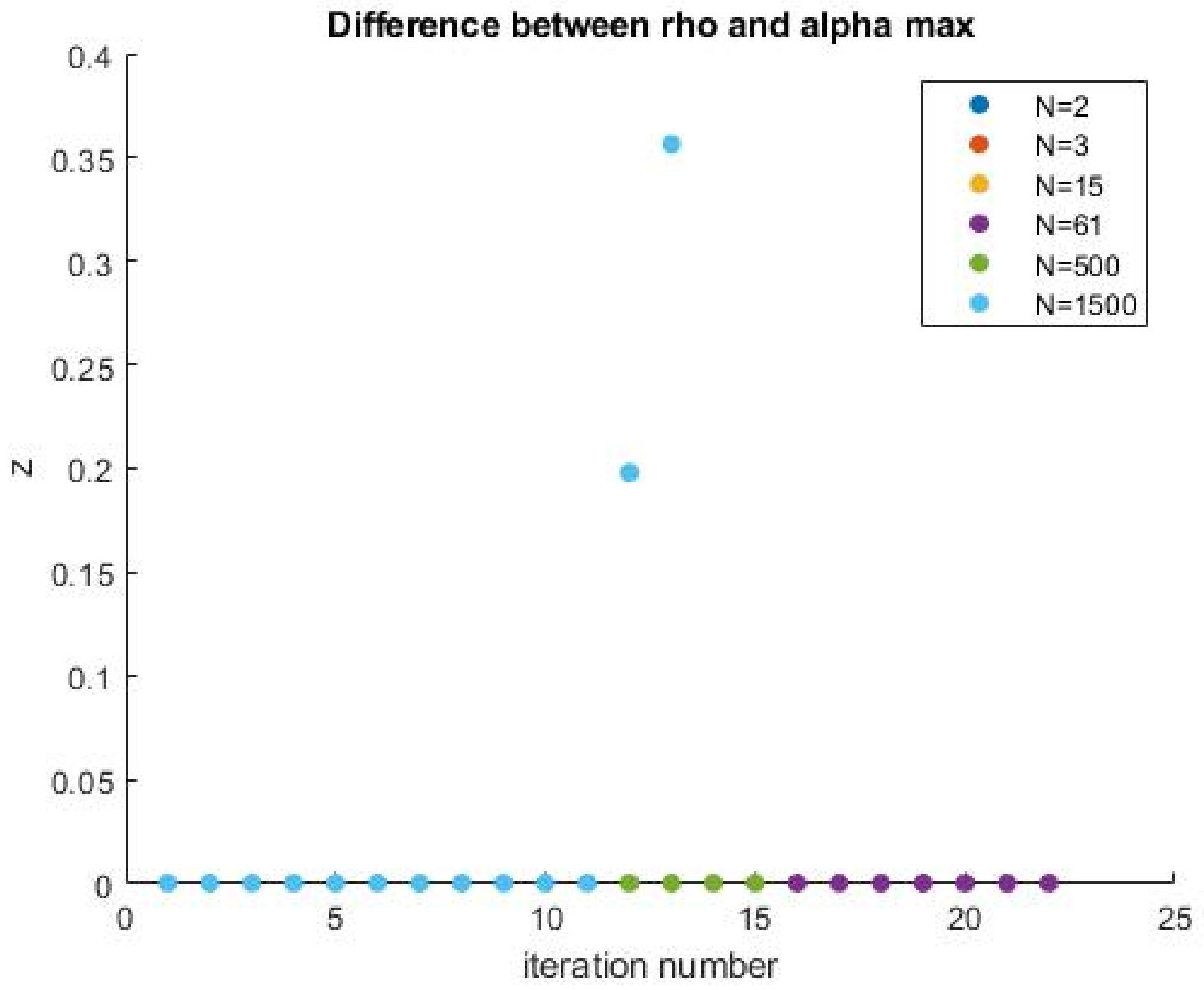

Figure 1.26: Difference between $\rho$ & $\alpha_{max}$ for various option N

***Now we can see the difference here****. Sometimes Armijo starts with* $t_0$ *and returns another value as can be seen in the case N = 1500.*

# Chapter 2

# A new approach for slack variable

In this chapter, we propose another place to put the slack variable which is more effective than the previous ones.

## 2.1 Interior-Point Method

The idea is to consider a sequence of approximate problems $P_\mu$ defined by

$$(\mathrm{P}_\mu)\begin{cases} (Au-b)_{3I-2} = \mathbf{0} \\ (Au-b)_{3I-1} = \mathbf{0} \\ u_{3I} \bullet (Au-b)_{3I} = \mu\mathbf{e} \end{cases} \tag{2.1}$$

for $1 \leqslant I \leqslant N$ and to drive $\mu > 0$ to 0 in some "smart" way.

To make things easier, let us introduce the slack variables.

$$w_I = (Au-b)_{3I} \in \mathbb{R}^N \tag{2.2}$$

system (2.1) then read

$$\begin{cases} (Au-b)_{3I-2} = \mathbf{0} \\ (Au-b)_{3I-1} = \mathbf{0} \\ (Au-b)_{3I} - w_I = \mathbf{0} \\ u_{3I} \bullet w_I = \mu\mathbf{e} \end{cases} \tag{2.3}$$

with 4N unknowns and 4N equations for a fixed $\mu$ for $1 \leqslant I \leqslant N$ .

**Ex. 3** — 1. Write down the Jacobian matrix and the Newton iteration for (2.1)

**Answer (Ex. 3)** — If we rearrange $(Au-b)_{3I-2}, (Au-b)_{3I-1}, (Au-b)_{3I} - w_I$ on top, I means the way we arrange the system (2.1) as follow

$$\begin{cases} (Au-b)_1 & = 0 \\ (Au-b)_2 & = 0 \\ (Au-b)_3 - w_1 & = 0 \\ (Au-b)_4 & = 0 \\ (Au-b)_5 & = 0 \\ (Au-b)_6 - w_2 & = 0 \\ ... & \\ ... & \\ (Au-b)_{3N-2} & = 0 \\ (Au-b)_{3N-1} & = 0 \\ (Au-b)_{3N} - w_N & = 0 \\ u_3 w_1 & = \mu \\ u_6 w_2 & = \mu \\ ... & \\ ... & \\ u_{3(N-1)} w_{N-1} & = \mu \\ u_{3N} w_N & = \mu \end{cases} \tag{2.4}$$

then we get the Jacobian as follow:

$$DF(X) = \begin{bmatrix} A & B \\ C & D \end{bmatrix} \in \mathbb{R}^{m+n,m+n} \tag{2.5}$$

where A $\in \mathbb{R}^{m,m}$

$$E = \begin{bmatrix} 0 & 0 & -1 & ... & 0 & 0 & & & \\ 0 & 0 & 0 & 0 & 0 & -1 & 0... & 0 & 0 \\ ... & & & & & & & & \\ ... & & & & & & & & \\ 0 & 0 & ... & 0 & 0 & 0 & 0 & 0 & 0 \\ 0 & 0 & ... & 0 & 0 & 0 & 0 & 0 & 0 \\ 0 & 0 & ... & 0 & 0 & 0 & 0 & 0 & -1 \end{bmatrix} \in \mathbb{R}^{n,m} \tag{2.6}$$

and **B = E'**, B $\in \mathbb{R}^{m,n}$(for simplification purpose only, because B is the same as C just take the transpose).

$$C = \begin{bmatrix} 0 & 0 & w_1 & ... & 0 & 0 & & & \\ 0 & 0 & 0 & 0 & 0 & w_2 & 0... & 0 & 0 \\ ... & & & & & & & & \\ ... & & & & & & & & \\ 0 & 0 & ... & 0 & 0 & 0 & 0 & 0 & 0 \\ 0 & 0 & ... & 0 & 0 & 0 & 0 & 0 & 0 \\ 0 & 0 & ... & 0 & 0 & 0 & 0 & 0 & w_N \end{bmatrix} \in \mathbb{R}^{n,m} \tag{2.7}$$

$$D = \begin{bmatrix} u_3 & 0 & 0 & ... & 0 & 0 \\ 0 & u_6 & 0 & ... & 0 & 0 \\ 0 & 0 & u_9 & 0 & ... & 0 \\ ... \\ ... \\ ... \\ 0 & 0 & ... & u_{3N-6} & 0 & 0 \\ 0 & 0 & ... & 0 & u_{3N-3} & 0 \\ 0 & 0 & ... & 0 & 0 & u_{3N} \end{bmatrix} \in \mathbb{R}^{n,n} \tag{2.8}$$

.

Newton iteration, firstly we choose $u^{(0)}, w^{(0)}, \mu^{(0)}$

.we have the system to solve is F(x) = 0

.compute the newton direction : $d_N(x) = -DF(x)^{-1}F(x)$

.update $x^{(k+1)} = x^{(k)} + \alpha_{max} d_N(x)$

**Remark 2.1.** *As in before we make a sensitivity analysis about the system 4N to choose a good way about safety parameter.*

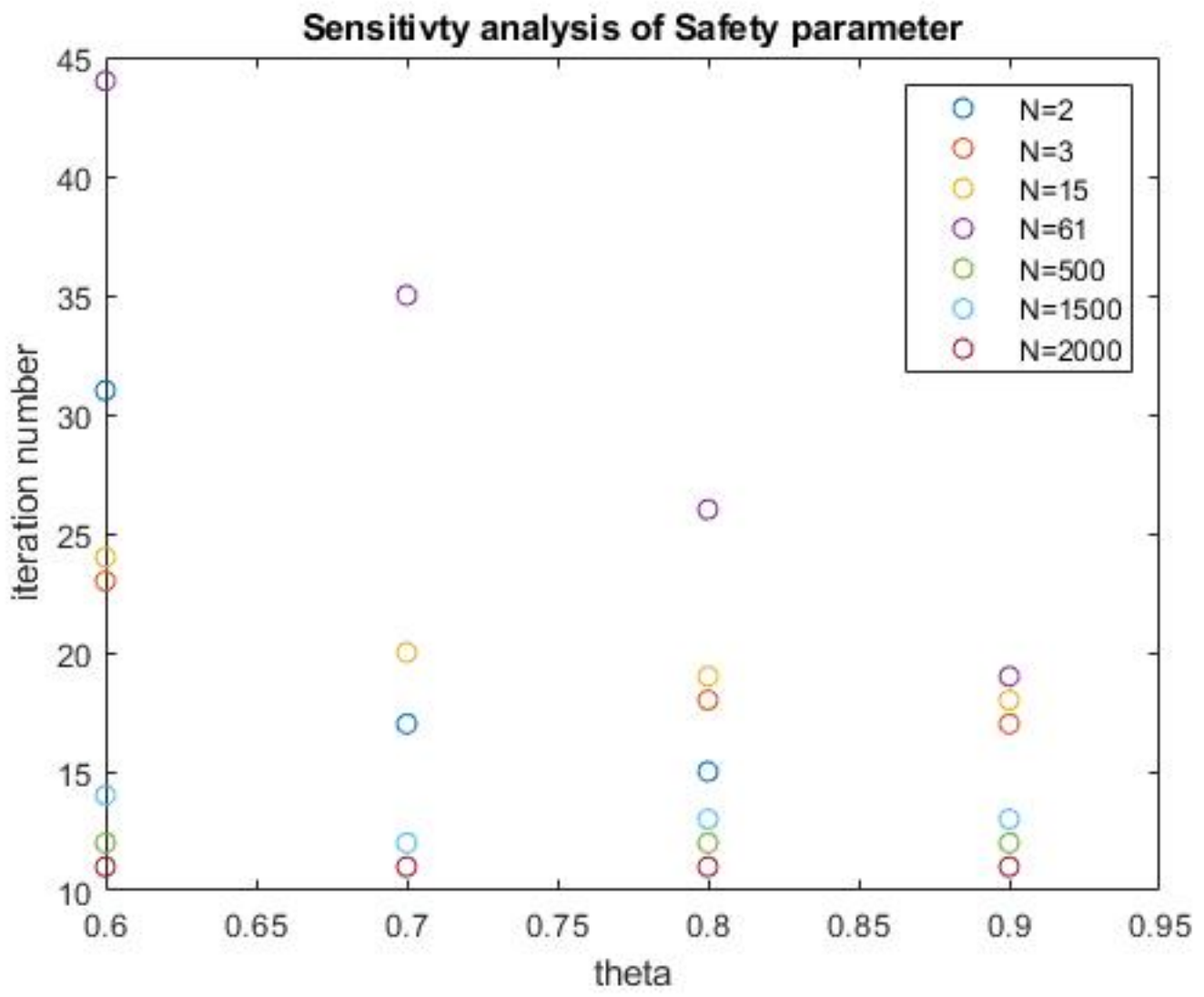

Figure 2.1: Iteration comparison of for various option N

*From this, we choose the safety parameter to be near 1 throughout this thesis.*

**Remark 2.2.** *To run IPM, we need a good initial point(I used here $\theta = 0.999$ (safety parameter)). Currently, There are three ways to choose an initial point.*

*The idea is when we compute F(X) = [Au-b-w;w]. Which make F has as much as possible number components of F is 0 (if we choose w as above when we run the code, approximation half components of F(X) = 0). It good because we can compute Newton direction $-DF(X)^{-1}F(X)$ faster.*

*Firstly, The same with the Standard method denoted*

$$\tilde{b_i} = \begin{cases} b_i & \text{if } i \neq 3I \\ max\{b_i, 0\} & \text{if } i = 3I \end{cases} \tag{2.9}$$

*for* $1 \leqslant I \leqslant N$ *. But after we solve Du=*$\tilde{b}$*. Some components of u is 0 ( because if* $b_i < 0$ *then* $\tilde{b}_i = 0$ *). At that time we touch the boundary* ***(violated the principle of IPM, we start with an initial guess inside the domain and then converge to a solution )***

*So if* $\exists j$ *s.t* $u_{3j} = 0$*, we need a safety scale to ensure that it don't touch the boundary. we take* $u_{3I} = max\{u_{3I}, 0.01\}$*. After that we take w = Au-b. Some components of w is negative ,suppose at the components j-th,* $w_j < 0$ *there are two ways*

1. $w_j = 0.01 w_j$ *which we make the negative value become small*

2. $w_j = 0.01$ *which make it a small positive number*

*tested two ways above, all of them is doing a good work also (make the first error become small ). But with* $w_j = 0.01$*, the error is smaller than* $w_j = 0.01w_j$ *because when components j-th is negative we then we replace this components w near to 0 in the vector F(X) = [Au-b-w;w], The second part of F(X) has w so we scale j-th components negative become small positive number to the second part have the components 0. To have as much as possible components of F(X) is zero. Also when we use* $w_j = 0.01$ *we can run non-parametric IPM. So i recommend use* $w_j = 0.01$ *instead of* $w_j = 0.01w_j$

*Secondly, we choose infeasible IPM point, that's is the point satisfies the complementarity equations.*

$$w_I u_{3I} = \mu, \forall 1 \leqslant I \leqslant N \tag{2.10}$$

*The idea of implementing it as follows:*

1. *choose initial* $\mu > 0$

2. *choose u randomly s.t* $u_{3I} > 0$

3. *choose w randomly, then setting* $w_{3I} = \frac{\mu}{u_{3I}}$

*then* $w_{3I}u_{3I} = \mu, \forall 1 \leqslant I \leqslant N$ *But it doesn't satisfy the others equations, but we don't need to think about it.*

*While the second choice satisfies the complementary conditions, the third way we choose initial point satisfying*

$$\begin{cases} (Au-b)_{3I-2} & = 0 \\ (Au-b)_{3I-1} & = 0 \end{cases} \tag{2.11}$$

$\forall 1 \leqslant I \leqslant N$*, except when* $u_{3I} < 0$ *then we scale it become* $u_{3I} = max\{u_{3I}, 0.01\}$

*The idea of the third way is almost the same with the first ways, but in the first way we don't need to compute* $u = A \setminus \tilde{b}$ *with is expensive than* $u = D \setminus \tilde{b}$*.*

*The idea of implementing it as follows:*

1. $w > 0$

2. *Calculate*

$$\tilde{b}_j = \begin{cases} b_{3I} + w_{3I} & \text{if } j = 3I \\ b_j & \text{if } j \neq 3I, \forall 1 \leqslant I \leqslant N \end{cases}$$

,

3. *Solving* $u = A \setminus \tilde{b}$

4. *Setting* $u_{3I} = max\{u_{3I}, 0.01\}$

*Fourthly, we choose* $u^{(0)}$ *s.t,* $u^{(0)}_{3I} > 0$,$(Au^{(0)} - b)_{3I} > 0$ *then we set* $w = (Au^{(0)} - b)_{3I}$. *Then F(X) = [Au-b;w]* ,$F(X)_{3I} = (Au - b)_{3I} - w_I = 0$. *Here we try to make components of F(X) is 0 as most as posible, but that point doesn't always exists, for example consider the system*

$$\begin{cases} x - y + 1 & = 0 \\ -2x + y - 2 & = 0 \end{cases} \tag{2.12}$$

*we want* $x, y > 0$ *and* $x - y + 1 > 0, -2x + y - 1 > 0$ *, Or*

$$\begin{cases} x & > y - 1 \\ -2x & > 2 - y \end{cases} \tag{2.13}$$

*Or*

$$\begin{cases} x & > y - 1 \\ x & < \frac{y}{2} - 1 \end{cases} \tag{2.14}$$

*Or*

$$\frac{y}{2} - 1 > x > y - 1$$

*so we get* $\frac{y}{2} > y$ *Or* $\frac{y}{2} < 0$ *Or* $y < 0$ *(a contradiction)*

*An idea as follow. Solve* $u = A \setminus b$ *, then setting*

$$v_i = \begin{cases} 1 & \text{if } i = 3I \\ 0 & \text{if } i = 3I - 2, 3I - 1 \end{cases} \tag{2.15}$$

*check if* $(Av)_{3I} > 0$, *then* $u_{new} = u+v$, *check* $u_{3I} > 0$ *if not increse* $u = u+2v$ *to make sure* $u_{3I} > 0$. *After that check* $(Au - b)_{3I} > 0$. *This method can work on the case N = 2, N = 3, N= 15, N = 61, N = 500, N = 2000. So i denoted way 4 in the cases this technique can work, and make a comparison of them.*

*I don't make a comparison of the way to choose initial point for the system 6N because the system 6N is just an idea to put the slack variable. When we introduce the system 4N, it takes less memory and efficient also (convergence, iteration, time consumption).*

*Here I denoted the first choice of the initial point (ways 1), the second choice of the initial point (ways 2), the third choice of the initial point (ways 3)*

*Firstly, let's talk about convergence, all of the ways above is convergence.*

*Secondly, about First error (the error after the first iteration)*

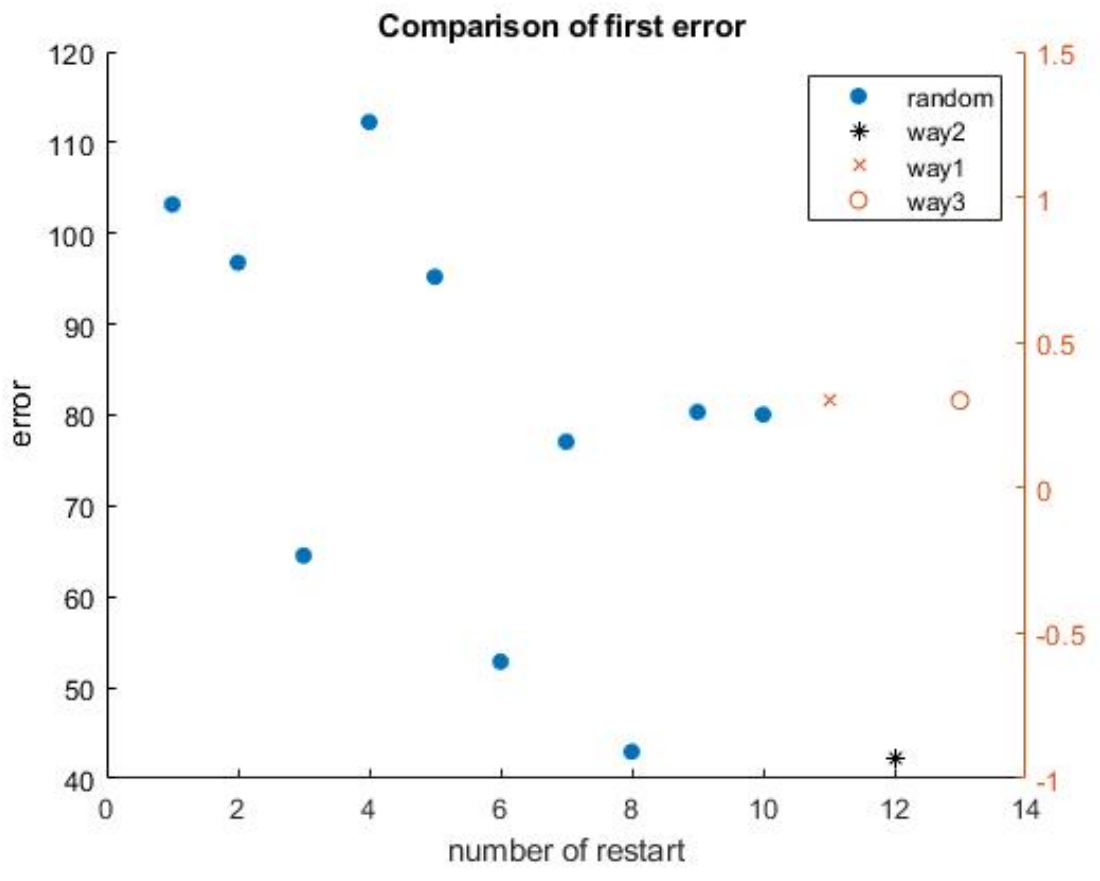

(a) $N = 2$

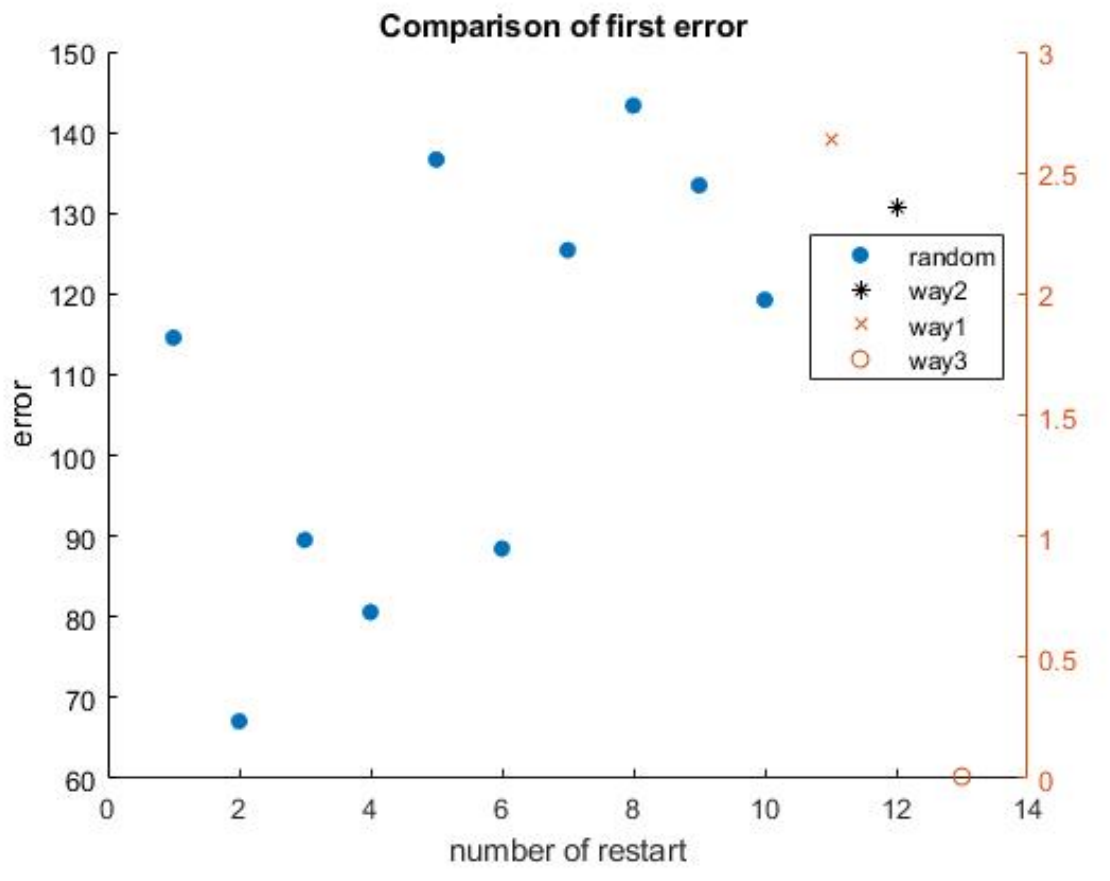

(b) $N = 3$

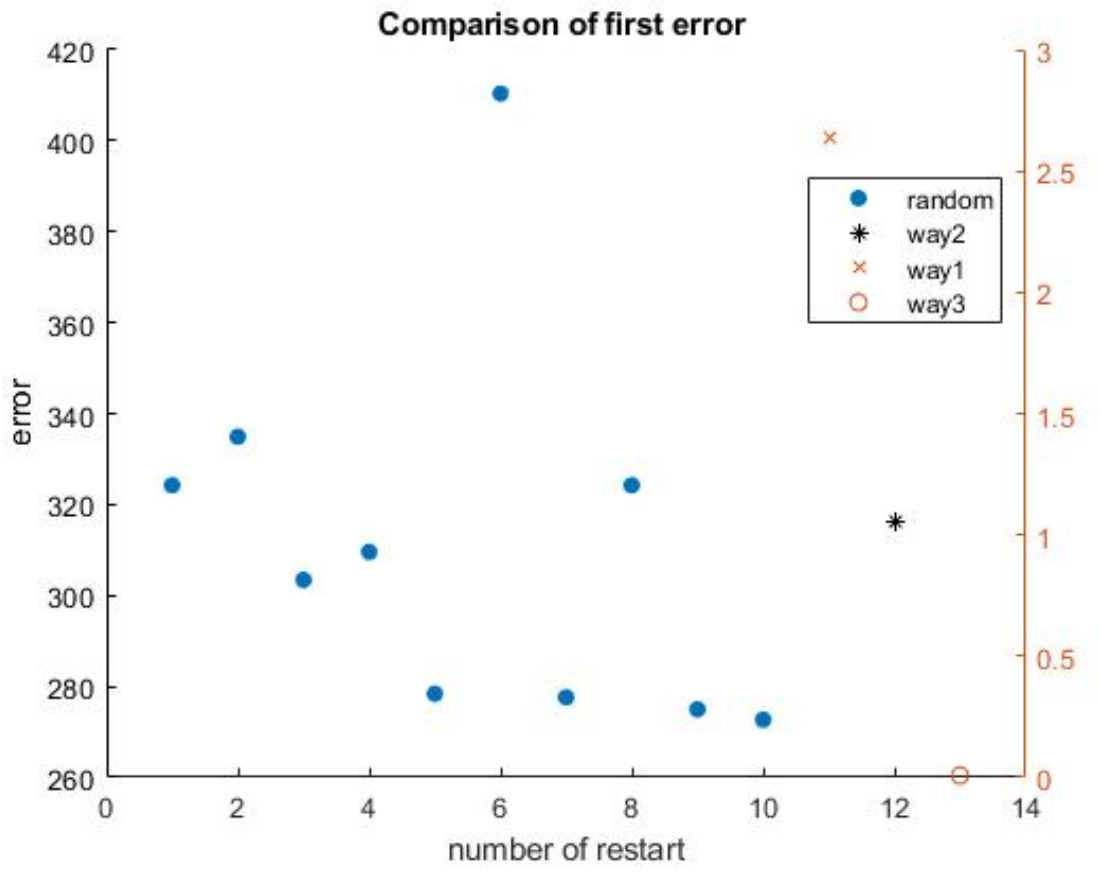

(c) $N = 15$

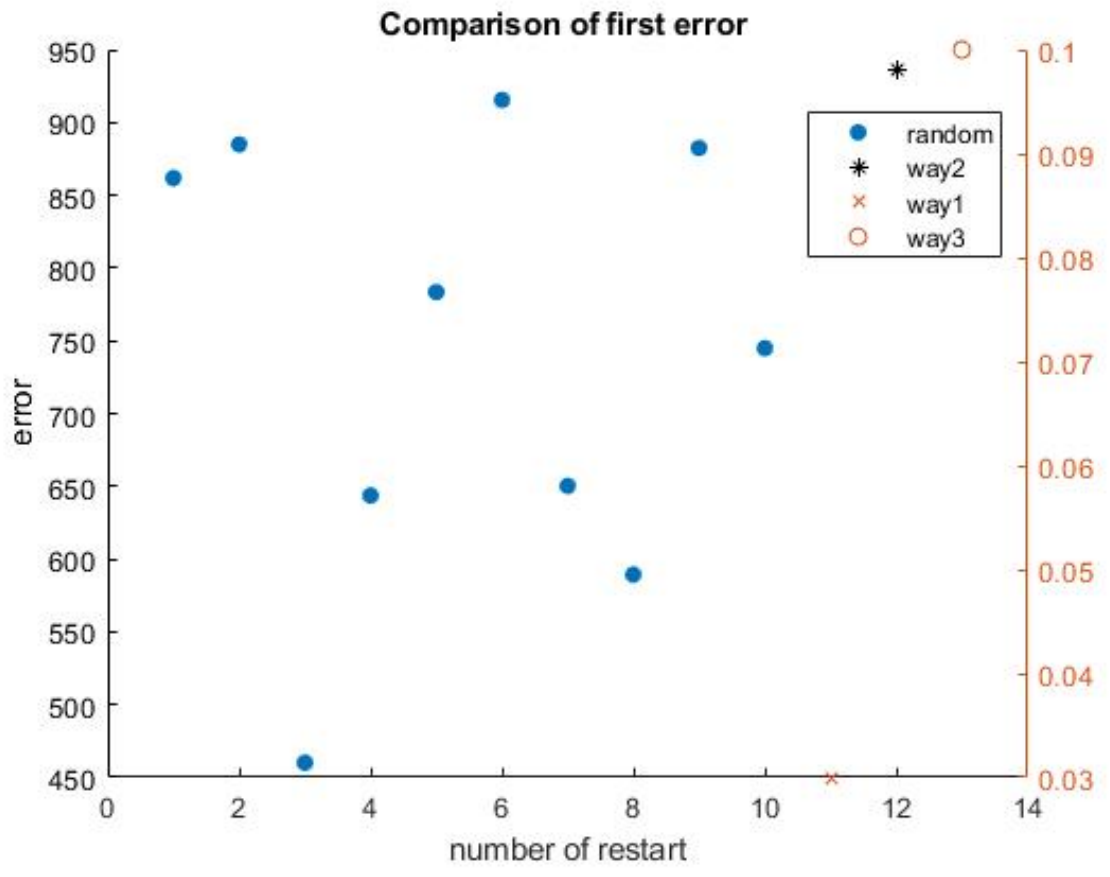

(d) $N = 61$

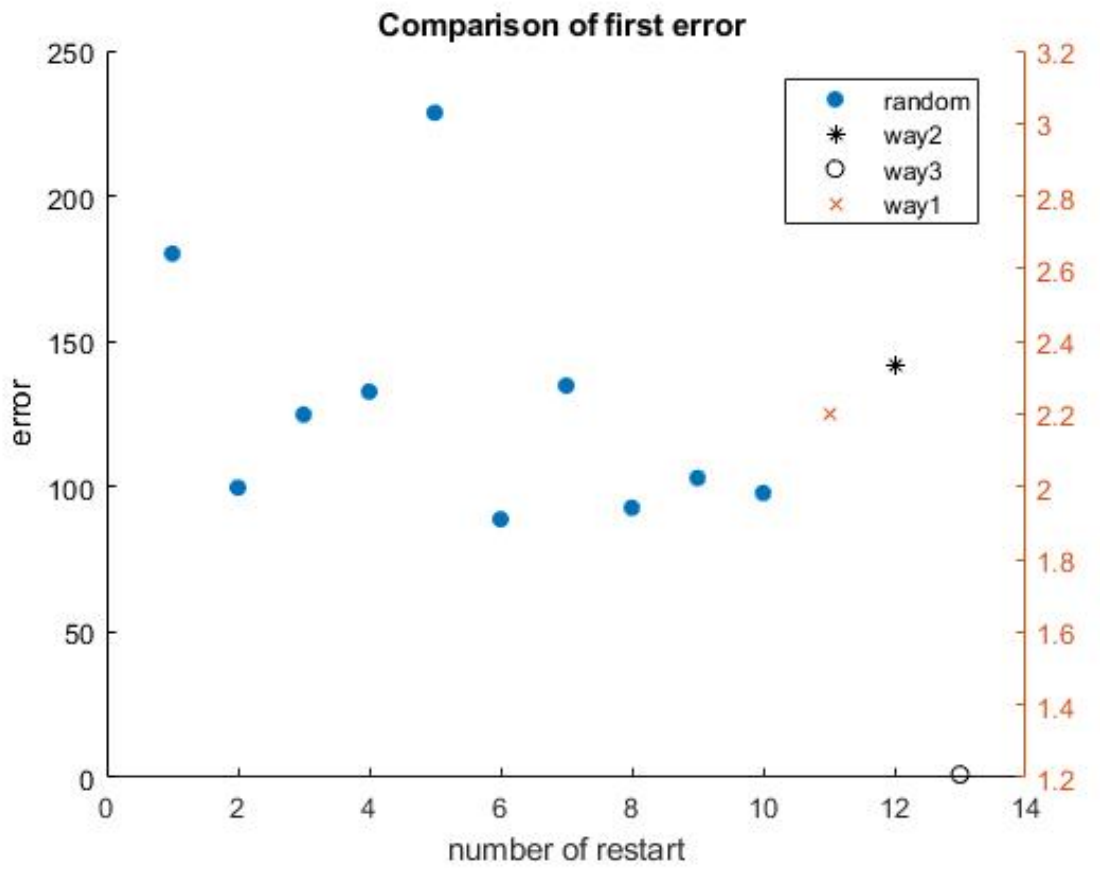

(e) $N = 500$

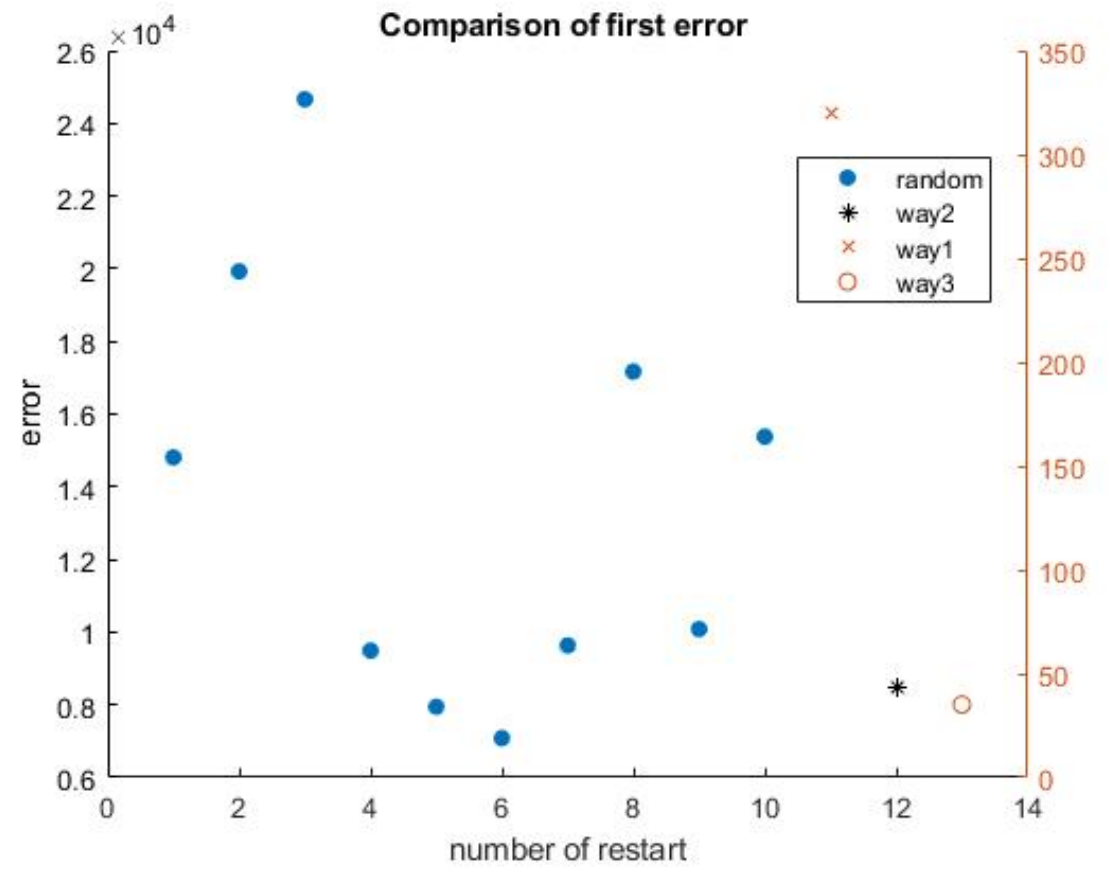

(f) $N = 1500$

Figure 2.2: First error Comparison

*As you can see in the picture, $theway3 > theway1 > theway2$, the way 3 wins about section first error, because the third way we solve $(Au - b)_{3I-2}, (Au - b)_{3I-1}$, so they almost satisfies the system of the equation except for the complementary condition, so they have the smallest error in 3 ways. The first way we just solve the approximation of the systems of it have the second*

*smallest. The ways 2 we just solve the "Perturbed" complementary problems so it has the error is largest in 3 ways but it better than the general ones, As can see in the picture, the way 2 met the condition "maximum iteration", it also decreases the norm of error, but slowly.*

*Thirdly, About iteration*

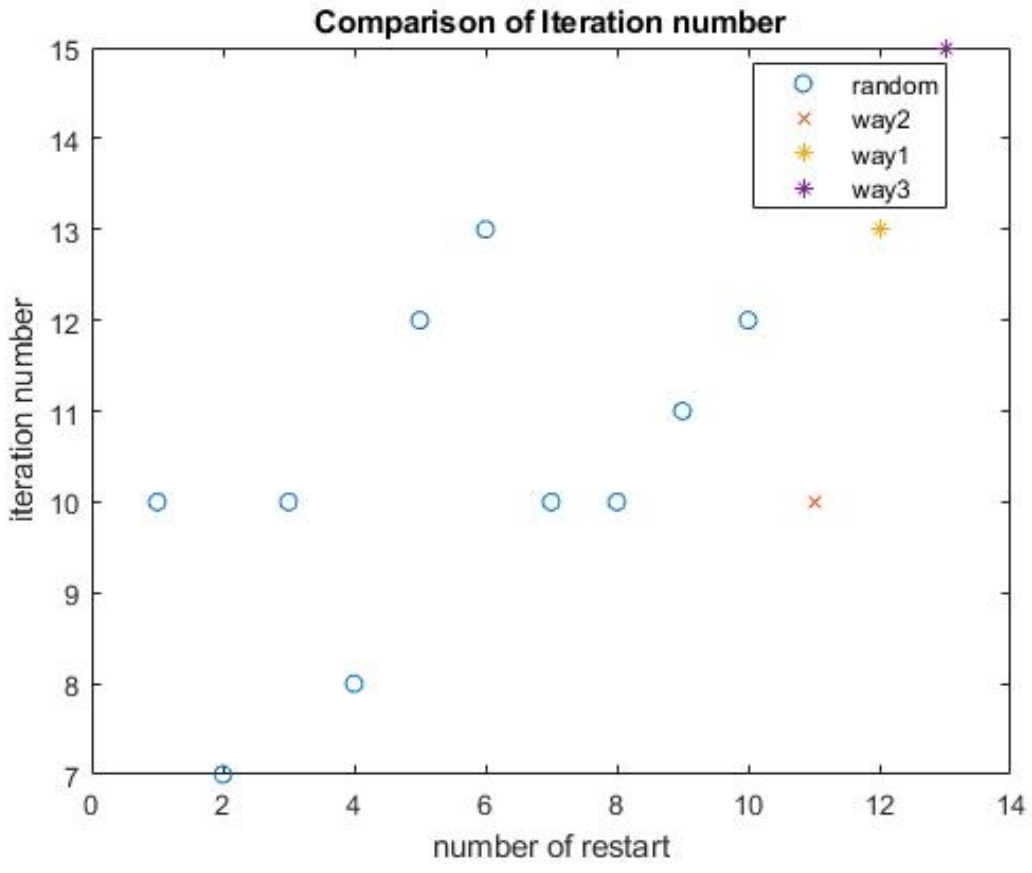

(a) $N = 2$

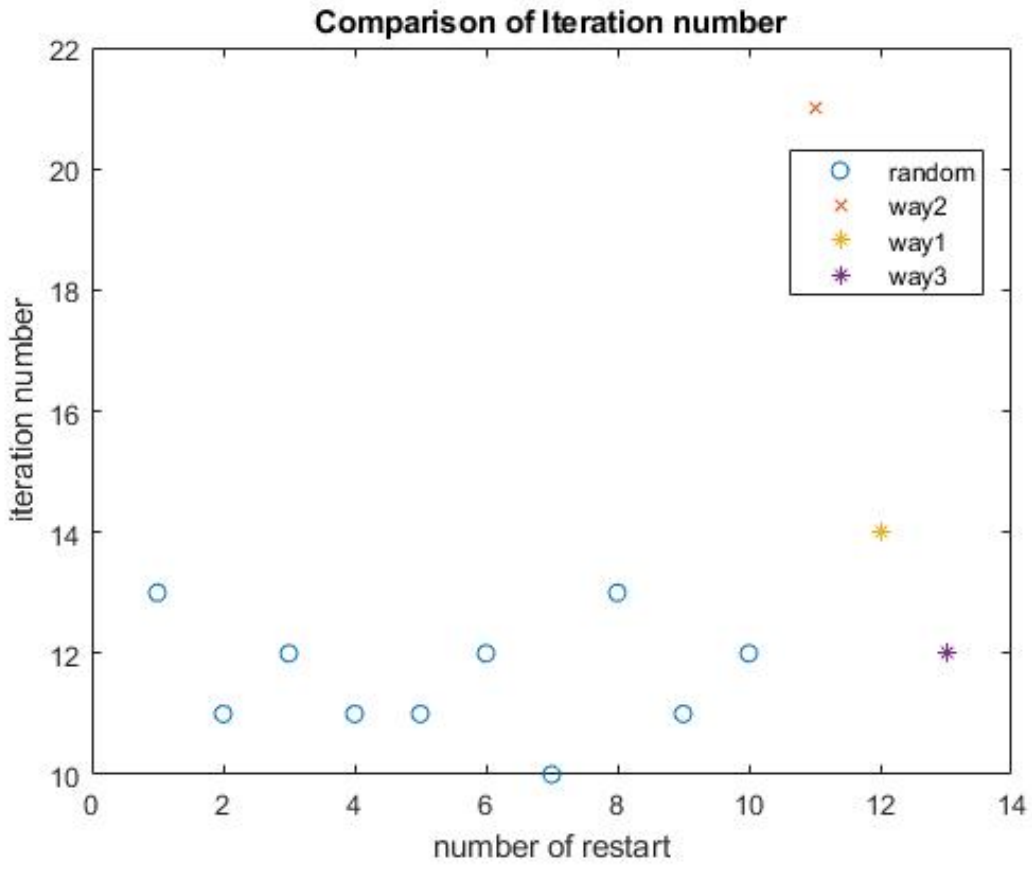

(b) $N = 3$

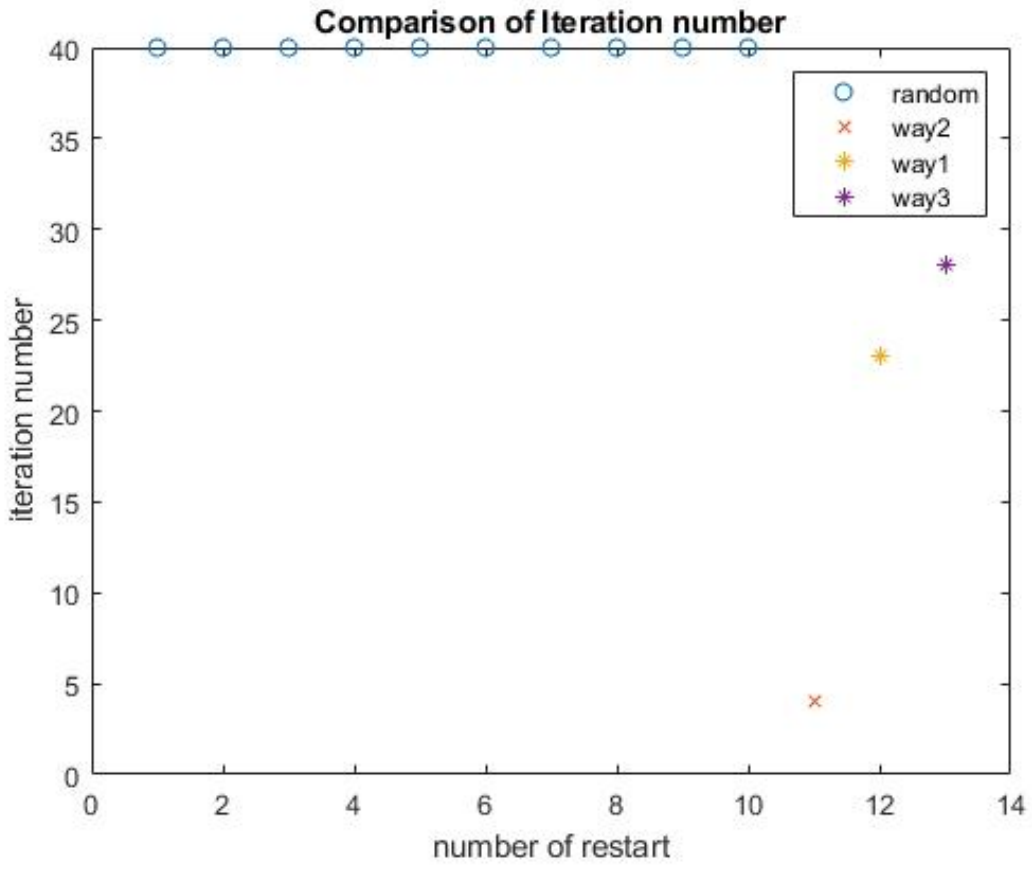

(c) $N = 15$

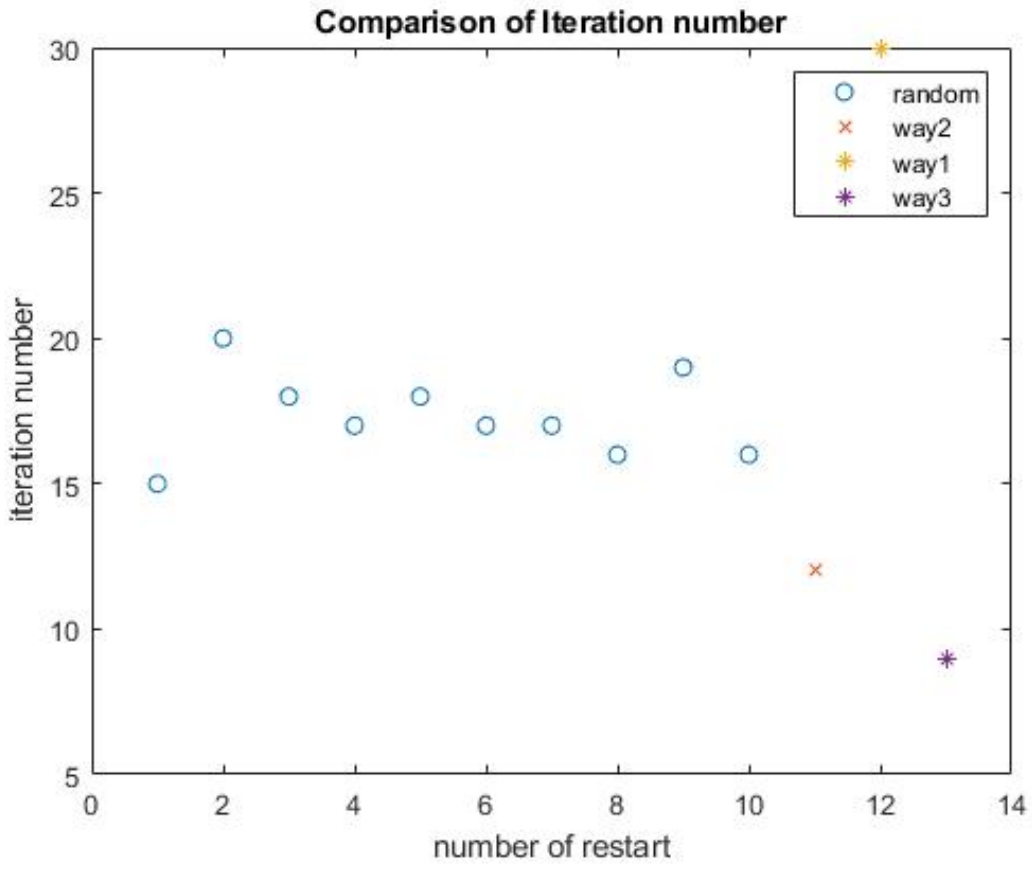

(d) $N = 61$

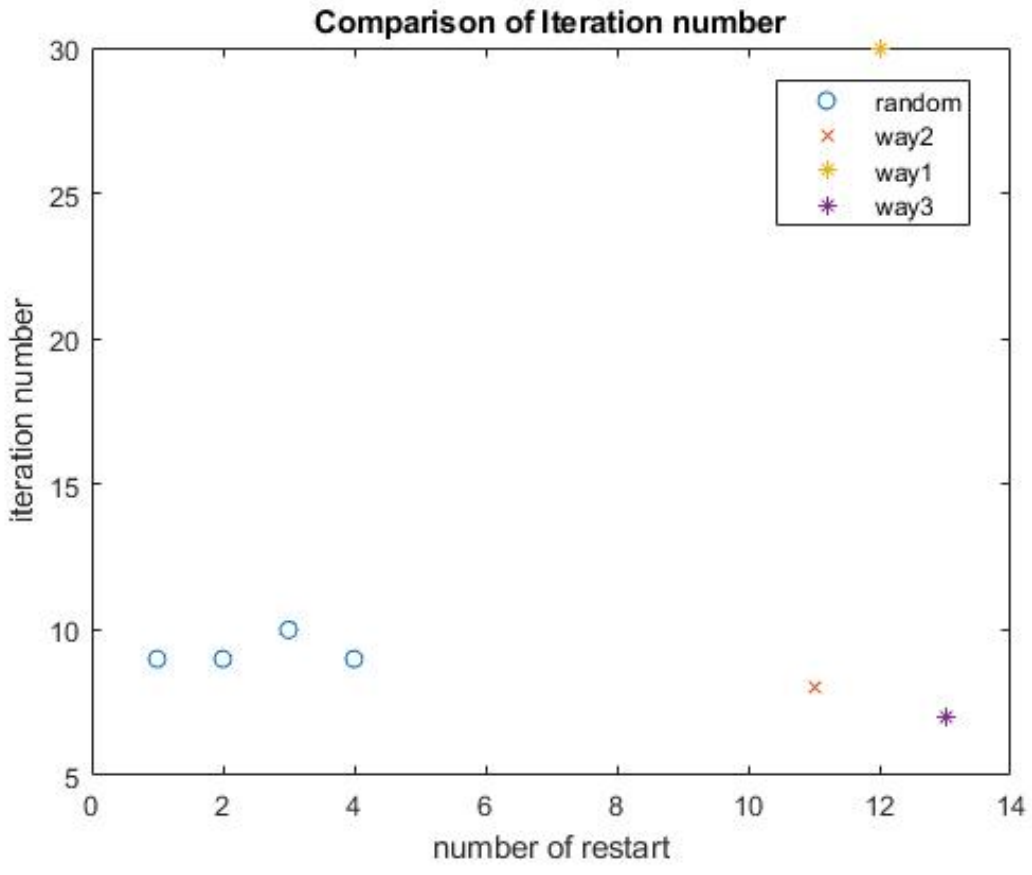

(e) $N = 500$

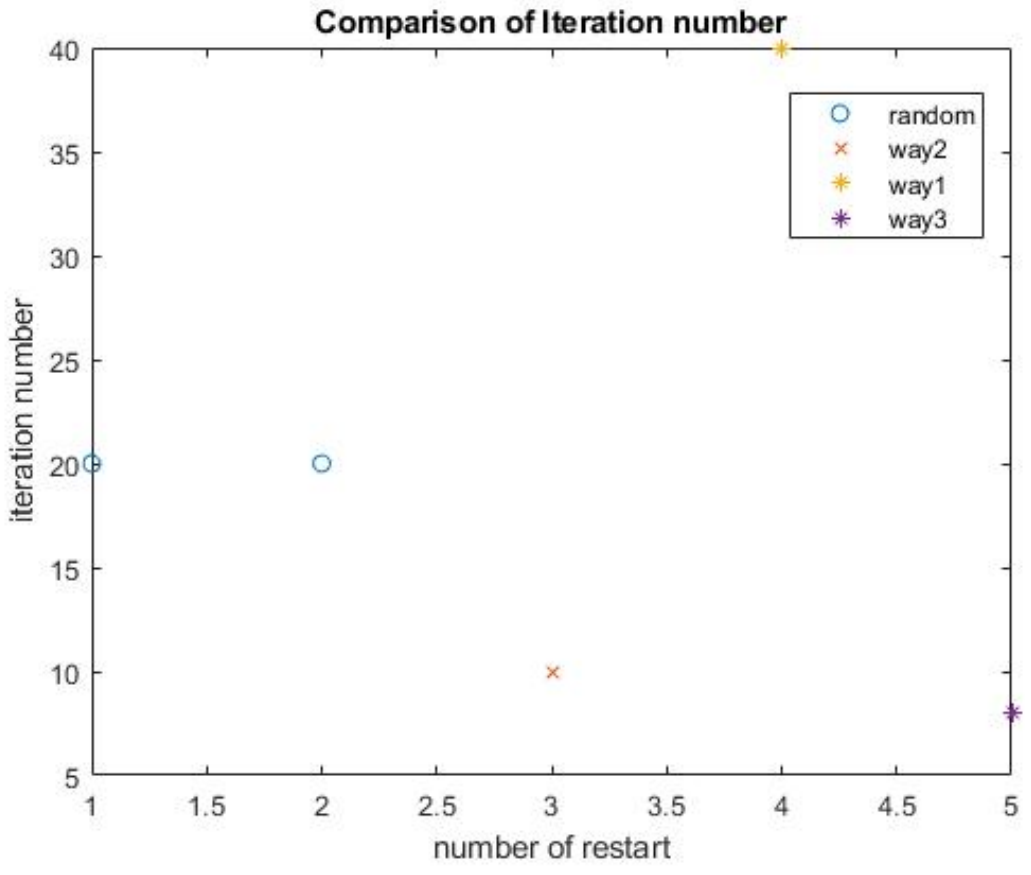

(f) $N = 1500$

Figure 2.3: Iteration Comparison

*easily see that* $theway3 > theway1 > theway2$ *at the iteration comparison. Moreover, It's better than when we choose* $u^{(0)}, w^{(0)} > 0$ *randomly.*

*The table below shows that way 3 also has less time consuming than the others.* ***Here I used CPU time (s)***

| N | 3 | 15 | 61 | 500 | 1500 | 2000 |
|---|---|---|---|---|---|---|
| randomly | 0.04 | 0.04 | 0.5 | 7.4 | 83 | 408.5 |
| the way 1 | 0.03 | 0.03 | 0.06 | 4.5 | 81 | 180 |
| the way 2 | 0.05 | 0.03 | 0.2 | 13.8 | 102 | 847 |
| the way 3 | 0.05 | 0.08 | 0.34 | 3.3 | 60 | 137 |
| the way 4 | 0.1 | 0.125 | 0.46 | 12.45 | ... | 508 |

*The table below shows that the way 3 also has less iteration that the way 4.*

| N | 2 | 3 | 15 | 61 | 500 | 1500 | 2000 |
|---|---|---|---|---|---|---|---|
| the way 3 | 15 | 12 | 27 | 9 | 7 | ... | 7 |
| the way 4 | 44 | 60 | 60 | 42 | 32 | ... | 29 |

*The table below shows that the way 3 also has the first error smaller than the way 4.*

| N | 2 | 3 | 15 | 61 | 500 | 1500 | 2000 |
|---|---|---|---|---|---|---|---|
| the way 3 | 0.37 | 0.25 | 0.25 | 0.1 | 1.2 | ... | 0.69 |
| the way 4 | 0.99 | 1.05 | 1.12 | 1.71 | 4.04 | ... | 3.36 |

*The ways 4 just has an idea but it decreases the value of error slowly, so it takes more time and iteration than the ways 3. So we choose the way 3 to be our initial point.*

### 2.1.1 Numerical results of IPM

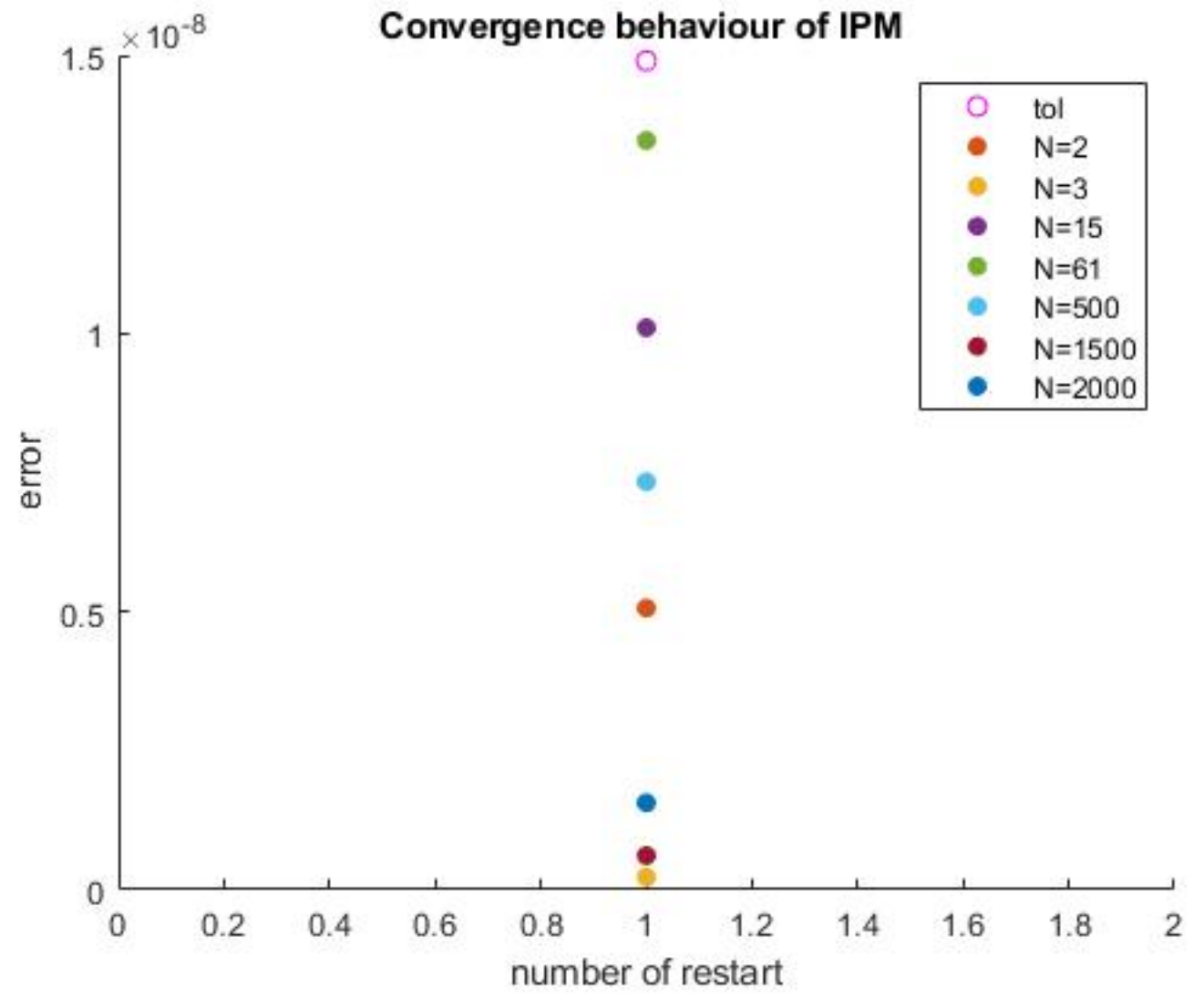

Figure 2.4: Convergence behavior of IPM for the various option N

### 2.1.2 Numerical resutls of IPM with Armijo line search

The same as before but we add Armijo line search to our code after computing $\alpha_{max}$.

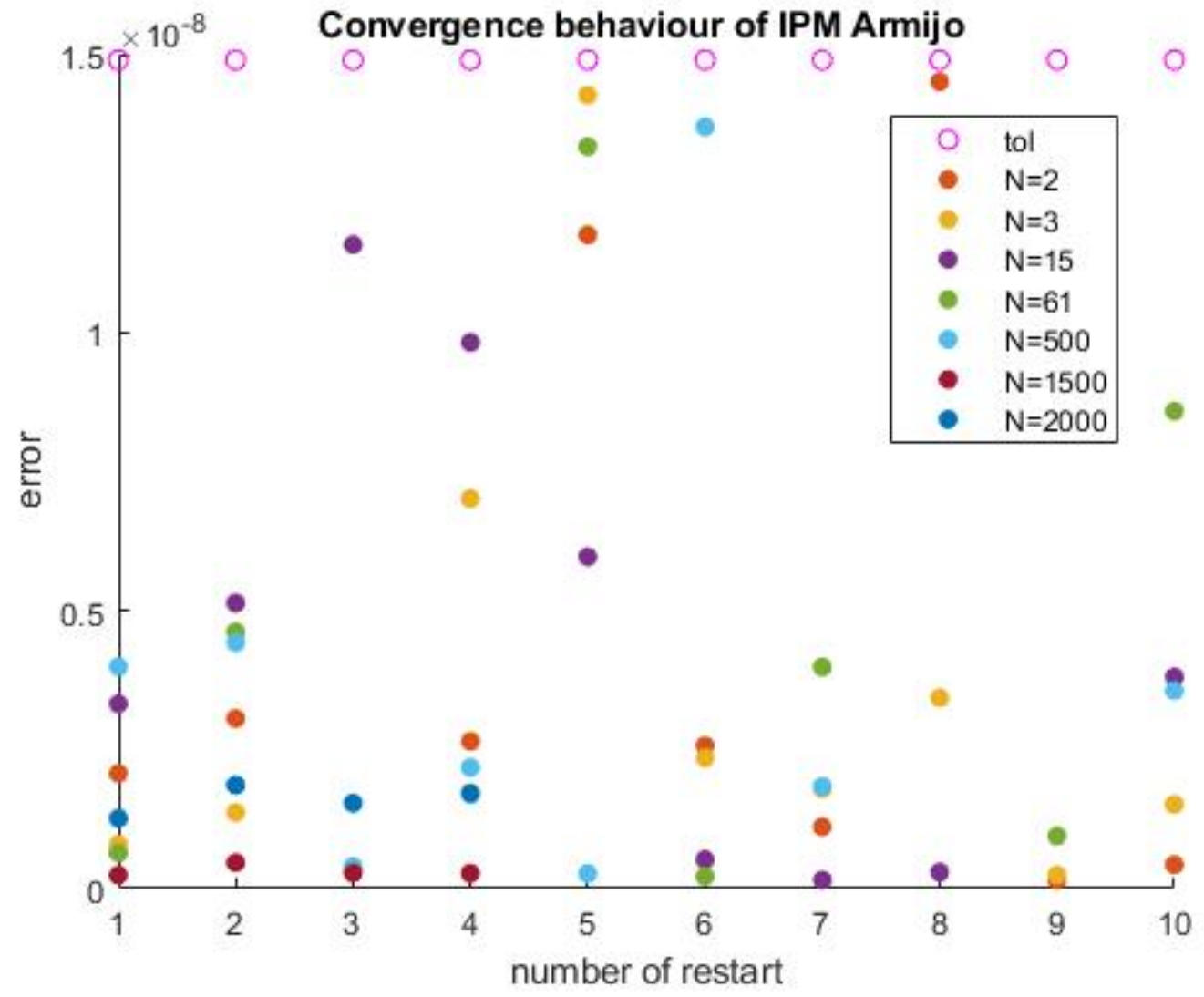

Figure 2.5: Convergence behavior of IPM for the various option N

**Remark 2.3.** *One observation that it's not difference too much between **"IPM "** and **"IPM with Armijo line search"**. Define $z = |\alpha_{max} - \rho|$, with rho is the step size find in Armijo function with $t_0 = \alpha_{max}$. Run our function we can see that*

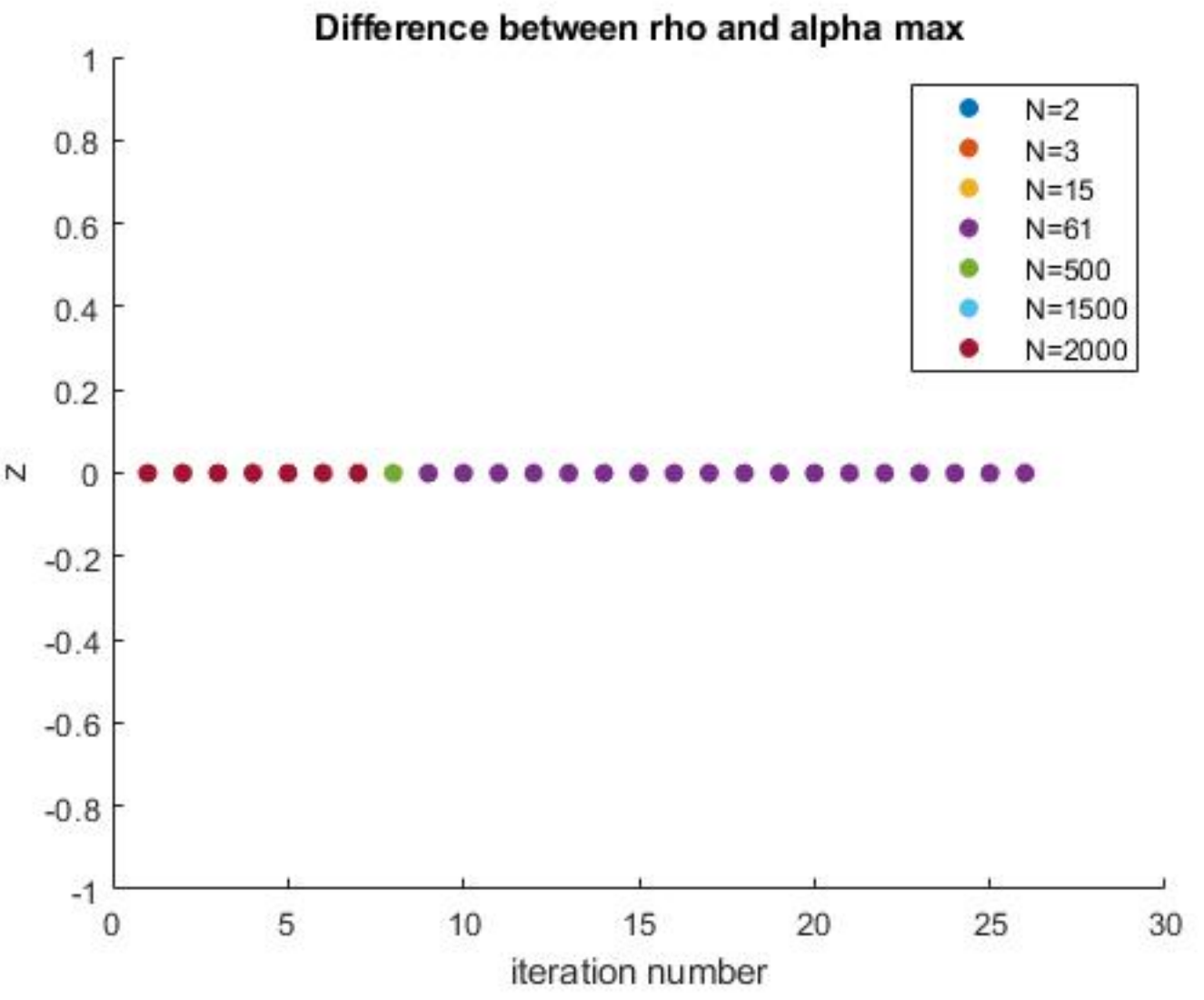

Figure 2.6: Difference between $\rho$ & $\alpha_{max}$ for various option N

*That means various option N. Armijo start with $t_0 = \alpha_{max}$ and return $t_0$. It's no meaning that Armijo doesn't always work. Sometimes it works sometimes not we don't know it.*

**Remark 2.4.** *We make a comparison of the iteration number between system 4N and system 6N. Here I used IPM methods. Most important, all the system 4N and 6N is convergence so we don't compare about that.*

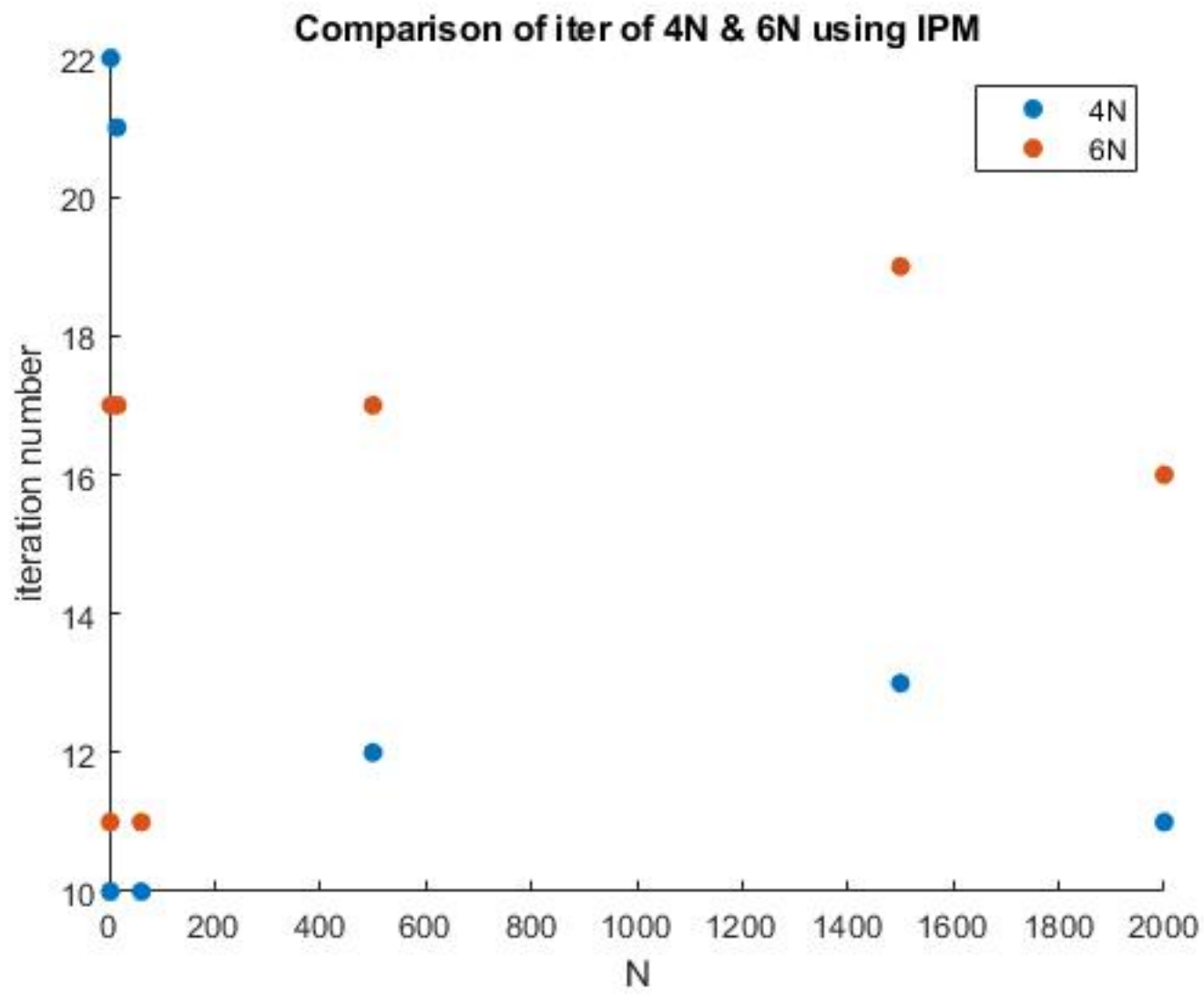

Figure 2.7: Iteration comparison between the system 4N & 6N

*When the size of the problems is small, it's not difference too much between 4N and 6N but when the size is large, we can see that 4N take less iteration than 6N*

*More attractive, we compare the time computation of IPM 6N and IPM 4N. Because it's not*

*difference too much when the size of the problems is small, so we can skip it, move to the size of the problems is large.*

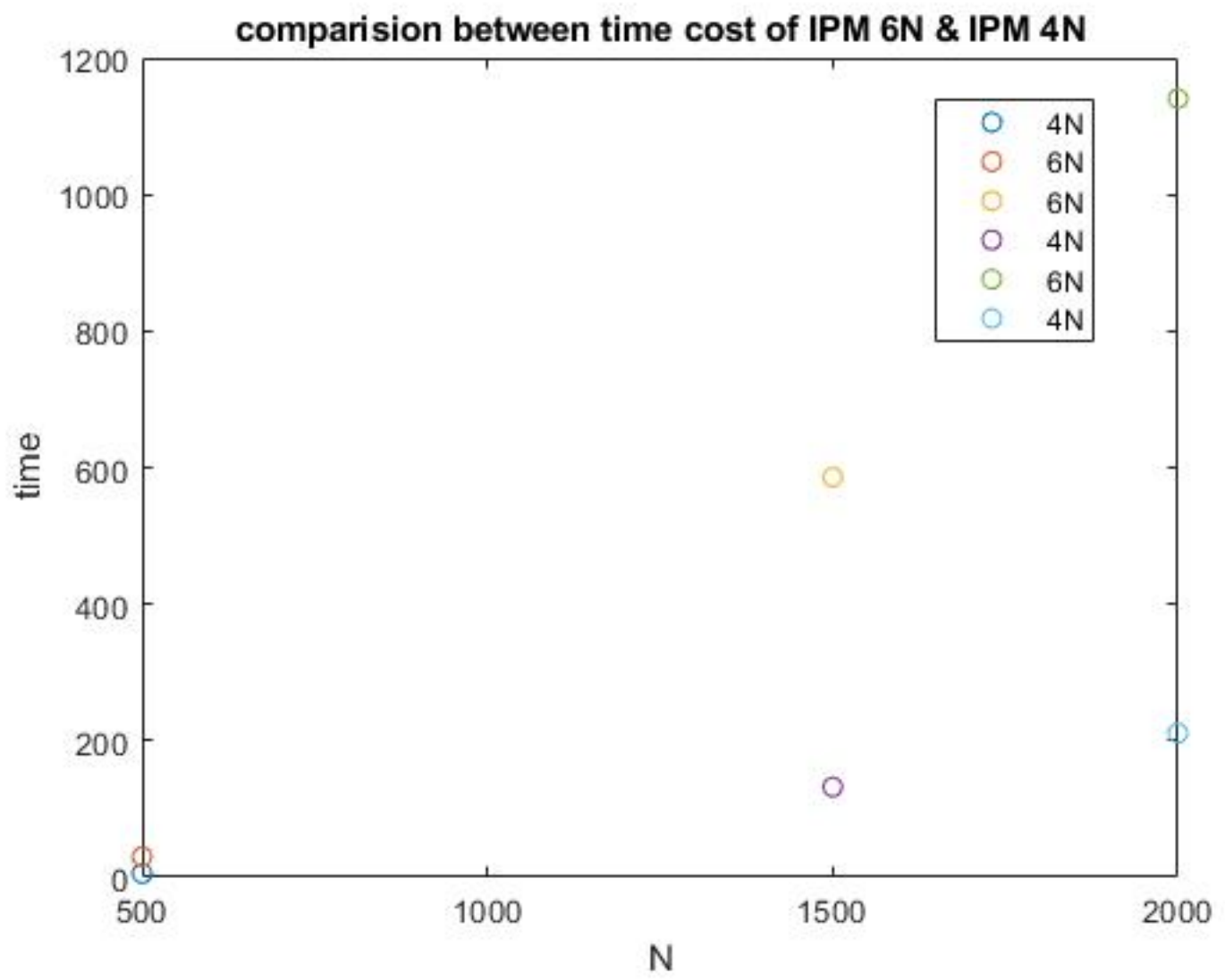

Figure 2.8: Time computation (s) comparison between the system 4N & 6N

*we can see that for the case 4N it take less time of computation than the case 6N.*

## 2.2 "Non-parametric IPM"

We put the slack variable just on the complementary condition, we are thus led to consider.

$$\begin{cases} (Au-b)_{3I-2} = \mathbf{0} \\ (Au-b)_{3I-1} = \mathbf{0} \\ (Au-b)_{3I} - w_I = \mathbf{0} \\ u_{3I} \bullet w_I = \mu \mathbf{e} \end{cases} \tag{2.16}$$

for $1 \leqslant I \leqslant 3N$, add the following equation

$$\frac{1}{2}\sum_{I=1}^{N} min\left\{u_{3I};0\right\}^2 + \frac{1}{2}\sum_{I=1}^{N} min\left\{w_I;0\right\}^2 + \mu^2 + \epsilon\mu = 0 \tag{2.17}$$

for a small value of $\epsilon > 0$ (unchanged during the iterations). System (2.16) - (2.17) consists of 4N+1 unknowns and 4N+1 equations. An ordinary Newton method is then applied to solve it.

**Remark 2.5.** *If all $u_{3I}, w_I \geqslant 0$ then Jacobi matrix of the system* (2.16) *as follow:*

$$DF(\chi) = \begin{bmatrix} DF(X) & -e \\ \boldsymbol{0} & 2\mu + \epsilon \end{bmatrix} \in \mathbb{R}^{m+n+1,m+n+1} \tag{2.18}$$

*where $m = 3N$, $N$ is the number of cells, $DF(X) \in \mathbb{R}^{m+n,m+n}$ , $e = (0,0,...,0,-1,-1,...,-1)^T \in \mathbb{R}^{m+n,1}$, first $m$ components is 0, components from m+1 to m+n is -1. $\boldsymbol{0} = (0,0,0,...,0) \in \mathbb{R}^{1,m+n}$*

**Proposition 2.1.** *Suppose we have the systems*

$$\begin{cases} F(x) = \mathbf{0} \\ x_i s_i - \mu = \mathbf{0}, \quad \forall i \in \{1, ..., N\} \\ f(x) + f(s) + \mu^2 + \epsilon\mu = \mathbf{0} \end{cases} \tag{2.19}$$

*where*

$$f(t) = \frac{1}{2}\sum_{i=1}^{N}(min(t_i, 0))^2, \quad \forall t \in \mathbb{R}^N$$

*then*

*(i) at any solution of the NPIPM system, we have*

$$x_i > -\frac{\epsilon}{\sqrt{2}}, \quad s_i > -\frac{\epsilon}{\sqrt{2}}, \quad \forall i \in \{1, ..., N\} \tag{2.20}$$

*(ii) if the NCP has at least a solution $(x^*, s^*)$ then $(x^*, s^*, 0)$ is a solution of the NPIPM system*

*(iii) If a solution admits negative component then the corresponding $\mu \in \left\{-\frac{\epsilon}{2} \pm \frac{\sqrt{\epsilon^2-4g}}{2}\right\}$*

*where $g = f(x) + f(s) \geqslant 0$ ,*

*Proof.* Suppose that $\exists i_0 / x_{i_0} < 0 \to f(x) > 0 \to g > 0$

- if $s_{i_0} = 0, \to \mu = 0 \to f(x) = f(s) = 0$ a contradiction.
- if $s_{i_0} < 0 \to \mu > 0$ contradict the fact that $f(x) + f(s) + \mu^2 + \epsilon\mu = 0$

So $s_{i_0} > 0$

Consider the equation:

$$\mu^2 + \epsilon\mu + g = 0 \tag{2.21}$$

we get $\Delta = \epsilon^2 - 4g \geqslant 0$ (because we suppose that (2.19) has a real solution), If $\Delta = \epsilon^2 - 4g < 0$, then we get complex solution of $\mu$ ( a contradiction )!.

Or we have

$$\epsilon^2 \geq 4g = 4(f(x) + f(s)) \geq 4f(x) = 2\sum_{i=1}^{N}(min(x_i, 0))^2 \tag{2.22}$$

Or

$$\epsilon \geqslant \sqrt{2}\sqrt{\sum_{i=1}^{N}(min(x_i, 0))^2} \tag{2.23}$$

We have that

$$\sqrt{\sum_{i=1}^{N}(min(x_i, 0))^2} + min(x, 0) \geqslant 0 \tag{2.24}$$

So that

$$\epsilon \geqslant -\sqrt{2}min(x, 0) \tag{2.25}$$

Or

$$x_i > min(x, 0) \geqslant -\frac{\epsilon}{\sqrt{2}}, \forall i \in\in \{1, ..., N\} \tag{2.26}$$

which complete the proof of (i)

Now we move to (ii), just replace $\mu = 0$ to the NPIPM system we got exactly the NCP system

Now we move to (iii), Consider the equation

$$\mu^2 + \epsilon\mu + g = 0 \tag{2.27}$$

by exactly the same argument above, we have $\Delta = \epsilon^2 - 4g > 0$ we got the solution is

$$\mu = \frac{-\epsilon - \sqrt{\epsilon^2 - 4g}}{2} \quad Or \quad \mu = \frac{-\epsilon + \sqrt{\epsilon^2 - 4g}}{2} \tag{2.28}$$

which complete the proof of 3 □

**Remark 2.6.** *The proposition above is the theorical results to say that $\epsilon$ shouldn't be too large, If $\epsilon$ is small then we have the lower bound of $x_i$ is negative but not a large negative number which means we are negative but not far from 0.*

### 2.2.1 Numerical results of Non-parametric IPM

Before going to convergence, we make the sensitivity analysis of $\epsilon \downarrow 0$ the have a good "$\epsilon$" for non-parametric IPM.

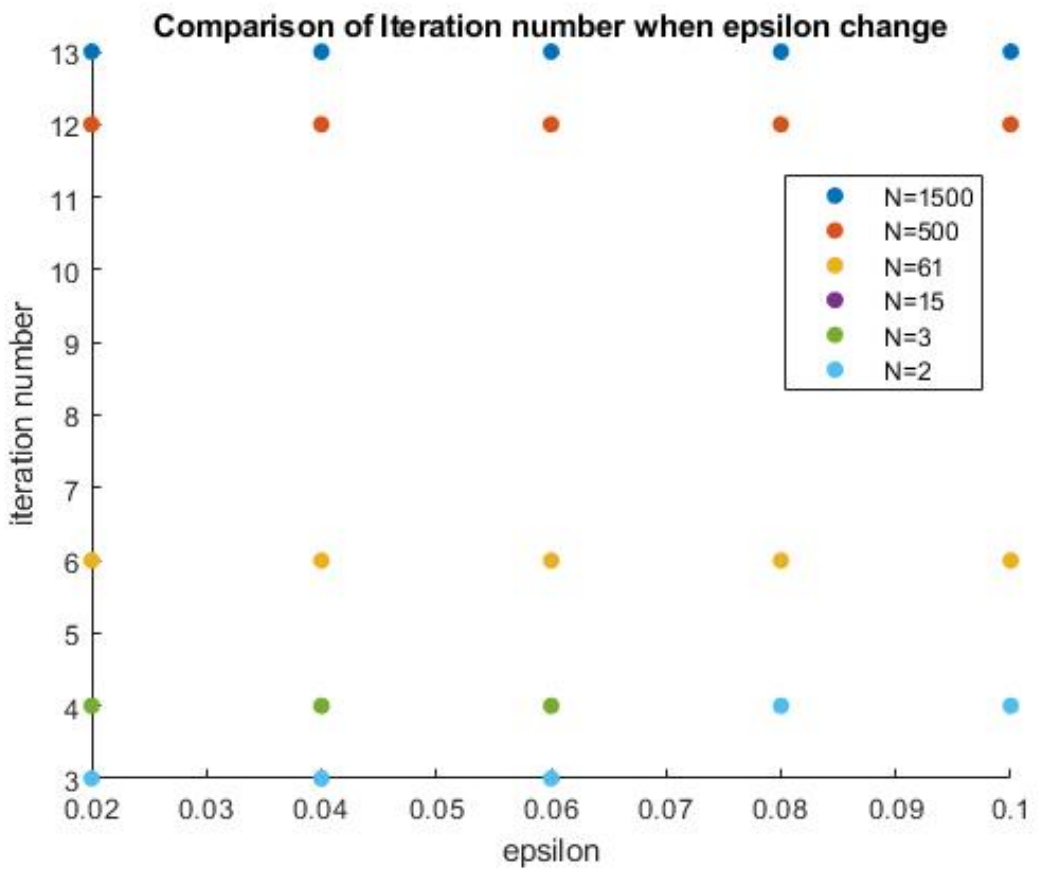

(a) the first choice of the initial point

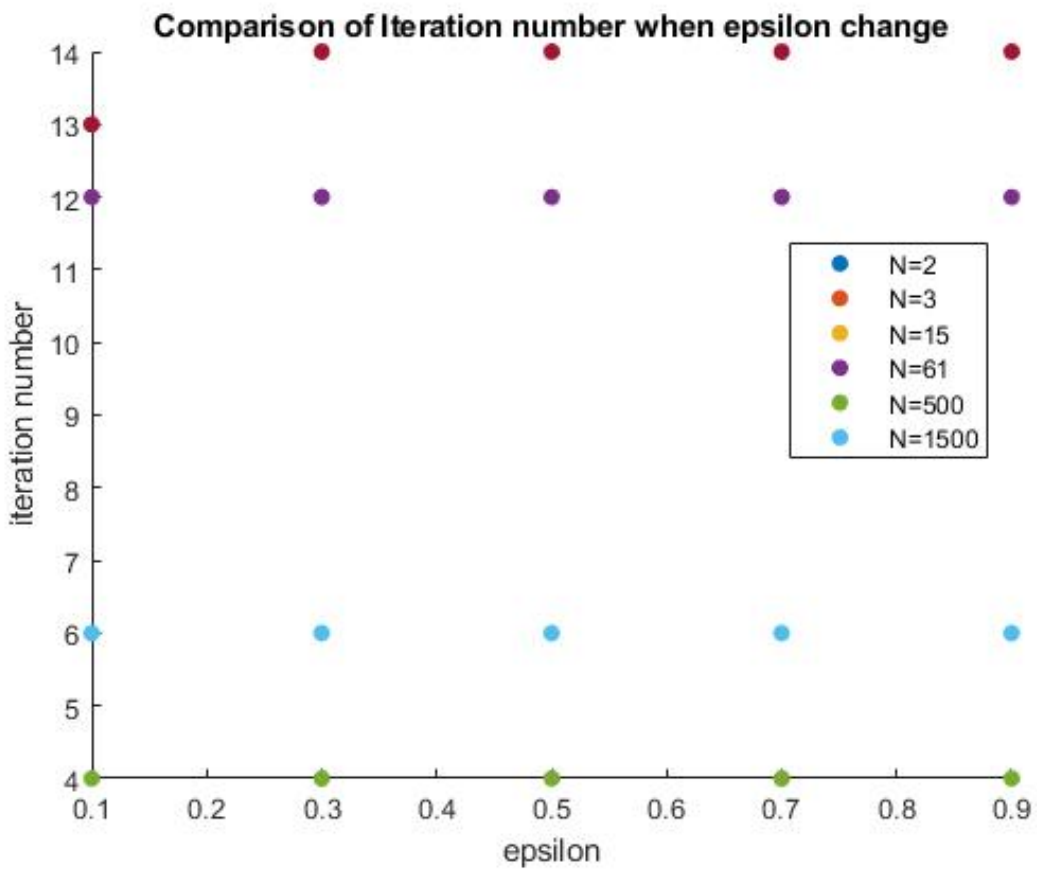

(b) the first choice of the initial point

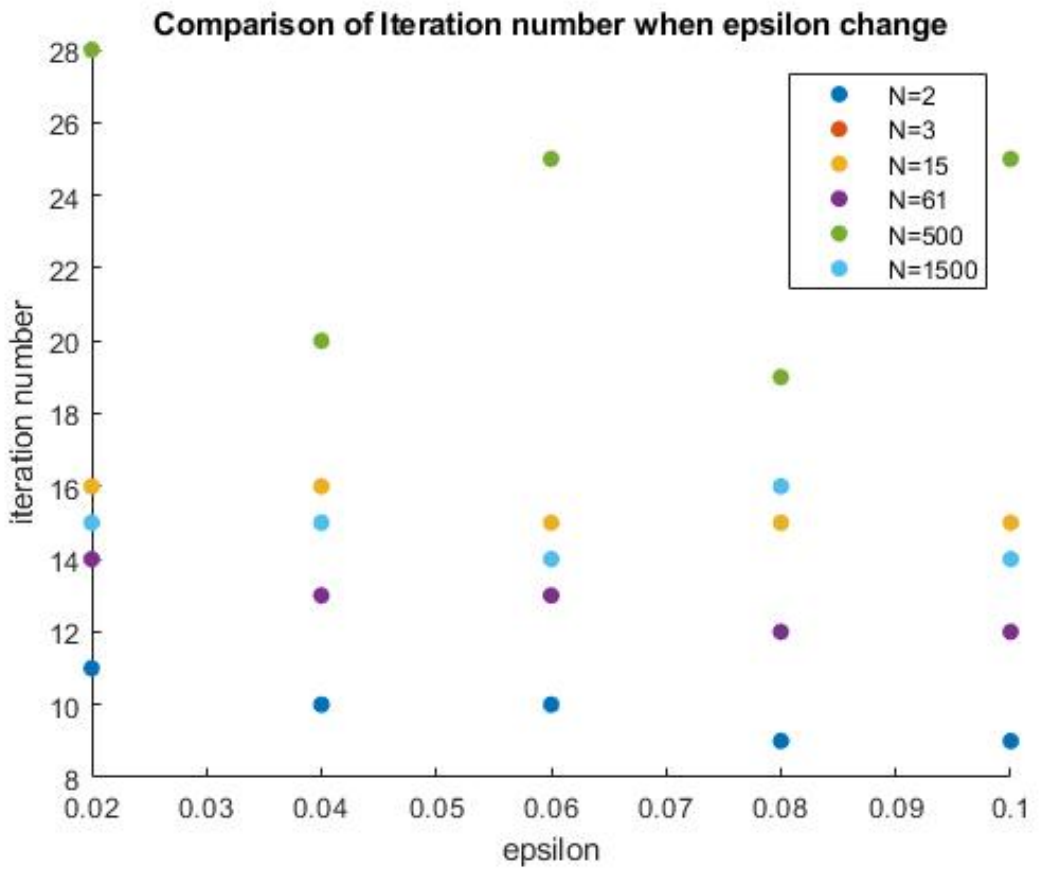

(c) For the second choice of the initial point

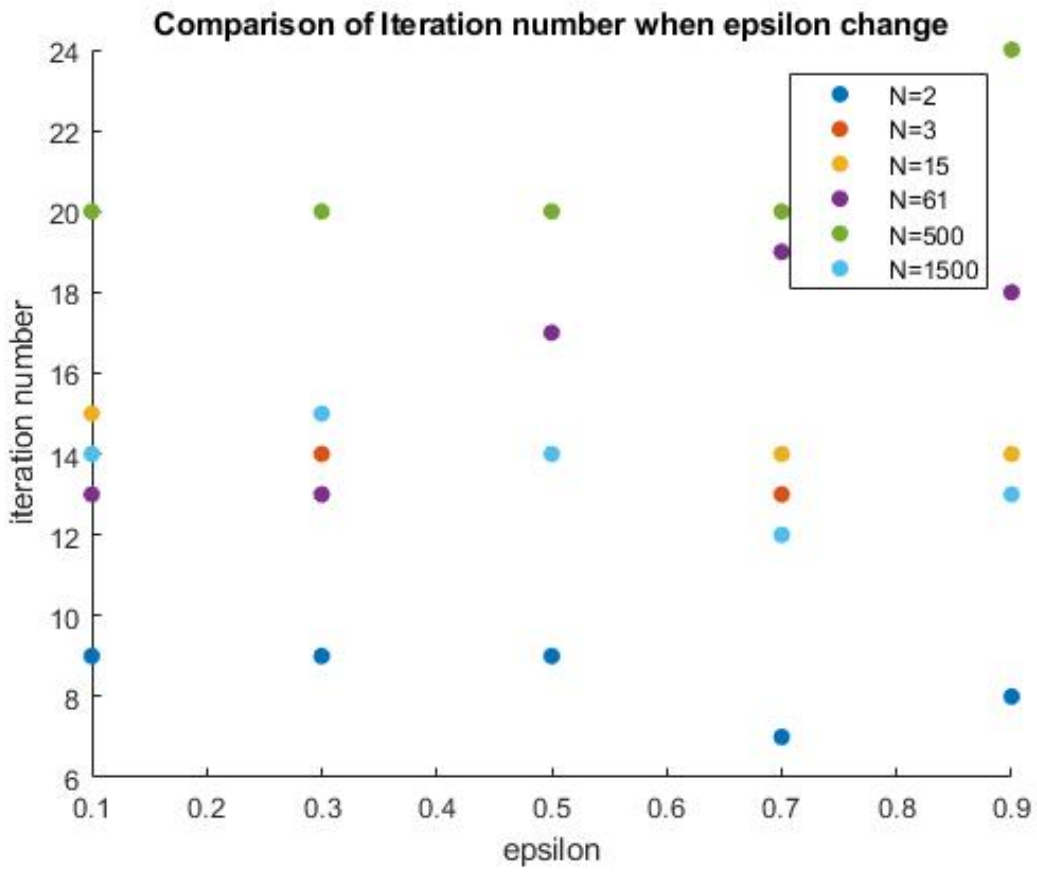

(d) For the second choice of the initial point

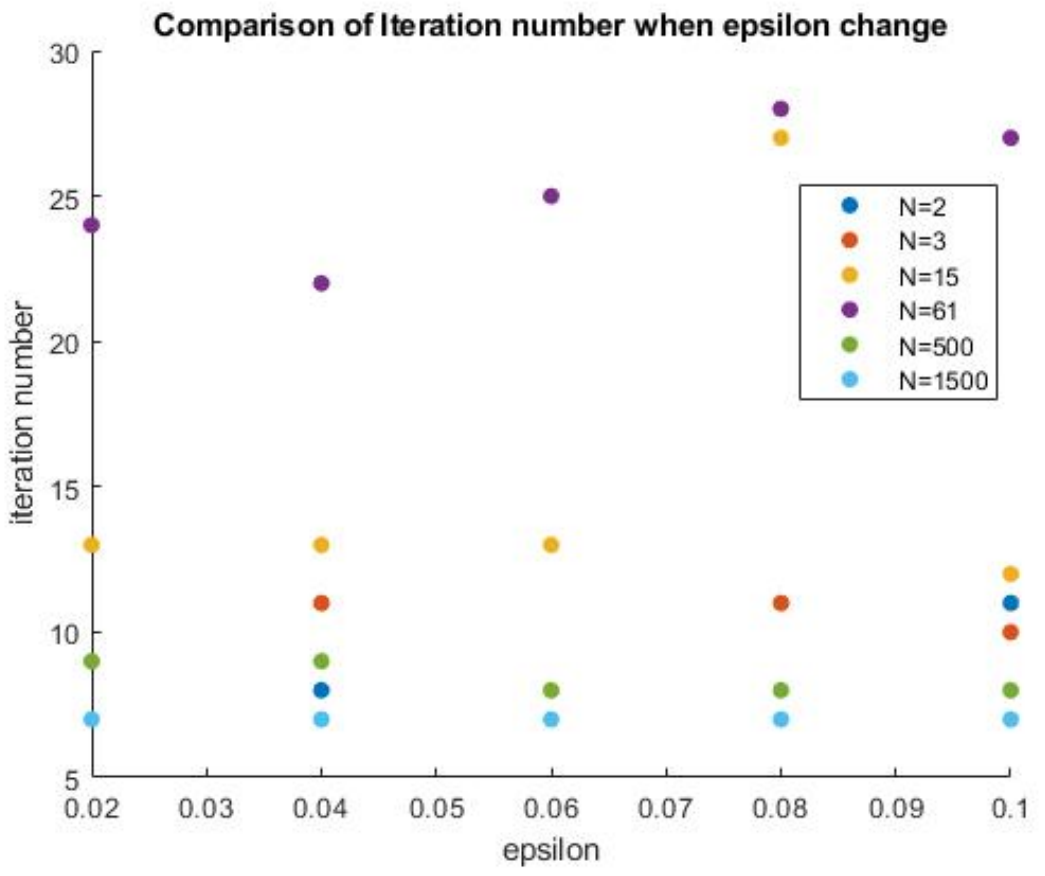

(e) For the third choice of the initial point

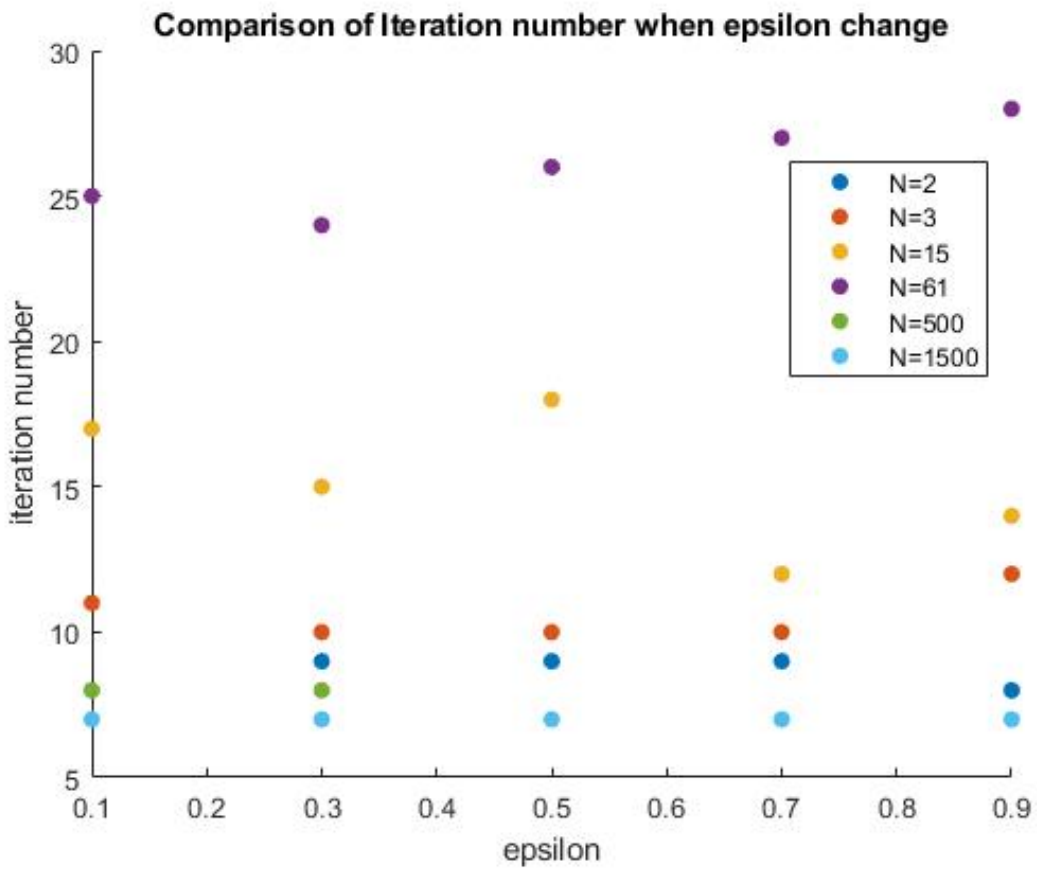

(f) For the third choice of the initial point

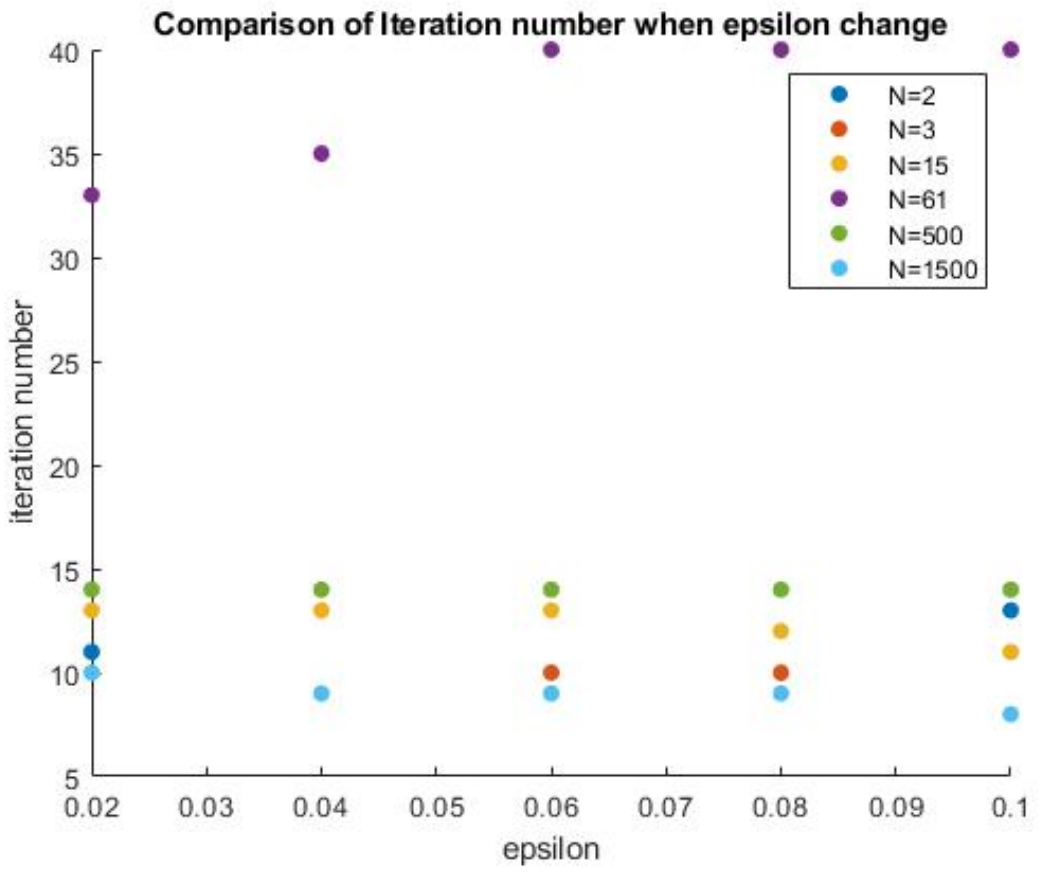

(g) Random of the initial point

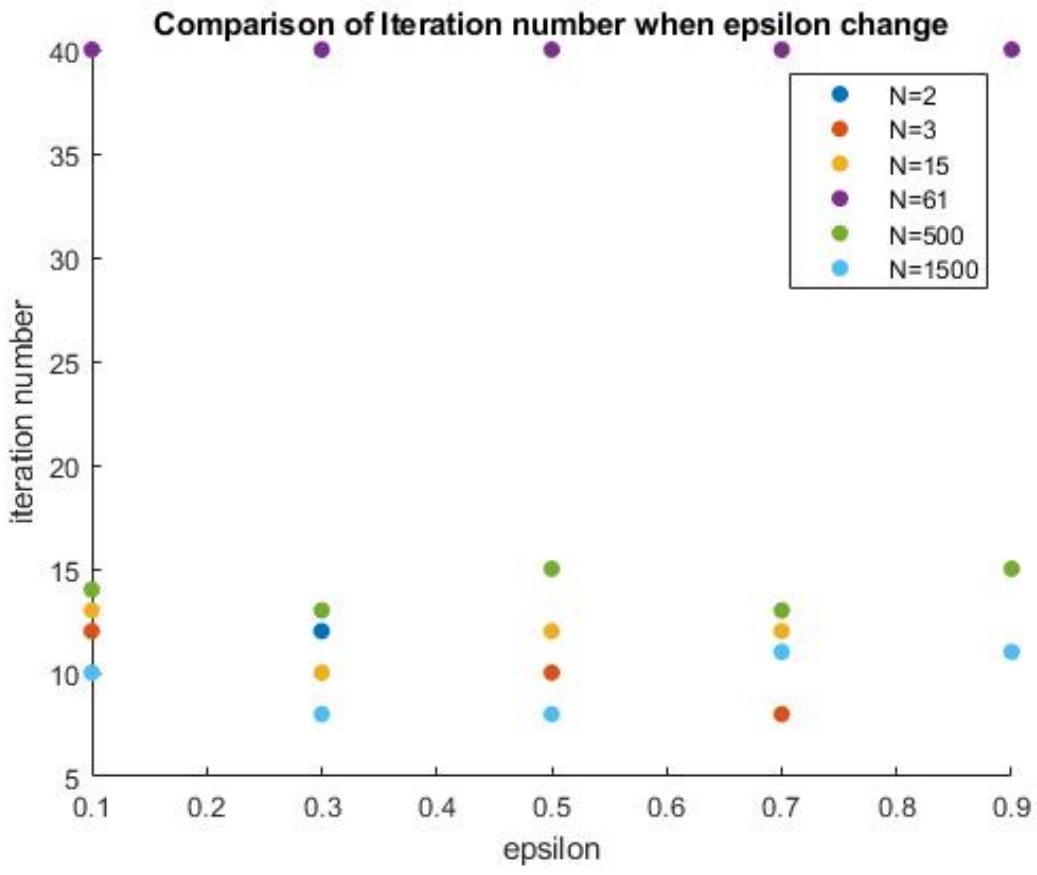

(h) Random of the initial point

we can conclude that depends on cases, we choose a suitable initial point and $\epsilon$. Especially for the cases N = 61, we can choose the first initial point and $\epsilon$ near 0. Other cases do not vary too much when $\epsilon$ change.Various option of $\epsilon$ is the same.

Below is the convergence of non-parametric IPM.

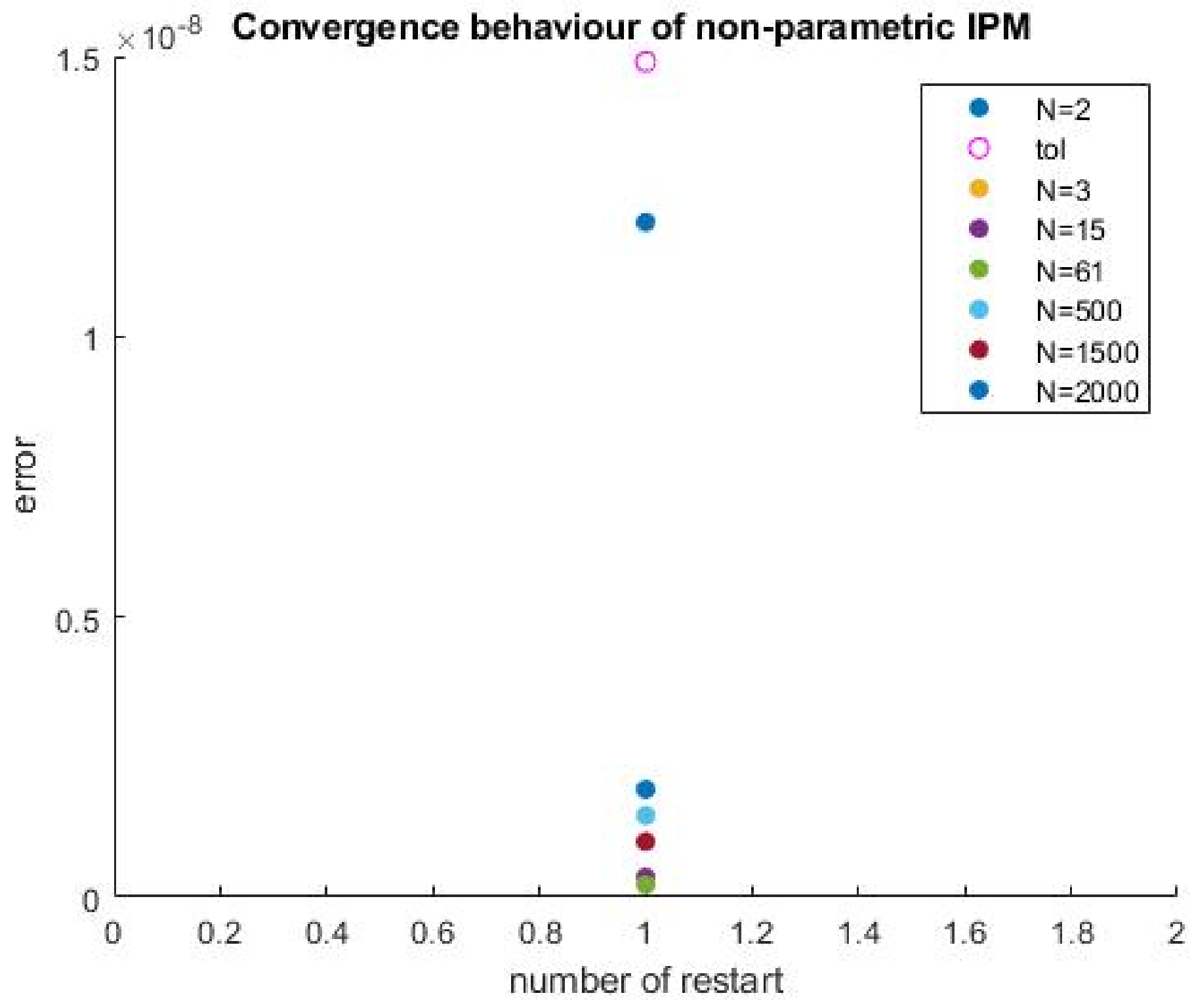

Figure 2.10: The convergence behavior of non-parametric IPM for the various option N

### 2.2.2 " Numerical resutls of Non-parametric IPM" with Armijo line search

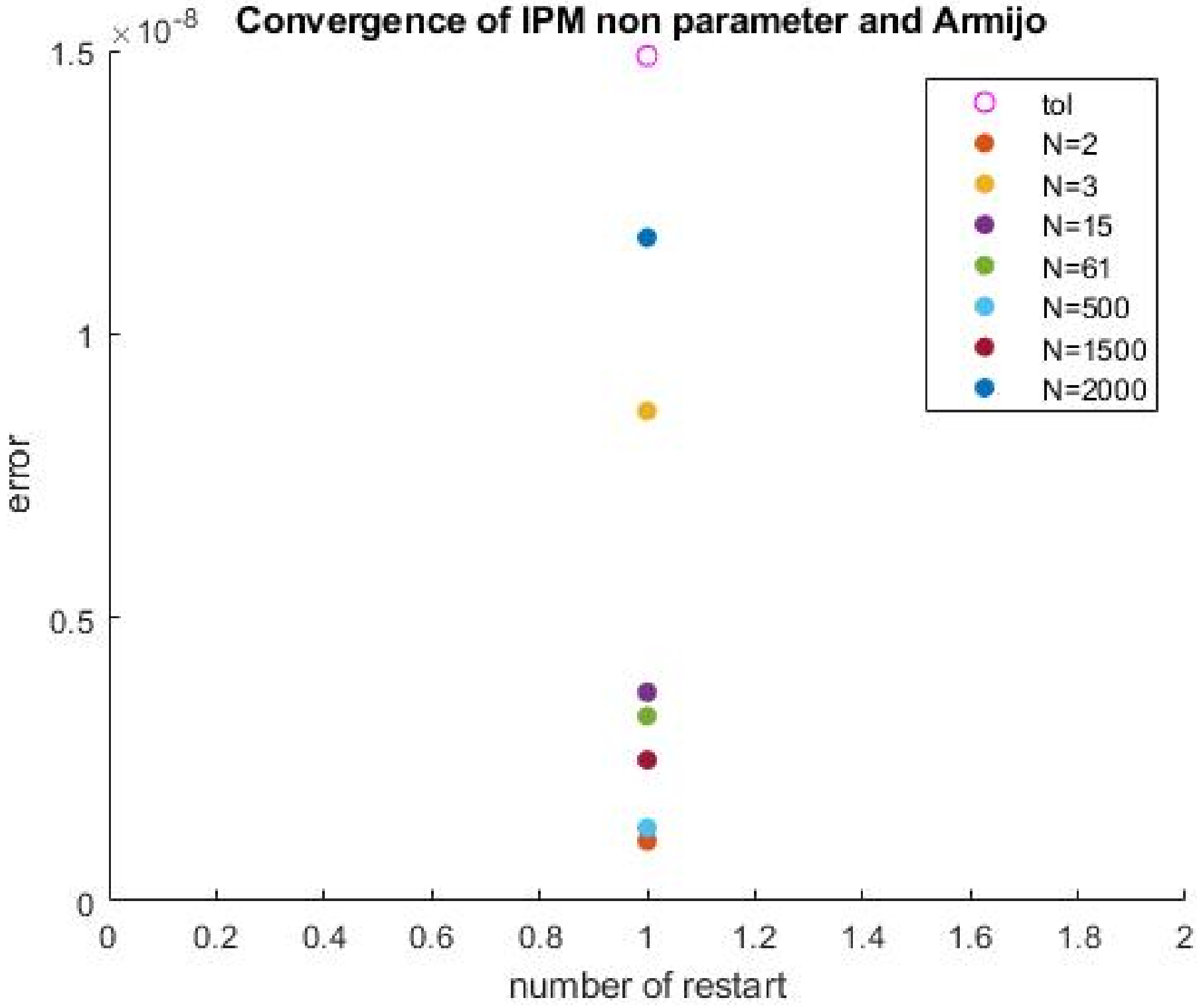

Figure 2.11: The convergence of non-parametric IPM with Armijo for the various option N

**Remark 2.7.** *One observation that it's not difference too much between* ***"IPM non parametric"*** *and* ***"IPM non parametric with Armijo line search"****. Define* $z = |\alpha_{max} - \rho|$*, with rho is the step size find in Armijo function with* $t_0 = \alpha_{max}$ *. Run our function with various N we can see that*

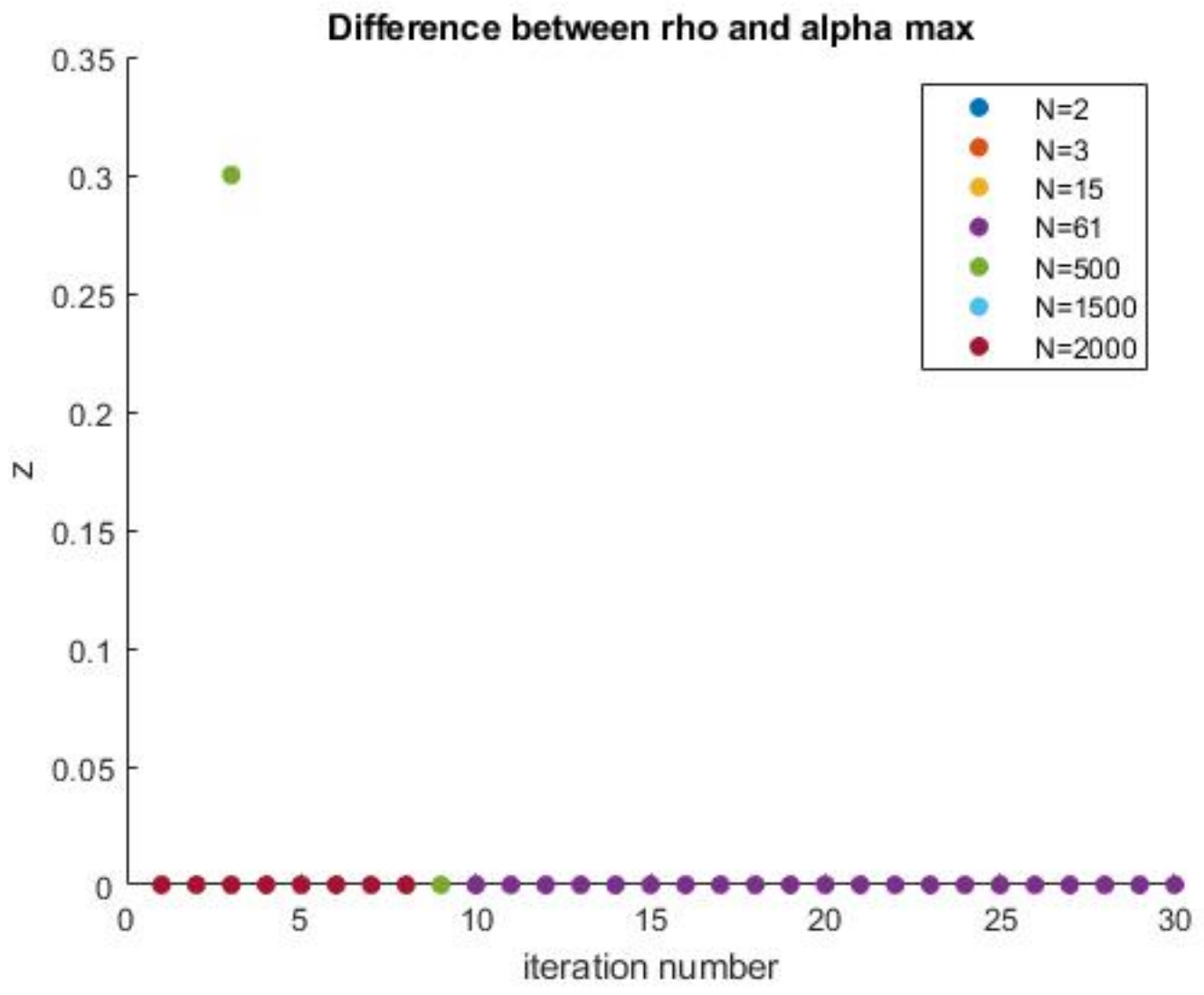

Figure 2.12: Difference between $\rho$ & $\alpha_{max}$ for various option N

***Now we can see the difference here.*** *Sometimes Armijo start with* $t_0$ *and return another value as can be seen in the case* $N = 500$

**Remark 2.8.** *The same as before ,we try to make non-parametric IPM still converge without computing* $\alpha_{max}$*, we change the merit function.*

$$\Phi(x) := \frac{1}{2}\|F(x)\|^2 + \frac{1}{2}\sum_{I=1}^{N} min\left\{u_{3I};0\right\}^2 + \frac{1}{2}\sum_{I=1}^{N} min\left\{w_I;0\right\}^2 \tag{2.29}$$

*The reason why we adding the term*

$$\frac{1}{2}\sum_{I=1}^{N} min\left\{u_{3I};0\right\}^2 + \frac{1}{2}\sum_{I=1}^{N} min\left\{w_{3I};0\right\}^2$$

*is*

*Consider the function*

$$f(y) = \frac{1}{2}\left(min\left\{y,0\right\}\right)^2 \tag{2.30}$$

*we have*

$$Df(y)(d) = \begin{cases} yd & if\ y \leqslant 0 \\ 0 & if\ y \geqslant 0 \end{cases} \tag{2.31}$$

*Then Armijo line search with the function f is*

$$f(z + td) \leqslant f(z) + mtDf(z)(d) \tag{2.32}$$

*where* $m \in \left(0, \frac{1}{2}\right)$ *, suppose* $x^{(k)} \geqslant 0$ *then* $f(x^{(k)}) = 0$ *,* $Df(x^{(k)})(d) = 0$*, apply to* $x^{(k)}$ *equation* (2.32) *we have,*

*Or*

$$0 \leqslant \frac{1}{2}\left(min\left\{x^{(k)} + td, 0\right\}\right)^2 \leqslant 0 \tag{2.33}$$

*Or*

$$\left(min\left\{x^{(k)}+td,0\right\}\right)^2=0 \tag{2.34}$$

*we get*

$$x^{(k+1)}=x^{(k)}+td\geqslant 0 \tag{2.35}$$

*So we have what we expected. The positivity of $x^{(k+1)}$. That's the idea how we maintain the positivity of $x$ when we don't compute the $\alpha_{max}$.*

*Simple computation leads to.*

$$D\Phi(x)(d)=(F(x),DF(x)(d))+\sum_{I=1}^{N}\mathbb{I}(u_{3I})d_i+\sum_{I=1}^{N}\mathbb{I}(w_I)d_i \tag{2.36}$$

*where*

$$\mathbb{I}(x)=\begin{cases}x & \text{if } x<0\\ 0 & \text{if } x\geqslant 0\end{cases} \tag{2.37}$$

*we still keep this condition on original Armijo line search*

$$\Phi\left(X^k+\beta^{j_k}d^k\right)\leq\Phi\left(X^k\right)+m\beta^{j_k}D\Phi\left(X^k\right)d^k \tag{2.38}$$

*Then run non-parametric IPM with the new merit function we obtain*

| **Cases** | **Convergence** | Lowest Iteration | Highest iteration |
|---|---|---|---|
| N = 2 | ✔ | 10 | 18 |
| N = 3 | ✔ | 6 | 8 |
| N = 15 | ✔ | 8 | 8 |
| N = 61 | ✔ | 22 | 32 |
| N = 500 | ✔ | 10 | 11 |
| N = 1500 | ✔ | 7 | 7 |
| N = 2000 | ✔ | 9 | 9 |

*when we compare with the original ways to compute $\alpha_{max}$ we obtain the following results.*

| **Cases** | **Old strategy** | Armijo |
|---|---|---|
| N = 2 | 14 | 13 |
| N = 3 | 7 | 22 |
| N = 15 | 8 | 22 |
| N = 61 | 27 | 13 |
| N = 500 | 9 | 11 |
| N = 1500 | 7 | 24 |
| N = 2000 | 8 | 9 |

***Condlude*** *We don't improve the results, the way above suggests to us another idea to compute $\alpha_{max}$.*

**Remark 2.9.** *New ways to compute* $\alpha_{max}$ *using Armijo line search. For example suppose* $d_i^{(k)} = -1$ *and* $u_i^{(k)} = 2$ *then apply Armijo line search with the function:*

$$f(x) = \frac{1}{2}\left(min\left\{x, 0\right\}\right)^2 \tag{2.39}$$

*we get* $\rho_i^{(k)} = 1.47$ *so if we move with this step size we have* $u_i^{(k+1)} = u_i^{(k)} + \rho_i^{(k)} d_i^{(k)} = 2 - 1.47 = 0.53 > 0$ *so we still have the positivity of* $u_i^{(k+1)}$*, work like this whenever* $d_j^{(k)} < 0$ *we have the* $\alpha_{max}$*. In practice, with this ways to compute* $\alpha_{max}$ *we don't need a safety parameter, we don't need to propose another form of* $\alpha_{max}$ *such as* $min\left\{1, \frac{x}{y}\right\} max\left\{1, \frac{y}{x}\right\} = 1$ *to avoid numeric problems, because the Armijo form* $t = \beta^{j_k}, \beta \in (0, 1)$ *where* $j_k$ *is a natural number, is ensure we don't touch the boundary, but sometimes* $\alpha_{max} > 1$ *so we need to scale it maximum is 1*

$$\alpha_{max} = min\left\{\alpha_{max}, 1\right\} \tag{2.40}$$

*Here we propose another idea is to implement the non-parametric IPM with Armijo line search. Fristly, we compute the* $\alpha_{max}$ *by proposing above. Then apply the original Armijo line search with function.*

$$\Phi(x) = \frac{1}{2}\left\|F(x)\right\|^2 \tag{2.41}$$

*But with this approach it's much more expensive because In the worst case we need 2N times run Armijo function to find a "good"* $\rho^k$ *for each 1,...,2N correspond to* $u_3, u_6, ..., u_{3N}, u_3, u_6, ..., u_{3N}, w_3, w_6, ..., w_{3N}$ *. In other words, we just the Armijo function to implement the non-parametric IPM .We call it with temporary name: "non-parametric IPM totally Armijo line search"*

***Numerical results of the Totally Armijo line search***

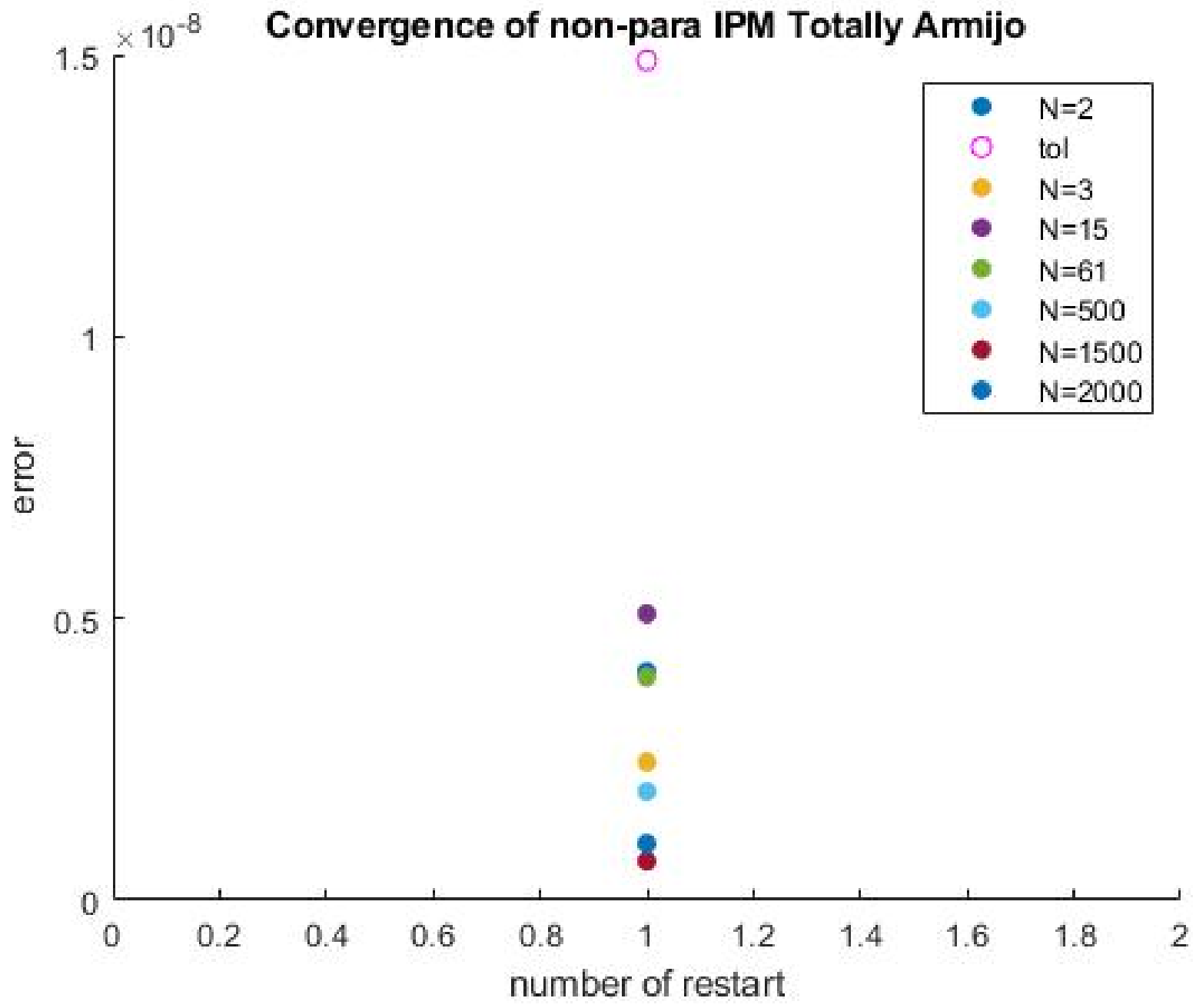

Figure 2.13: The convergence of non-parametric IPM Totally Armijo for the various option N

*From this, we have to compare 2 ways of computing to have a better ones.*

*(i) Convergence*

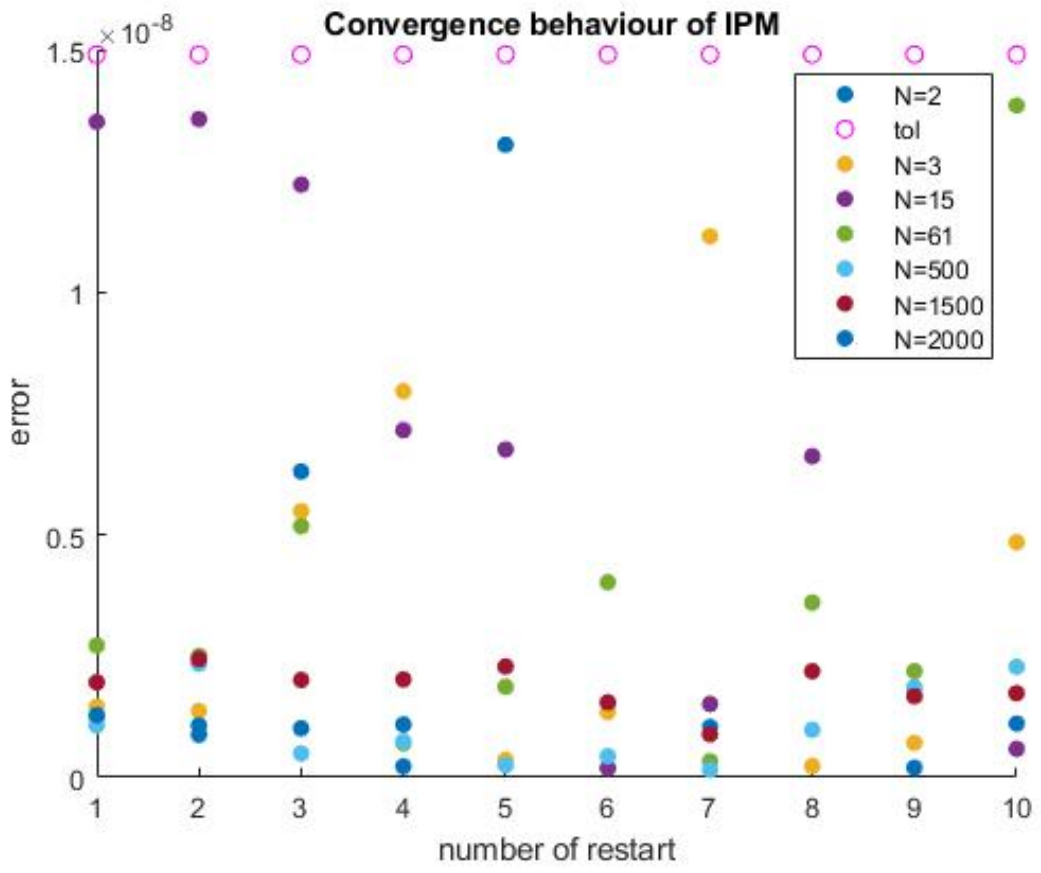

(a) The convergence of IPM with various N

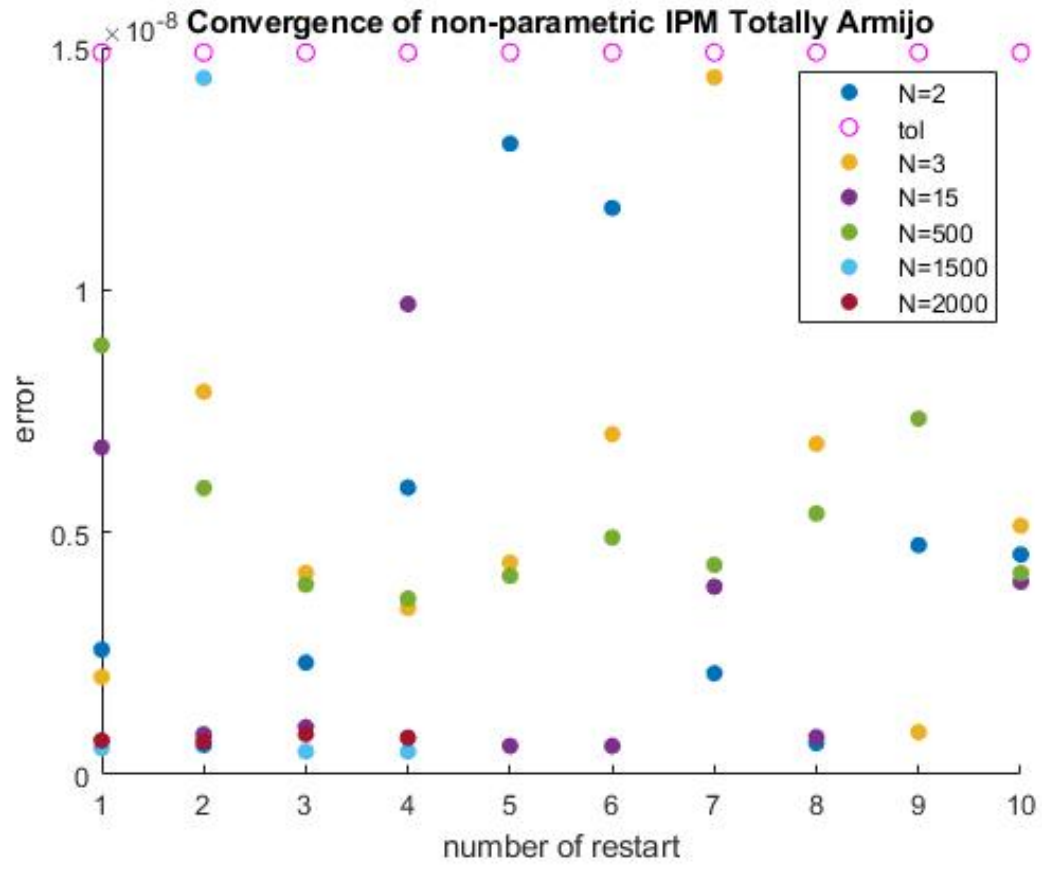

(b) The convergence of IPM totally Armijo with various N

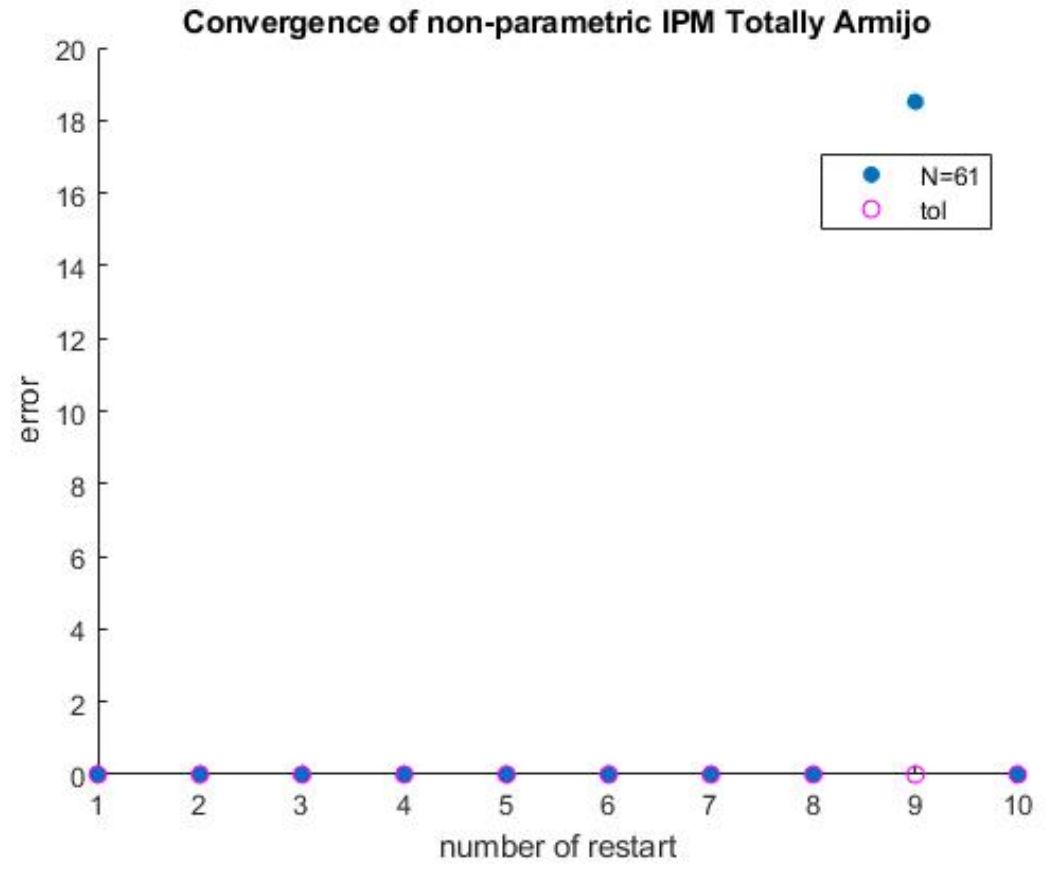

(c) The convergence of IPM totally Armijo with various N

Figure 2.14: Compare the convergence of IPM and IPM totally Armijo for various N

*(ii) Iteration*

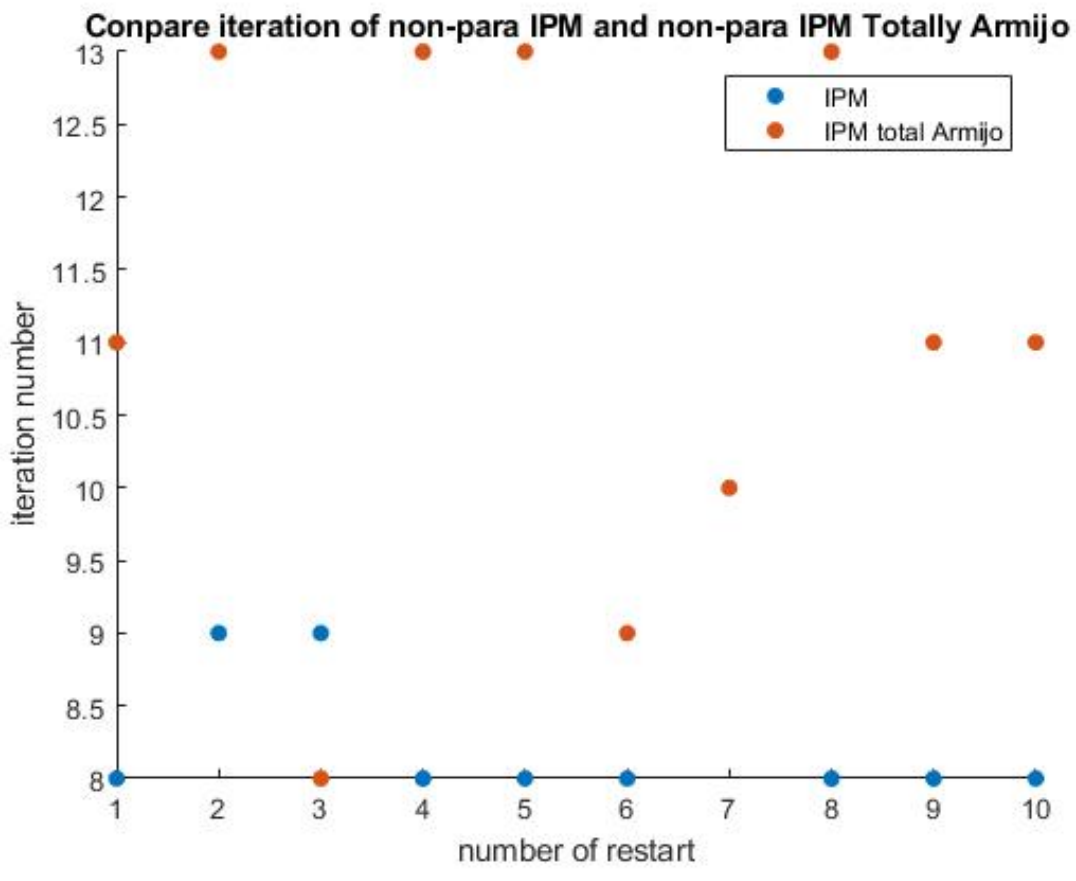

(a) $N = 2$

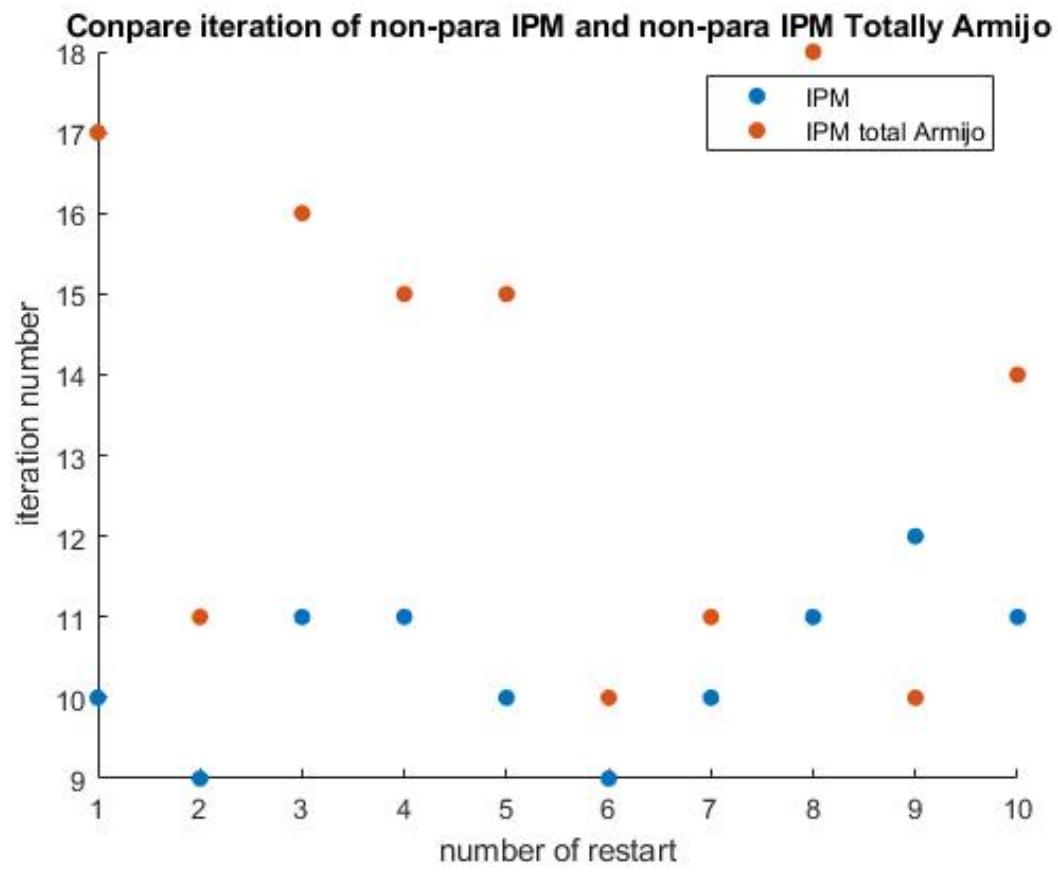

(b) $N = 3$

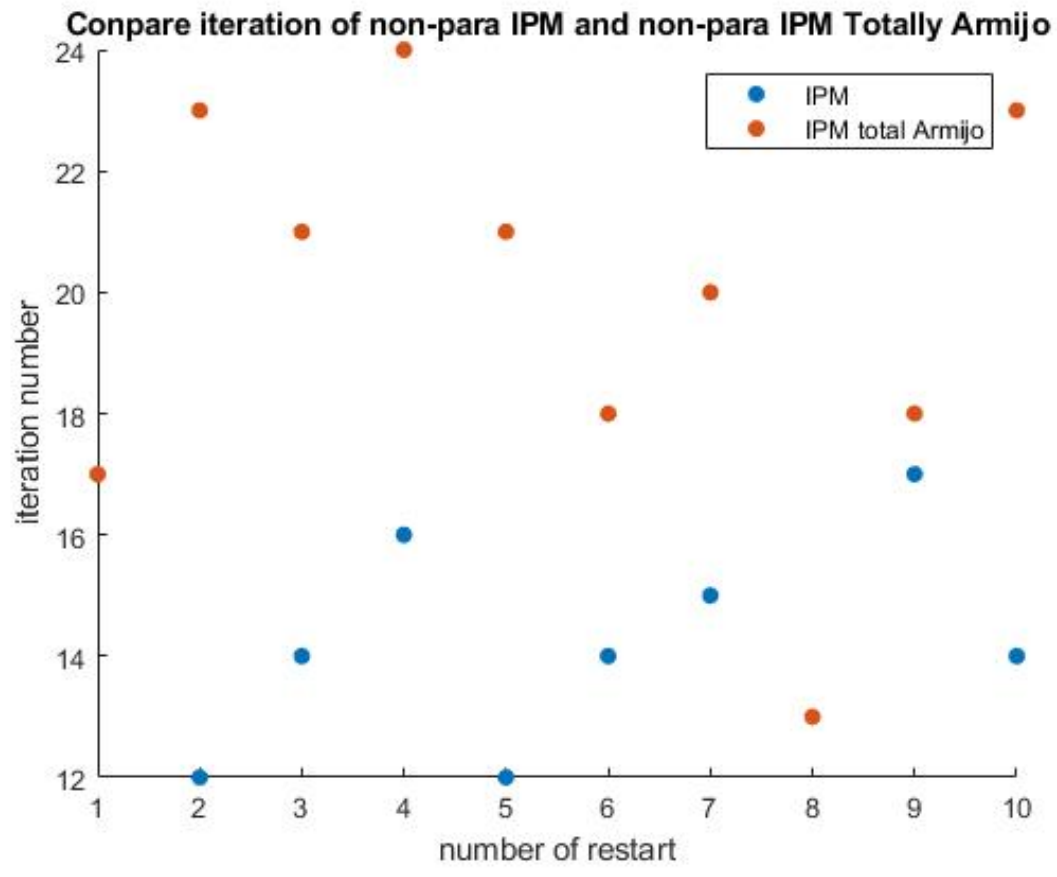

(c) $N = 15$

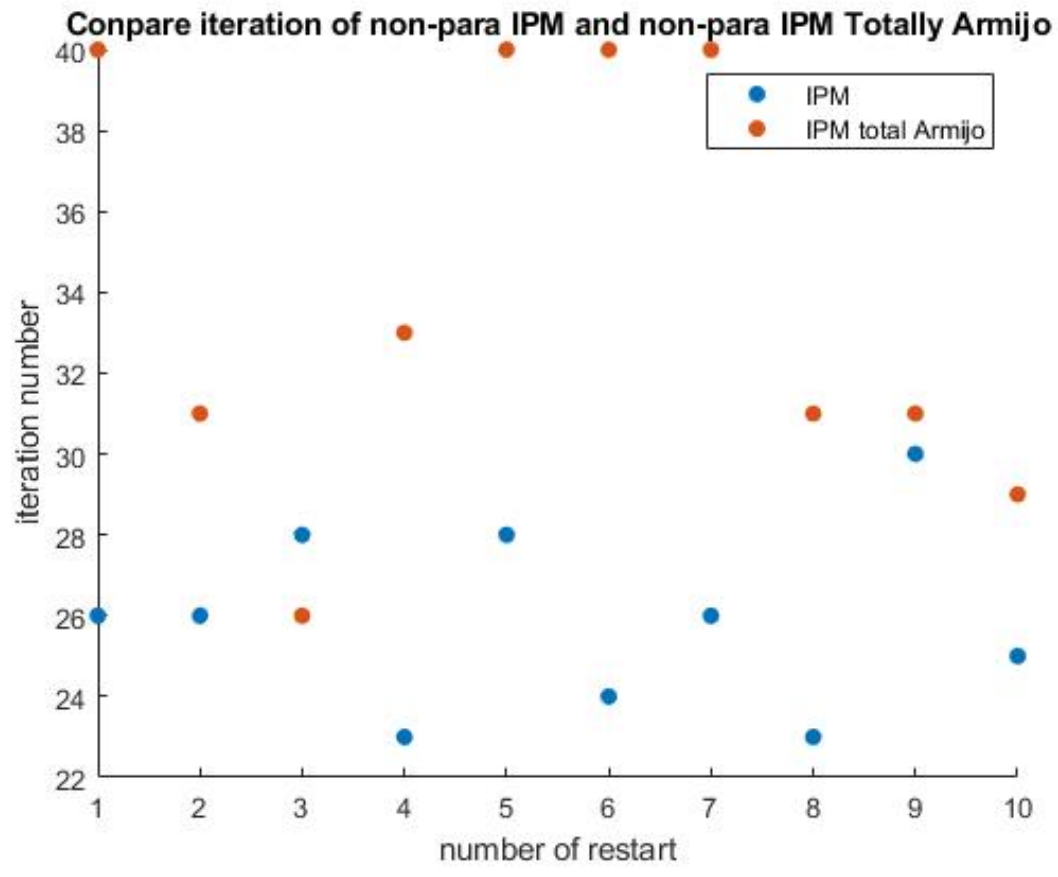

(d) $N = 61$

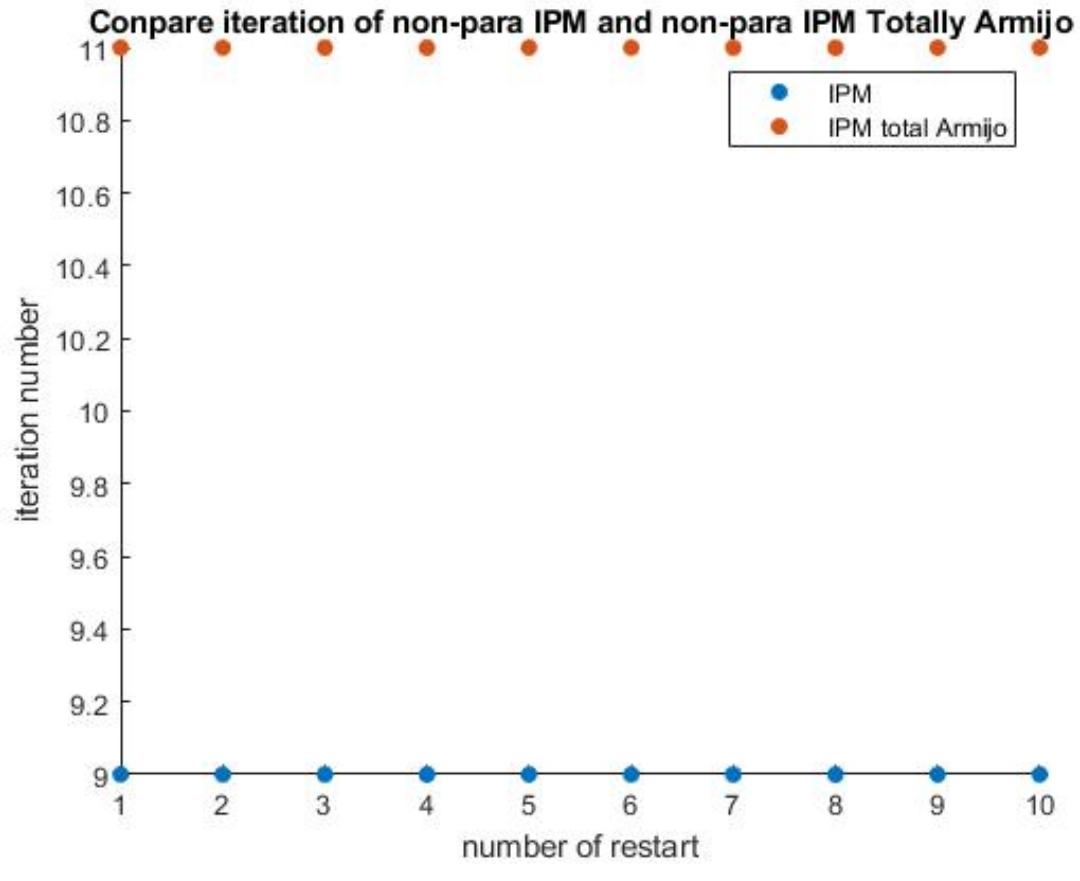

(e) $N = 500$

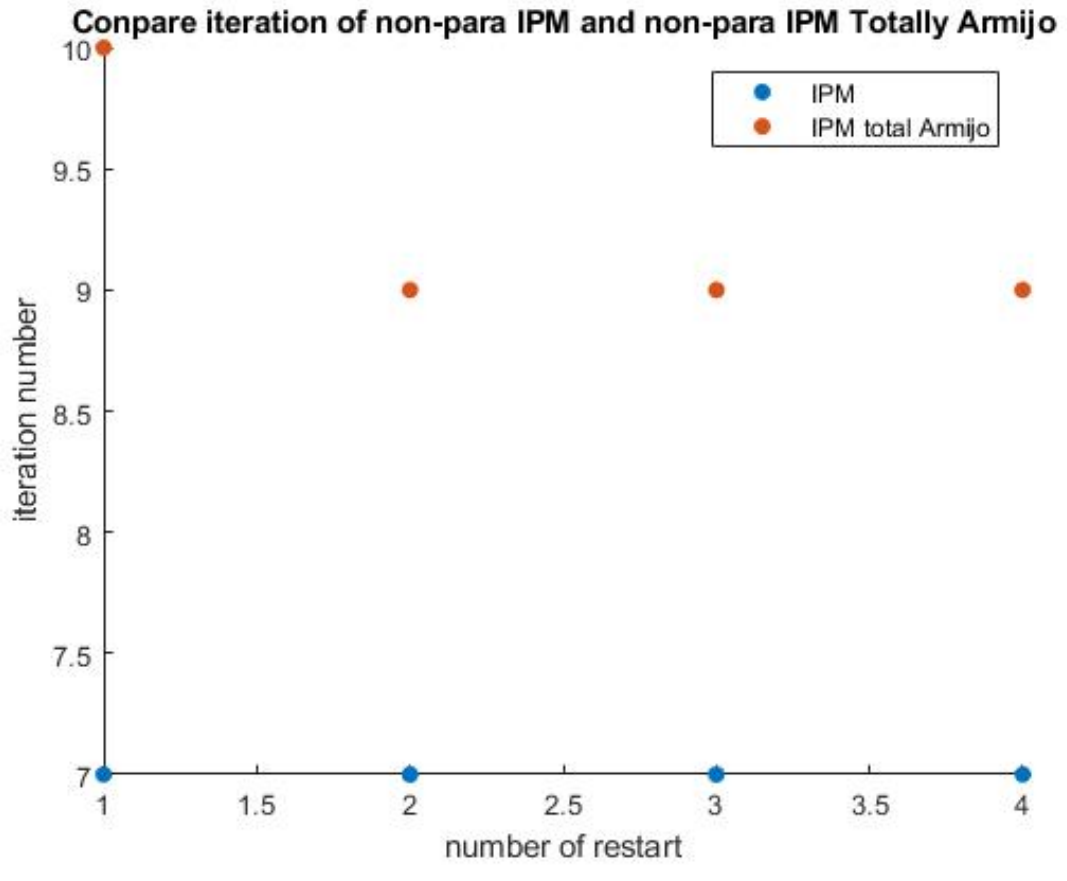

(f) $N = 1500$

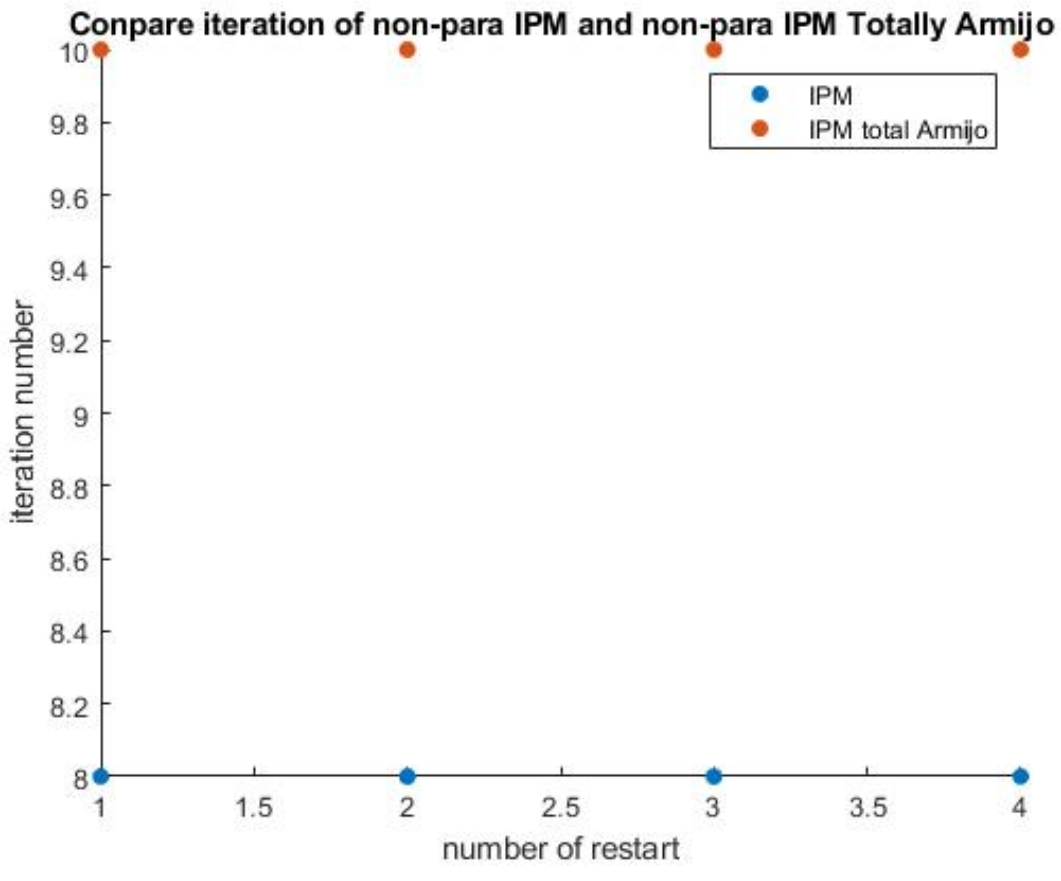

(g) $N = 2000$

*(iii) Time computing (s)*

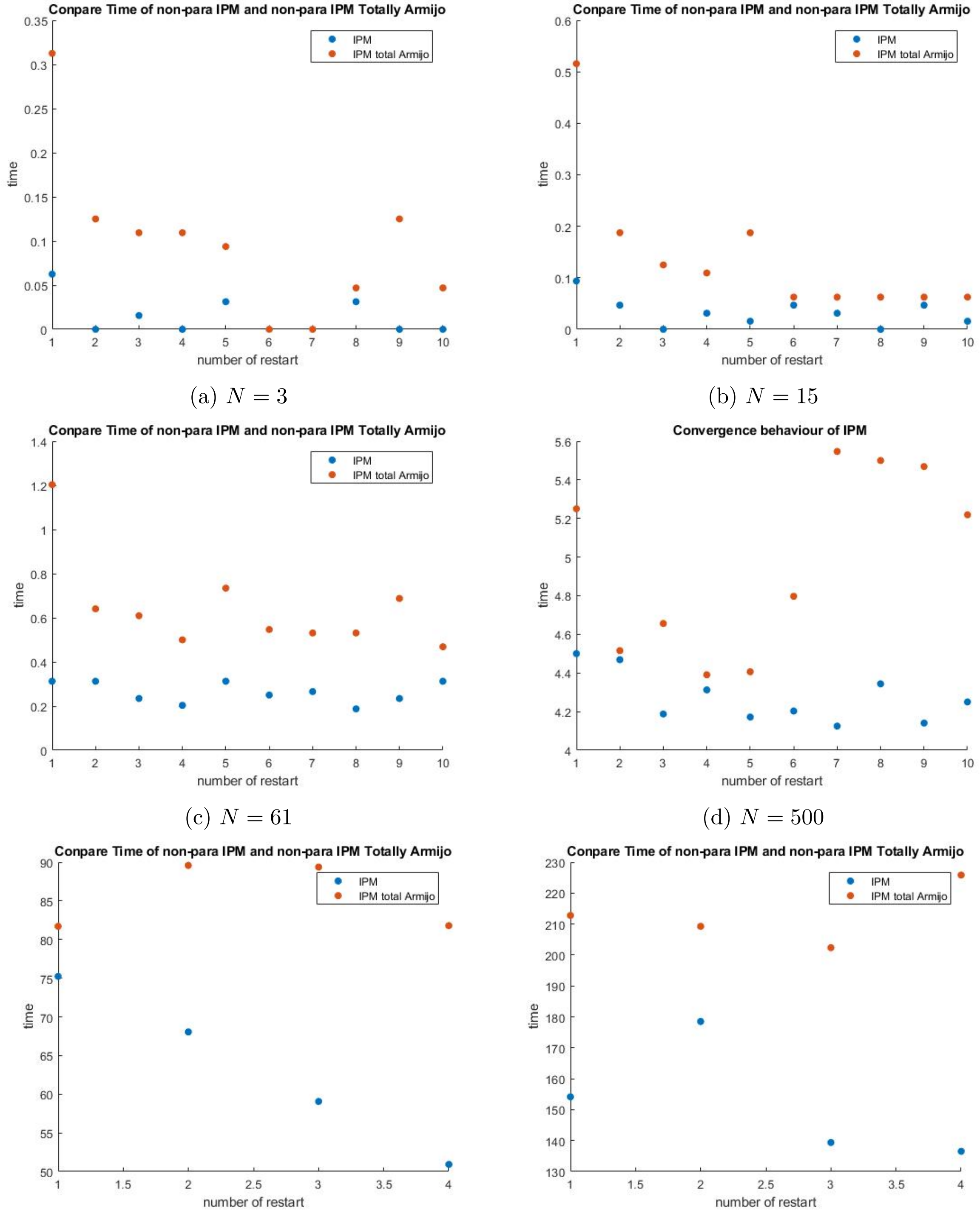

(a) $N = 3$

(b) $N = 15$

(c) $N = 61$

(d) $N = 500$

(e) $N = 1500$

(f) $N = 2000$

Figure 2.16: Time(s) comparison of non-parametric IPM and non-parametric IPM totally Armijo with various N

***Conclude: Some case N = 61 non-para IPM Totally Armijo does not converges***

***(In practice, in that case I can see that the $\alpha_{max}$ I found is too small and that's a reason why we can't decrease the value of error, it's seem like we don't move), Iteration and Time computing is higher than non-para IPM (because we can see that the original way we compute $\alpha_{max}$ is more simple than using non-parametric IPM Totally Armijo, with Totally Armijo when the size of $d < 0$ is large we have to compute large number of Armijo function to obtained $\alpha_{max}$). After all numerical results, we can see that the IPM Totally Armijo is worst than non-para IPM. So we dont' choose non-para IPM with $\alpha_{max}$ computed by Totally Armijo.***

# Chapter 3

# Compare the behavior of Standard methods, IPM with parameter and non-parametric IPM

## 3.1 Convergence

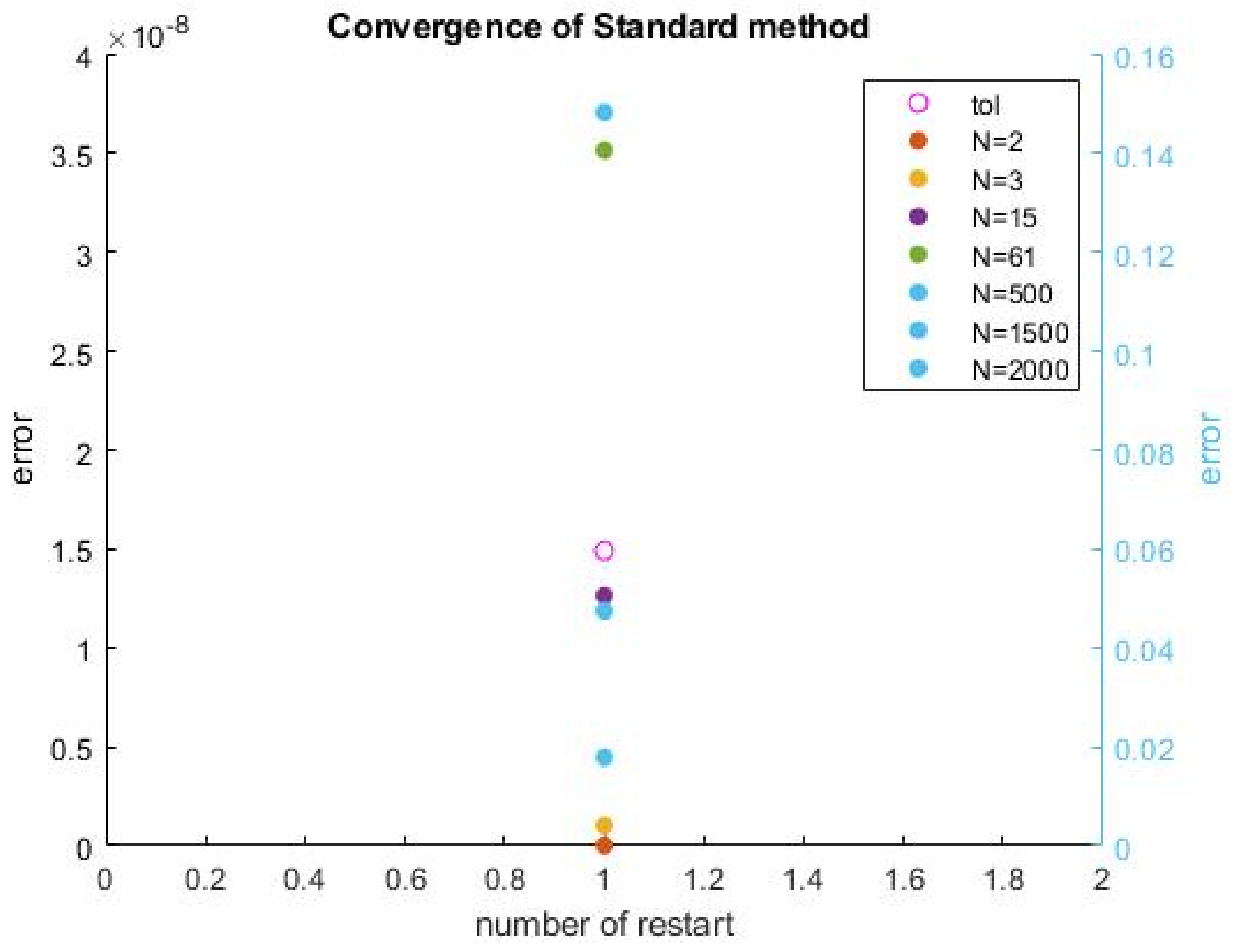

Figure 3.1: The convergence of Standard methods for the various option N

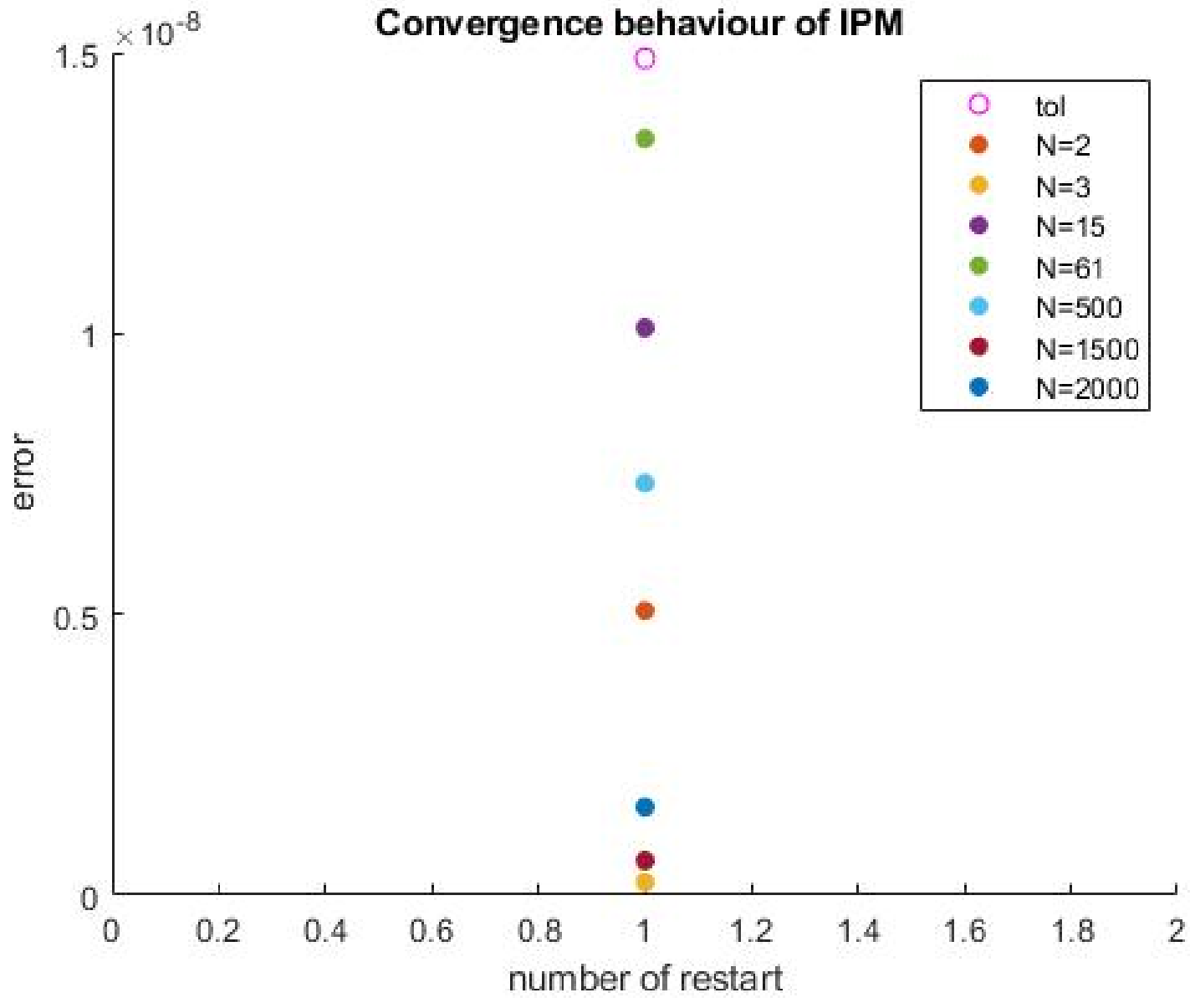

Figure 3.2: The convergence of IPM for the various option N

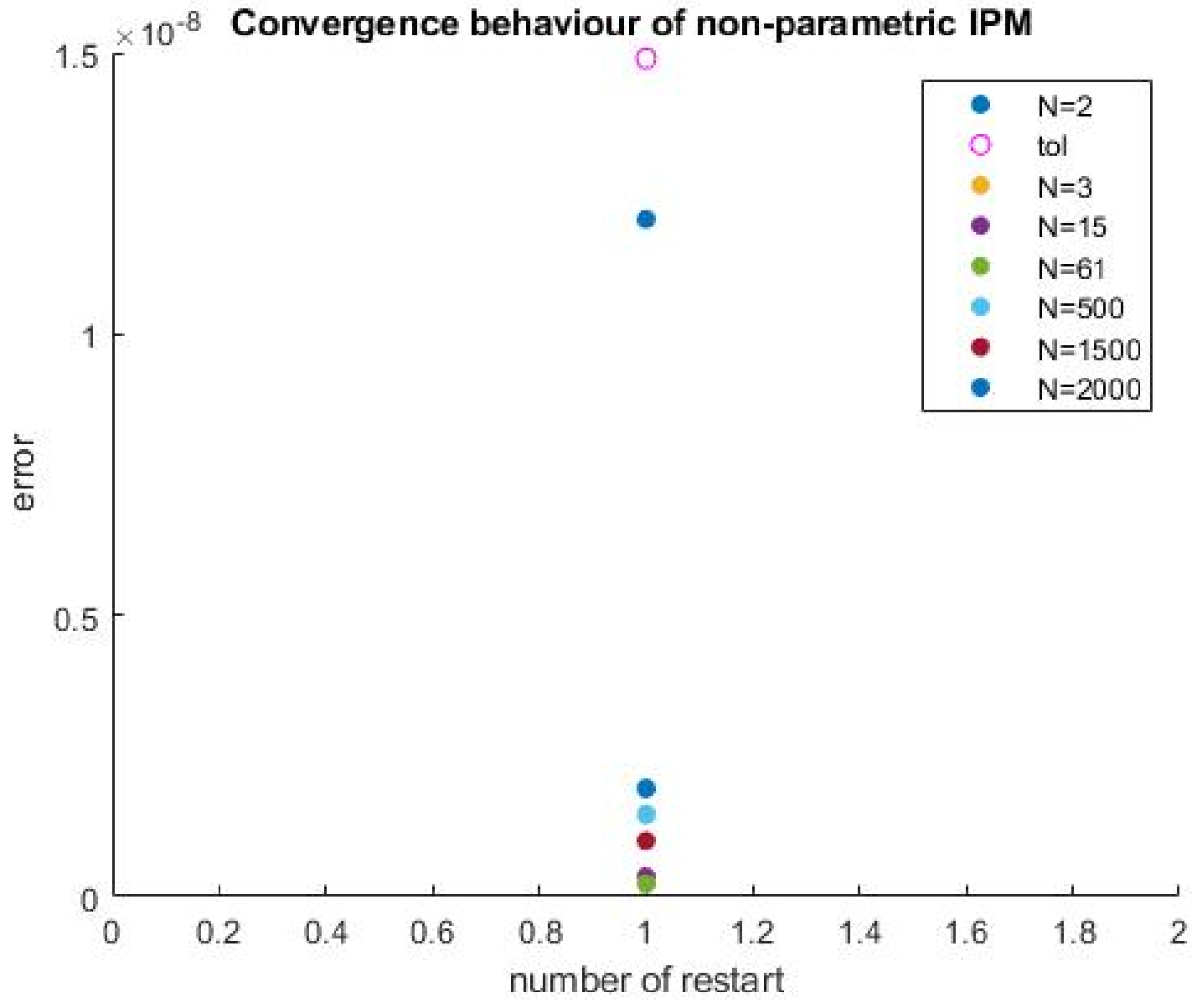

Figure 3.3: The convergence of non-parametric IPM for the various option N

## 3.2 Iteration

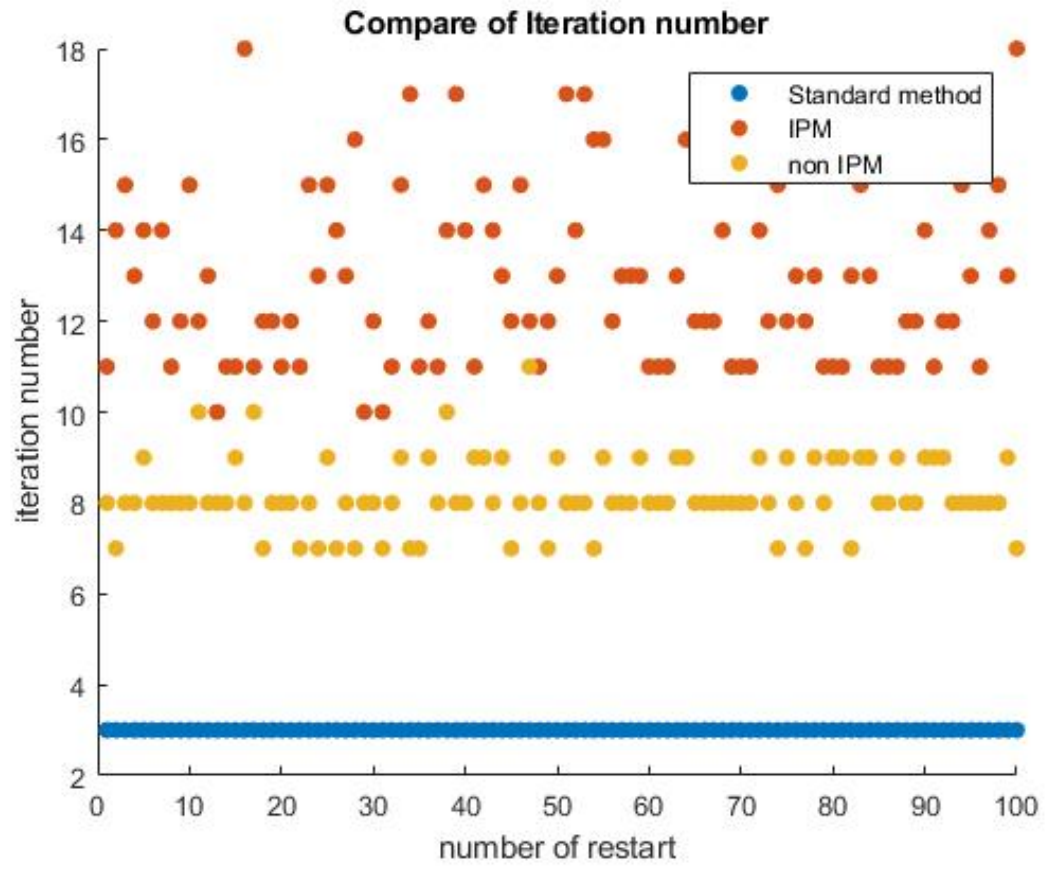

(a) $N = 2$

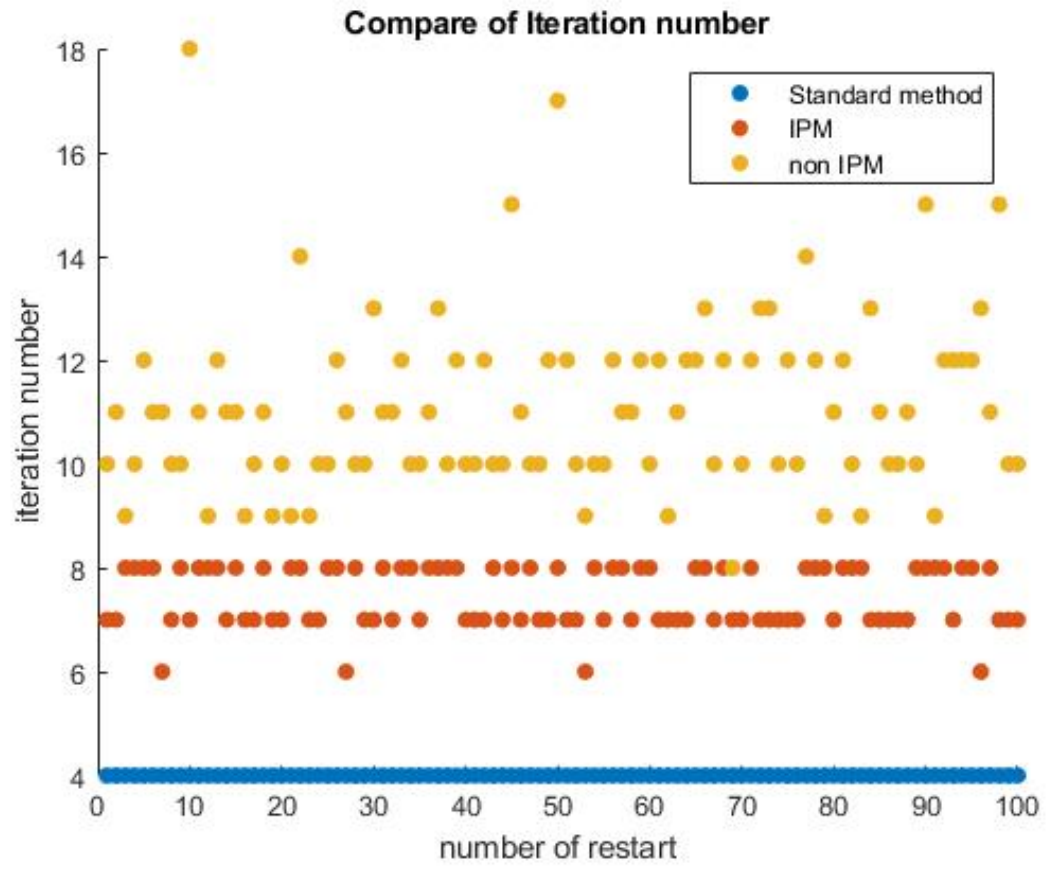

(b) $N = 3$

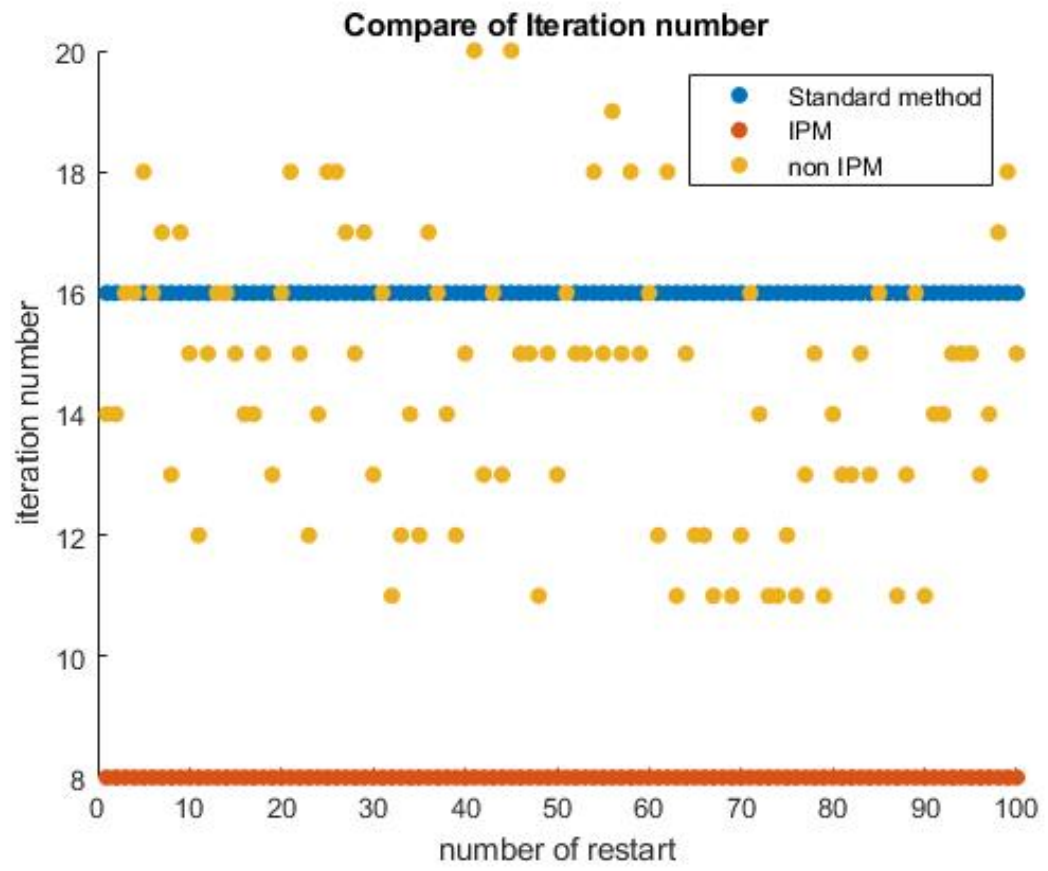

(c) $N = 15$

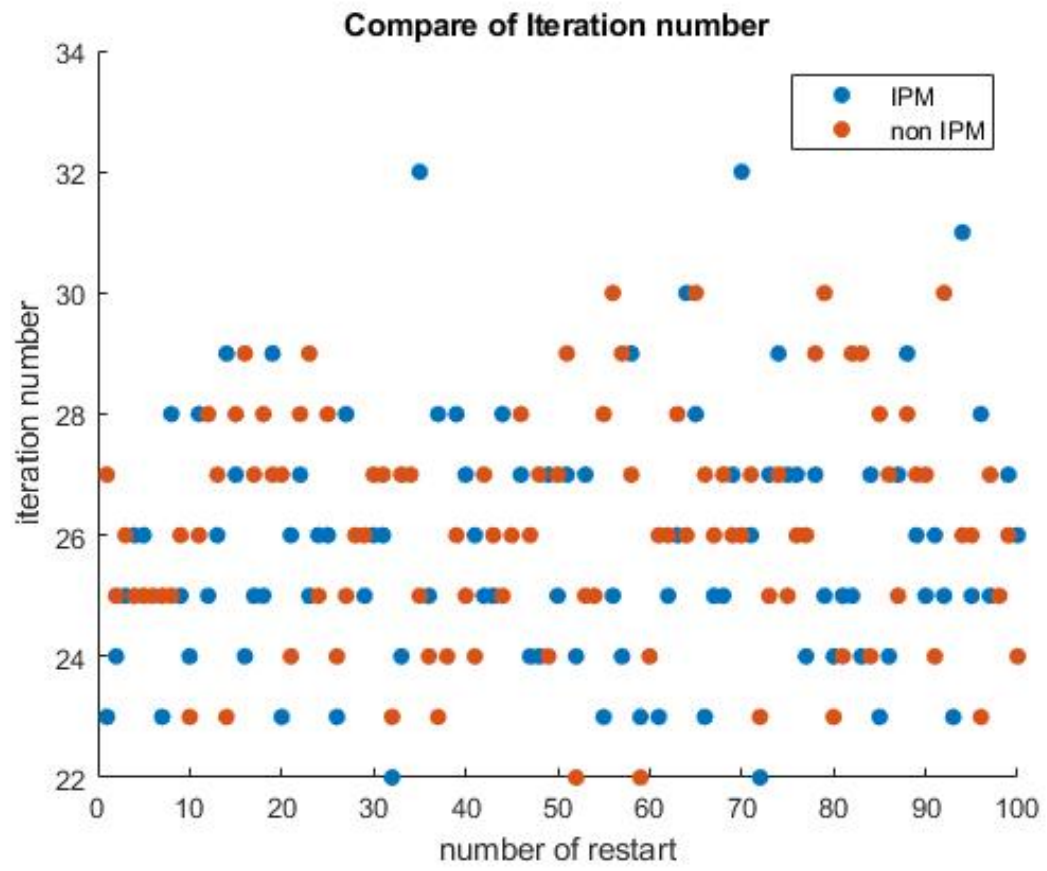

(d) $N = 61$

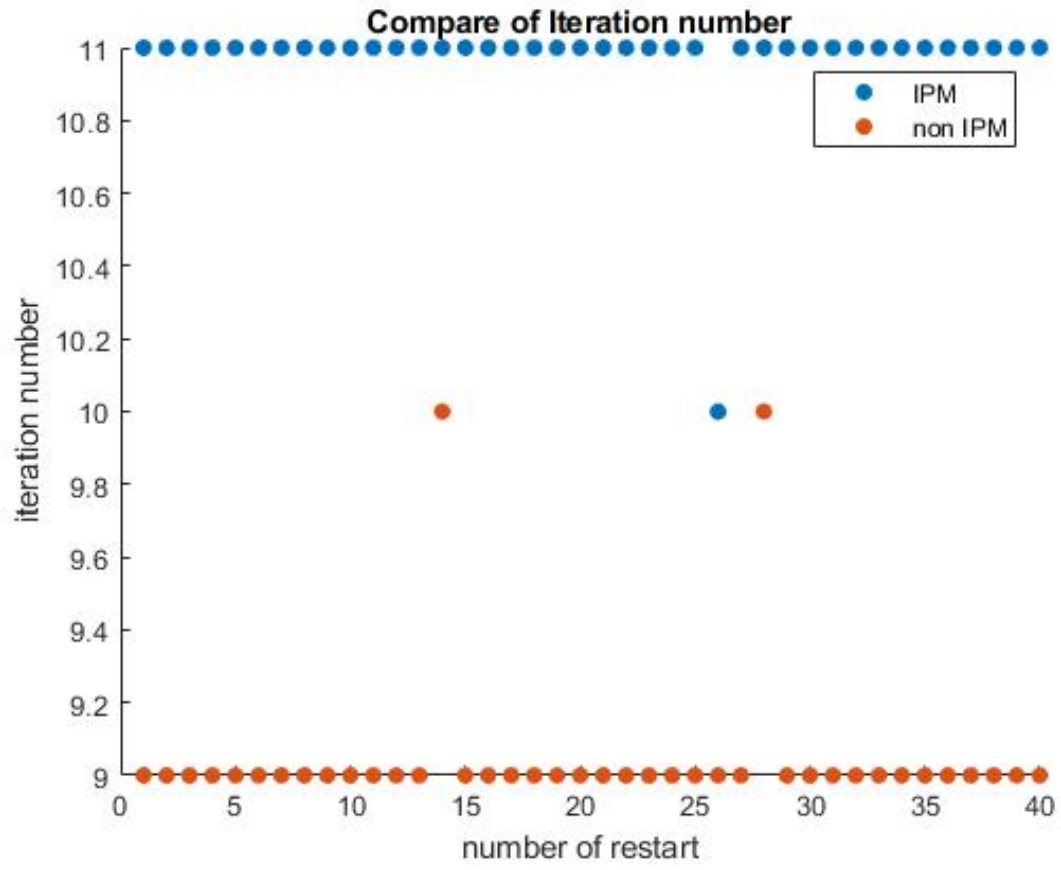

(e) $N = 500$

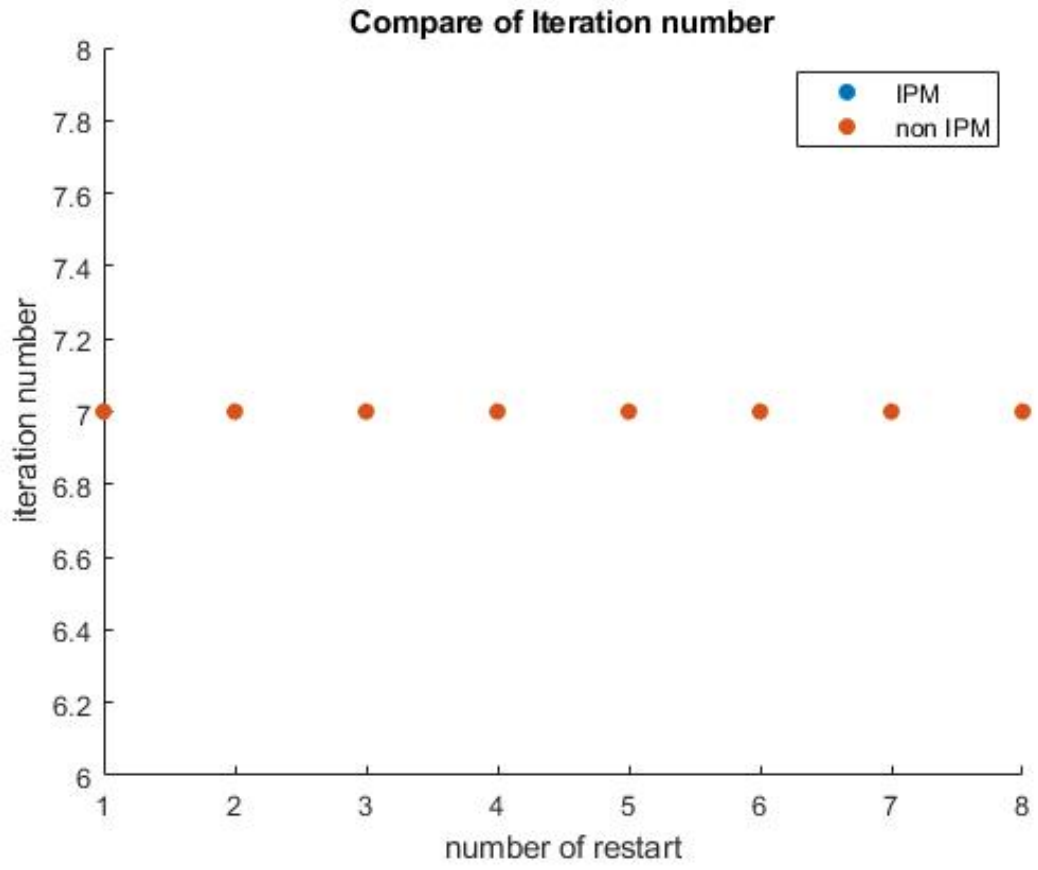

(f) $N = 1500$

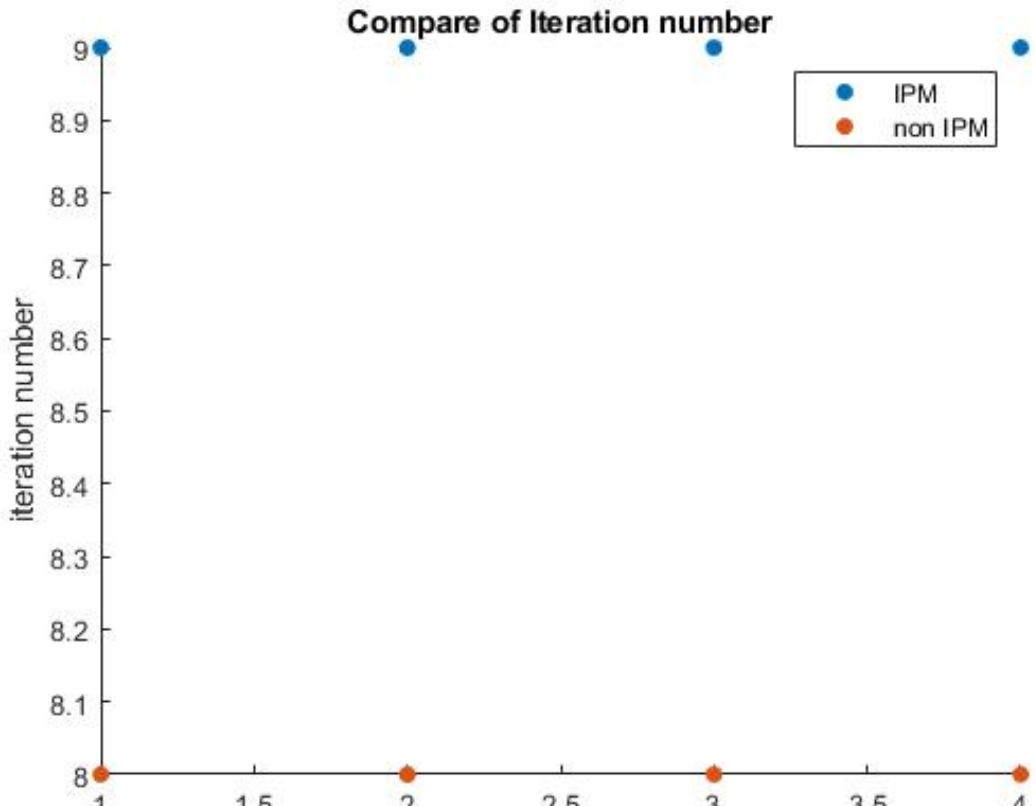

## 3.3 Time

**The time to run small cases is negligible so we skip it. We move to the cases when the size of the system is large. Here I don't compare Standard Methods because it diverge.**

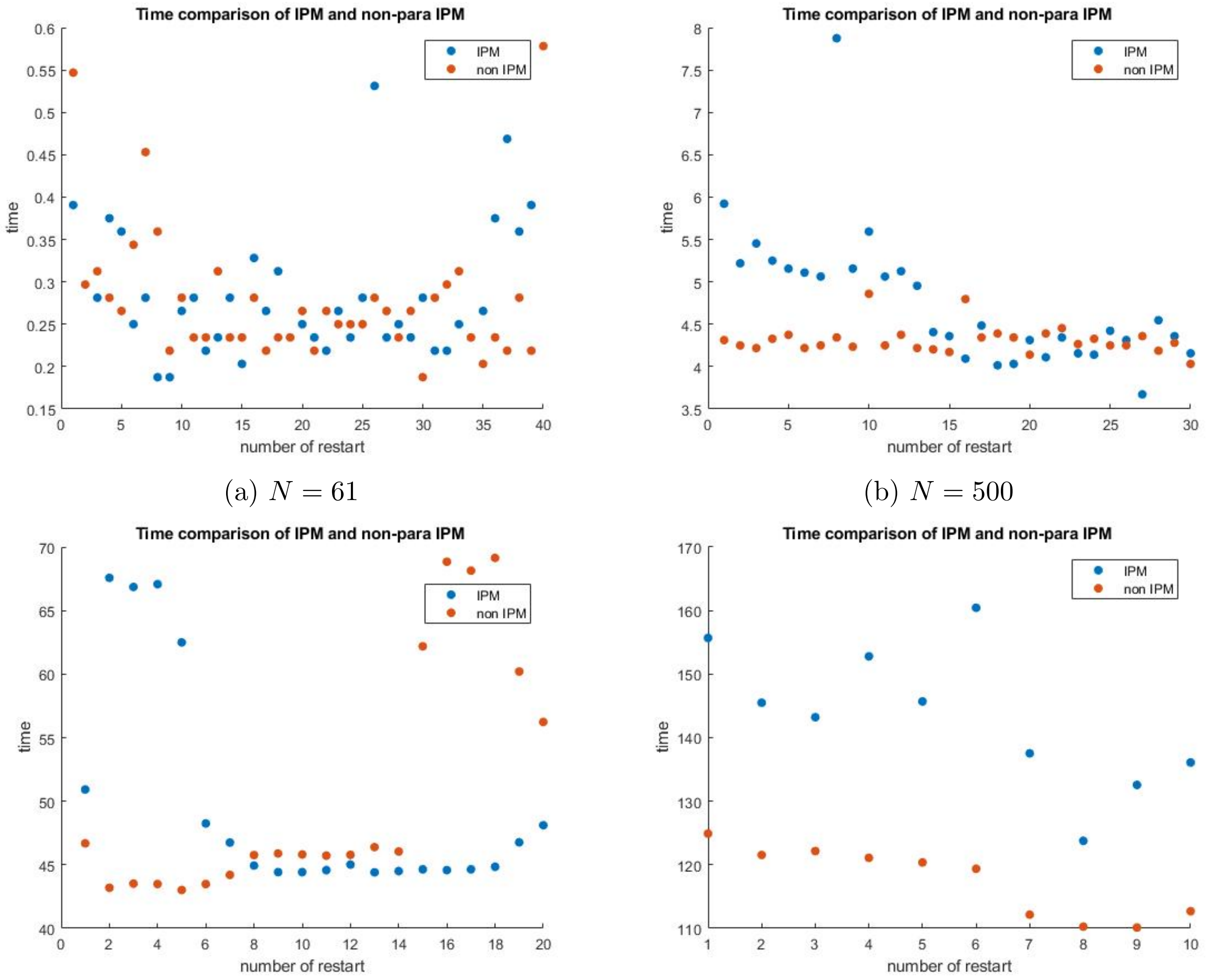

(a) $N = 61$ (b) $N = 500$

(c) $N = 1500$ (d) $N = 2000$

Figure 3.5: Compare Time of IPM and Non-Parametric IPM for various N (s)

# Chapter 4

# Conclusion

1. After all theoretical and numerical results,we can conclude that new algorithm **"non-parametric Interior Point Method"** is a promising method in solving system of equations containing complementarity conditions. It has a remarkable advantage when comparing to Interior Point Methods that we do not mind the strategy to enforce relax parameter for complementarity equations. We can say that Standard Methods are not converged the case N = 61,500,1500 ,2000(in the sense that they decrease the value of error slowly so that they meet the "Maximum iteration" condition ), but it still convergence to a true solution when we increase the value of $iter_{max}$. So as a remedy, IPM, and non-parametric IPM will converge to a true solution (in the sense that the number of iteration it takes to convergence is less than the maximum iteration).

2. we can see that when N is small, Standard methods work less iteration than IPM and non-parametric IPM, It has an advantage than IPM and non-para IPM. We can use this alongside with IPM or non-parametric IPM when N is small.

3. we can see that in case N = 61, not differ too much between IPM and non-parametric IPM. Some cases when N is small, IPM work with less iteration than Non-parametric IPM, but when N is large, then non-parametric IPM work more efficient than IPM

4. we can see that in most cases, non-parametric IPM work with less time consuming than IPM.

5. After all The things I listed above, in practice, we can conclude that NPIPM is better than Standard methods, IPM methods.

6. Some problem remain unsolved as of today, such as computational time.

7. In the remaining of the thesis, we will apply our new algorithms to several problems and verify its efficient.

# Chapter 5

# New Problem Application of non-parametric IPM

## 5.1 Pareto Eigenvalue

We define the set of Pareto eigenvalue

$$\lambda \in P_{sp}(A) \Leftrightarrow \{\exists x \geqslant 0/Ax - \lambda x \geqslant 0, < Ax - \lambda x, x >= 0\} \tag{5.1}$$

**We define new problems: Given a matrix A, compute one or several** $\lambda \in P_{sp}(A)$

We will solve the following form:

$$\begin{cases} \text{Find } \lambda > 0 \text{ and x } \in \mathbb{R}^n \setminus \{0\}\, s.t \\ x \geqslant 0, \lambda x - Ax \geqslant 0, < x, \lambda x - Ax >= 0 \end{cases} \tag{5.2}$$

Using non-parametric IPM we put slack variable

$$w = \lambda x - Ax \in \mathbb{R}^n \tag{5.3}$$

we have the system to solve is:

$$\begin{cases} \lambda x - Ax - w & & = 0 \\ x_i w_i & = \mu & \forall 1 \leqslant I \leqslant N \\ \|x\|_2 & & = 1 \\ f(x) + f(w) + \mu^2 + \epsilon\mu & & = 0 \end{cases} \tag{5.4}$$

where

$$f(t) = \frac{1}{2}\sum_{i=1}^{n} [min\,\{t_i, 0\}]^2\,, \forall t \in \mathbb{R}^n$$

**Remark 5.1.** *Let $\alpha > 0$, we notice that if $x$ a solution of (5.1) $\Leftrightarrow$ $\alpha x$ is also a solution of (5.1), choose $\alpha = \frac{1}{\|x\|_1}$, we can assume that $\|x\|_1 = 1$. Here we want to find $x \geqslant 0$, so the 1-norm $\|x\|_1 = |x_1| + |x_2| + ... + |x_n| = x_1 + x_2 + ... + x_n$ which is simpler than the 2-norm.*

*But with that convenience, there are many shortcomings. We need to find $x \geqslant 0/Ax - \lambda x \geqslant 0, < Ax - \lambda x, x >= 0$, In theory sometimes we have $x_i = 0$ and $(Ax - \lambda x)_i > 0$, In this cases in practice we have $x_i$ near to 0 and we are in danger, our function is not differentiable. So to avoid it, we use 2-norm instead of 1-norm which make the function always differentiable.*

then we get the Jacobian as follow :

$$DF(X) = \begin{bmatrix} \widehat{A} & -I \\ C & D \end{bmatrix} \in \mathbb{R}^{2m+2,2m+2} \tag{5.5}$$

m is the size of A , $\widehat{A} \in \mathbb{R}^{m,m}$ , $-I \in \mathbb{R}^{m,m}$, $I$ is the identity matrix .

$$\widehat{A} = \begin{bmatrix} \lambda - a_{11} & a_{12} & \cdots & a_{1m} \\ a_{21} & \lambda - a_{22} & \cdots & a_{2n} \\ \vdots & \vdots & \ddots & \vdots \\ a_{m1} & a_{m2} & \cdots & \lambda - a_{mm} \end{bmatrix} \in \mathbb{R}^{m,m} \tag{5.6}$$

$$C = \begin{bmatrix} w_1 & 0 & 0 & ... & 0 & 0 \\ 0 & w_2 & 0 & ... & 0 & 0 \\ 0 & 0 & w_3 & 0 & ... & 0 \\ ... \\ ... \\ ... \\ 0 & 0 & ... & w_{3N-2} & 0 & 0 \\ 0 & 0 & ... & 0 & w_{3N-1} & 0 \\ 0 & 0 & ... & 0 & 0 & w_{3N} \end{bmatrix} \in \mathbb{R}^{m,m}, D = \begin{bmatrix} u_1 & 0 & 0 & ... & 0 & 0 \\ 0 & u_2 & 0 & ... & 0 & 0 \\ 0 & 0 & u_3 & 0 & ... & 0 \\ ... \\ ... \\ ... \\ 0 & 0 & ... & u_{3N-2} & 0 & 0 \\ 0 & 0 & ... & 0 & u_{3N-1} & 0 \\ 0 & 0 & ... & 0 & 0 & u_{3N} \end{bmatrix} \in \mathbb{R}^{m,m} \tag{5.7}$$

Let $D\Im(\chi)$ be the Jacobian matrix of $\Im$. If $w_I \geqslant 0$ and $u_I \geqslant 0, \forall \quad 1 \leqslant I \leqslant m$, then we have Jacobian matrix of F

$$DF(\chi,\mu) = \begin{bmatrix} D\Im(\chi) \\ l \\ ll \end{bmatrix} \in \mathbb{R}^{2m+2,2m}, JacobiF = \begin{bmatrix} DF(\chi,\mu) & e & ee \end{bmatrix} \in \mathbb{R}^{2m+2,2m+2} \tag{5.8}$$

l $= (1, 1, ..., 1, 0, ..., 0) \in \mathbb{R}^{1,2m}$ first m components is 1, components from m+1 to end = 0, ll $=(0, 0, ..., 0, 0, ..., 0) \in \mathbb{R}^{1,2m}$, e $= (x_1, x_2, ..., x_m, 0, 0, 0, 0, ..., 0, ..., 0, 0, 0)^T \in \mathbb{R}^{2m+2,1}$ with m component above is the vector x, component m+1 to 2m+2 = 0, ee $= (0, ..., 0, 0, 0, -1, -1, ..., -1, ..., -1, 0, 2\mu + \epsilon)^T \in \mathbb{R}^{2m+2,1}$ with m component above is zeros, component m+1 to 2m = -1. $D\Im(\chi) \in \mathbb{R}^{2m,2m}$ is as in exercise before.

### 5.1.1 Numerical test of non-parametric IPM to find Pareto Eigenvalue

We choose the following matrix A which is known to have 23 Pareto eigenvalue.

$$A = \begin{bmatrix} 179 & -179 & -52 & 72 \\ -160 & 216 & -44 & 61 \\ -97 & -92 & 341 & 37 \\ -77 & -73 & -21 & 397 \end{bmatrix} \tag{5.9}$$

We then test the NPIPM method by applying to (5.9) using a sample of 900 random initial points. Table 2.1 presents the results obtained. The matrix (5.9) which is known to have 23 eigenvalues of Pareto. NPIPM finds each of them.

| Pareto Eigenvalue | Vector x | | | | Vector dual | | | |
|---|---|---|---|---|---|---|---|---|
| | $x_1$ | $x_2$ | $x_3$ | $x_4$ | $w_1$ | $w_2$ | $w_3$ | $w_4$ |
| $\lambda_1 = 0.5523$ | 0.4234 | 0.3589 | 0.2176 | 0 | 0 | 0 | 0 | 63.3766 |
| $\lambda_2 = 27.2583$ | 0.5412 | 0.4588 | 0 | 0 | 0 | 0 | 94.7060 | 75.1648 |
| $\lambda_3 = 29.9013$ | 0.3611 | 0.3061 | 0.1856 | 0.1472 | 0 | 0 | 0 | 0 |
| $\lambda_4 = 56.6292$ | 0.4433 | 0.3758 | 0 | 0.1809 | 0 | 0 | 70.8821 | 0 |
| $\lambda_5 = 152.2735$ | 0.6605 | 0 | 0.3395 | 0 | 0 | 120.6197 | 0 | 57.9888 |
| $\lambda_6 = 179.000$ | 1 | 0 | 0 | 0 | 0 | 160.000 | 97.000 | 77.000 |
| $\lambda_7 = 181.6219$ | 0.5204 | 0 | 0.2675 | 0.2121 | 0 | 82.0927 | 0 | 0 |
| $\lambda_8 = 189.3134$ | 0 | 0.6225 | 0.3775 | 0 | 131.0532 | 0 | 0 | 53.3682 |
| $\lambda_9 = 208.3947$ | 0.7101 | 0 | 0 | 0.2899 | 0 | 95.9312 | 58.1529 | 0 |
| $\lambda_{10} = 216.000$ | 0 | 1.000 | 0 | 0 | 179.0000 | 0 | 92.000 | 73.000 |
| $\lambda_{11} = 218.6354$ | 0 | 0.4791 | 0.2906 | 0.2303 | 84.2901 | 0 | 0 | 0 |
| $\lambda_{12} = 245.3669$ | 0 | 0.6750 | 0 | 0.3250 | 97.4314 | 0 | 50.0783 | 0 |
| $\lambda_{13} = 341$ | 0 | 0 | 1.000 | 0 | 52.000 | 44.000 | 0 | 21.000 |
| $\lambda_{14} = 366.3542$ | 0 | 0 | 0.5934 | 0.4066 | 1.5796 | 1.3053 | 0 | 0 |
| $\lambda_{15} = 367.6053$ | 0.2763 | 0 | 0 | 0.7237 | 0 | 0.0579 | 0.0215 | 0 |
| $\lambda_{16} = 367.6331$ | 0 | 0.2869 | 0 | 0.7131 | 0.0067 | 0 | 0.0074 | 0 |
| $\lambda_{17} = 367.6343$ | 0.0260 | 0.2599 | 0 | 0.7141 | 0 | 0 | 0.0052 | 0 |
| $\lambda_{18} = 367.6542$ | 0.0188 | 0.2363 | 0.0628 | 0.6821 | 0 | 0 | 0 | 0 |
| $\lambda_{19} = 367.6601$ | 0 | 0.2464 | 0.0819 | 0.6717 | 0.0032 | 0 | 0 | 0 |
| $\lambda_{20} = 367.6789$ | 0.1919 | 0 | 0.1772 | 0.6309 | 0 | 0.0184 | 0 | 0 |
| $\lambda_{21} = 367.6992$ | 0.1714 | 0 | 0.2204 | 0.6083 | 0 | 0.0100 | 0 | 0 |
| $\lambda_{22} = 367.7045$ | 0 | 0.1863 | 0.2036 | 0.6101 | 0.0020 | 0 | 0 | 0 |
| $\lambda_{23} = 367.7094$ | 0.0316 | 0.1463 | 0.2180 | 0.6041 | 0 | 0 | 0 | 0 |

Table 5.1: non-parametric IPM applied to the matrix (5.9) with 900 random initial points

**Here we test with another definition of Pareto Eigenvalue.**

Let $\mathbb{M}_n$ denote the space of real matrices of order $n$. To avoid trivialities we always assume that $n$ is greater than or equal to 2. A scalar $\lambda \in \mathbb{R}$ is called a Pareto eigenvalue of $A \in \mathbb{M}_n$ if the complementarity problem

$$x \geq \mathbf{0}_n, \quad Ax - \lambda x \geq \mathbf{0}_n, \quad \langle x, Ax - \lambda x \rangle = 0 \tag{5.10}$$

***We test non-parametric IPM on another concerns a small size matrix of special interest:***

$$A = \begin{bmatrix} 100 & 106 & -18 & -81 \\ 92 & 158 & -24 & -101 \\ 2 & 44 & 37 & -7 \\ 21 & 38 & 0 & 2 \end{bmatrix} \tag{5.11}$$

This matrix constructed by Y.-C. Qi (personal communication) which is known to have 23 Pareto eigenvalues. We tested the non-parametric IPM method on (5.11) with a sample of 900 random initial points. We get 23 Pareto eigenvalue when run NPIPM.

| Pareto Eigenvalue | vector x | | | | Vector dual | | | |
|---|---|---|---|---|---|---|---|---|
| | $x_1$ | $x_2$ | $x_3$ | $x_4$ | $w_1$ | $w_2$ | $w_3$ | $w_4$ |
| $\lambda_1 = 26.2823$ | 0.4314 | 0.0762 | 0 | 0.4924 | 0 | 0 | 0.7693 | 0 |
| $\lambda_2 = 26.4149$ | 0.4558 | 0.0368 | 0.0581 | 0.4493 | 0 | 0 | 0 | 0 |
| $\lambda_3 = 28.7114$ | 0.4527 | 0 | 0.1913 | 0.3559 | 0 | 1.1099 | 0 | 0 |
| $\lambda_4 = 29.1341$ | 0.2265 | 0.2491 | 0 | 0.5242 | 0 | 0 | 7.7457 | 0 |
| $\lambda_5 = 32.6080$ | 0 | 0.4461 | 0 | 0.5538 | 2.4260 | 0 | 15.7525 | 0 |
| $\lambda_6 = 32.8635$ | 0 | 0.2844 | 0.2897 | 0 | 3.0862 | 0 | 0 | 0 |
| $\lambda_7 = 37.5767$ | 0.2238 | 0 | 0.7761 | 0 | 0 | 1.9626 | 0 | 4.7001 |
| $\lambda_8 = 41.0162$ | 0.1241 | 0.0680 | 0.8077 | 0 | 0 | 0 | 0 | 5.194 |
| $\lambda_9 = 46.4681$ | 0 | 0.1770 | 0.8229 | 0 | 3.9579 | 0 | 0 | 6.7290 |
| $\lambda_{10} = 49.1435$ | 0.1561 | 0.1588 | 0.4874 | 0.1975 | 0 | 0 | 0 | 0 |
| $\lambda_{11} = 66.9700$ | 0 | 0.3428 | 0.4565 | 0.2005 | 11.8834 | 0 | 0 | 0 |
| $\lambda_{12} = 77.4251$ | 0.7814 | 0 | 0.0009 | 0.2175 | 0 | 49.8943 | 0 | 0 |
| $\lambda_{13} = 77.4575$ | 0.7822 | 0 | 0 | 0.2177 | 0 | 49.9815 | 0.0405 | 0 |
| $\lambda_{14} = 99.4233$ | 0.9689 | 0 | 0.0310 | 0 | 0 | 88.3988 | 0 | 20.3480 |
| $\lambda_{15} = 100.0000$ | 1 | 0 | 0 | 0 | 0 | 92 | 2 | 21 |
| $\lambda_{16} = 107.5010$ | 0 | 0.5132 | 0.3019 | 0.1848 | 33.9921 | 0 | 0 | 0 |
| $\lambda_{17} = 127.3920$ | 0 | 0.7674 | 0 | 0.2325 | 62.5094 | 0 | 32.1389 | 0 |
| $\lambda_{18} = 148.5319$ | 0 | 0.7170 | 0.2829 | 0 | 70.9203 | 0 | 0 | 27.2497 |
| $\lambda_{19} = 158.0000$ | 0 | 1 | 0 | 0 | 106 | 0 | 44 | 38 |
| $\lambda_{20} = 197.1730$ | 0.3414 | 0.4237 | 0.1154 | 0.1192 | 0 | 0 | 0 | 0 |
| $\lambda_{21} = 204.5836$ | 0.3874 | 0.4820 | 0 | 0.1305 | 0 | 0 | 21.06940 | 0 |
| $\lambda_{22} = 226.2813$ | 0.3934 | 0.4887 | 0.1177 | 0 | 0 | 0 | 0 | 26.8356 |
| $\lambda_{23} = 231.9223$ | 0.4455 | 0.5544 | 0 | 0 | 0 | 0 | 25.2880 | 30.4260 |

Table 5.2: non-parametric IPM applied to the matrix (5.11) with 900 random initial points

**Remark 5.2.** *we note that in [6], they use Newton methods with* $10^3$ *ramdom initial point and they just get some of Pareto Eigenvalues, but with non-parametric IPM, we got all of thems. It very important remark for effective of non-parametric IPM.*

***We test non-parametric IPM on another concerns a small size matrix of special interest: we take our matrix A with N = 61. More precisely,***

$$A \in \mathbb{R}^{183,183} \tag{5.12}$$

In our problems. we can find one Pareto eigenvalue is: $\lambda_1 = 494.9844$

***Conclude:*** with NPIPM and with 2 example we found from the paper, we can get 23 Pareto eigenvalue and with the cases N = 61 we can get some of thems. It show that the effective of our algorithm, It make our algorithm is a promising algorithm.

## 5.2 Application of NPIPM to the new Problems

### 5.2.1 adding one more equation in 1-norm

In this section we add the following equaion:

$$\sum_{i=1}^{N} |u_i| - 1 = 0$$

We define a new problems is
Let $N \geqslant 1$ be an integer. Given

. a 3N $\times$ 3N - matrix A with real entries

we want to find
a 3N-vector u with real components and $\lambda \in \mathbb{R}$ s.t

$$\begin{cases} (Au-\lambda u)_1 & = 0 \\ (Au-\lambda u)_2 & = 0 \\ min\,\{u_3; (Au-\lambda u)_3\} & = 0 \\ (Au-\lambda u)_4 & = 0 \\ (Au-\lambda u)_5 & = 0 \\ min\,\{u_6; (Au-\lambda u)_6\} & = 0 \\ \dots & \\ \dots & \\ \dots & \\ (Au-\lambda u)_{3N-2} & = 0 \\ (Au-\lambda u)_{3N-1} & = 0 \\ min\,\{u_{3N}; (Au-\lambda u)_{3N}\} & = 0 \end{cases} \tag{5.13}$$

where $(Au - \lambda u)_i$ denotes the *i-th* components of $Au - \lambda u$. In other words

$$(Au - \lambda u)_i = \sum_{j=1}^{3N} A_{ij} u_j - \lambda u_i \quad \forall 1 \leqslant i \leqslant 3N \tag{5.14}$$

**Remark 5.3.** *System (5.13) can be seen as a modification of the linear system* $\mathbf{Au} - \lambda\mathbf{u} = \mathbf{0}$*, in which every thirds equation has been replaced by a complementarity equation.*

**Ex. 4 —** 1. Write down the Jacobian matrix and the Newton iteration for (5.13)

**Answer (Ex. 4) —** If we rearrange $(Au-b)_{3I-2}, (Au-b)_{3I-1}, (Au-b)_{3I} - w_I$ on top , I means the way we arrange the system (5.13) as follow

$$\begin{cases} (Au-\lambda u)_1 = 0 \\ (Au-\lambda u)_2 = 0 \\ (Au-\lambda u)_3 - w_1 = 0 \\ (Au-\lambda u)_4 = 0 \\ (Au-\lambda u)_5 = 0 \\ (Au-\lambda u)_6 - w_2 = 0 \\ ... \\ ... \\ (Au-\lambda u)_{3N-2} = 0 \\ (Au-\lambda u)_{3N-1} = 0 \\ (Au-\lambda u)_{3N} - w_N = 0 \\ u_3 w_1 = \mu \\ u_6 w_2 = \mu \\ ... \\ ... \\ u_{3(N-1)} w_{N-1} = \mu \\ u_{3N} w_N = \mu \\ \sum_{i=1}^{N} |u_i| - 1 = 0 \\ \frac{1}{2}\sum_{I=1}^{N} min\left\{u_{3I};0\right\}^2 + \frac{1}{2}\sum_{I=1}^{N} min\left\{w_I;0\right\}^2 + \mu^2 + \epsilon\mu = 0 \end{cases} \tag{5.15}$$

**Remark 5.4.** *Let $\alpha > 0$, we note that if $u$ is a solution then $\alpha u$ is also a solution of the system (5.13), choose $\alpha = \frac{1}{\|u\|}$, we can suppose that $\|u\| = 1$.*
The Jacobian as follow:

$$DF(X) = \begin{bmatrix} \widehat{A} & B \\ C & D \end{bmatrix} \in \mathbb{R}^{m+n,m+n} \tag{5.16}$$

m is the size of A, $\widehat{A} \in \mathbb{R}^{m,m}$, $-I \in \mathbb{R}^{m,m}$, $-I$ is the identity matrix.

$$\widehat{A} = \begin{bmatrix} a_{11}-\lambda & a_{12} & \cdots & a_{1m} \\ a_{21} & a_{22}-\lambda & \cdots & a_{2n} \\ \vdots & \vdots & \ddots & \vdots \\ a_{m1} & a_{m2} & \cdots & a_{mm}-\lambda \end{bmatrix} \in \mathbb{R}^{m,m} \tag{5.17}$$

$$E = \begin{bmatrix} 0 & 0 & -1 & ... & 0 & 0 & & & \\ 0 & 0 & 0 & 0 & 0 & -1 & 0... & 0 & 0 \\ ... & & & & & & & & \\ ... & & & & & & & & \\ 0 & 0 & ... & 0 & 0 & 0 & 0 & 0 & 0 \\ 0 & 0 & ... & 0 & 0 & 0 & 0 & 0 & 0 \\ 0 & 0 & ... & 0 & 0 & 0 & 0 & 0 & -1 \end{bmatrix} \in \mathbb{R}^{n,m} \tag{5.18}$$

and **B = E'**, B $\in \mathbb{R}^{m,n}$(for simplification purpose only, because B is the same as C just take the transpose).

$$C = \begin{bmatrix} 0 & 0 & w_1 & ... & 0 & 0 & & & \\ 0 & 0 & 0 & 0 & 0 & w_2 & 0... & 0 & 0 \\ ... & & & & & & & & \\ ... & & & & & & & & \\ 0 & 0 & ... & 0 & 0 & 0 & 0 & 0 & 0 \\ 0 & 0 & ... & 0 & 0 & 0 & 0 & 0 & 0 \\ 0 & 0 & ... & 0 & 0 & 0 & 0 & 0 & w_N \end{bmatrix} \in \mathbb{R}^{n,m} \tag{5.19}$$

$$D = \begin{bmatrix} u_3 & 0 & 0 & ... & 0 & 0 \\ 0 & u_6 & 0 & ... & 0 & 0 \\ 0 & 0 & u_9 & 0 & ... & 0 \\ ... & & & & & \\ ... & & & & & \\ ... & & & & & \\ 0 & 0 & ... & u_{3N-6} & 0 & 0 \\ 0 & 0 & ... & 0 & u_{3N-3} & 0 \\ 0 & 0 & ... & 0 & 0 & u_{3N} \end{bmatrix} \in \mathbb{R}^{n,n} \tag{5.20}$$

.

Let $D\Im(\chi)$ be the Jacobian matrix of $\Im$. If $w_I \geqslant 0$ and $u_I \geqslant 0, \forall \quad 1 \leqslant I \leqslant m$, then we have Jacobian matrix of F

$$DF(\chi,\mu) = \begin{bmatrix} DF(X) \\ l \\ ll \end{bmatrix} \in \mathbb{R}^{m+n+2,m+n}, JacobiF = \begin{bmatrix} DF(\chi,\mu) & e & ee \end{bmatrix} \in \mathbb{R}^{m+n+2,m+n+2} \tag{5.21}$$

l $= (1,1,...,1,0,...,0) \in \mathbb{R}^{1,m+n}$ first m components is 1, components from m+1 to end = 0, ll $=(0,0,...,0,0,...,0) \in \mathbb{R}^{1,m+n}$, e $= (-u_1,-u_2,...,-u_m,0,0,0,0,...,0,...,0,0,0)^T \in \mathbb{R}^{m+n+2,1}$ with m component above is the vector x, component m+1 to m+n+2 = 0, ee $= (0,...,0,0,0,-1,-1,...,-1,...,-1,0,2\mu+\epsilon)^T \in \mathbb{R}^{m+n+2,1}$ with m component above is zeros, component m+1 to m+n = -1.

### 5.2.2 Numerical test of New Problems with adding 1-norm

In this section, we are not sure to get the convergence. $Tol$ is set to be 1e-6.

| Option N | Number of Initial points ( data in the box is the error ) | | | | | | | |
|---|---|---|---|---|---|---|---|---|
| | 1 | 2 | 3 | 4 | 5 | 6 | 7 | 8 |
| $N = 2$ | 7.7e-7 | 1.1e-7 | 2.1e-7 | 9.8e-7 | 4.1e-7 | 2.6e-7 | 9.3e-7 | 2.5e-7 |
| $N = 3$ | 4.5e-8 | 5.8 e-7 | 9.6e-8 | 3.7e-7 | 6.8e-7 | 4.5e-7 | 9.2e-7 | 9.0e-7 |
| $N = 15$ | 2.1e-7 | 1.7e-7 | 1.8e-8 | 2.2e-7 | 1.6e-7 | 2.4e-7 | 2.1e-7 | 1.9e-7 |
| $N = 61$ | 2.6e-08 | 8.8e-7 | 7.4e-7 | 1.0e-7 | 8.8e-7 | 5.7e-8 | 5.5e-7 | 3.17e-7 |
| $N = 500$ | 6.6e-7 | 5.5e-7 | 5.1e-7 | 3.2e-7 | 1.4e-7 | 1.0e-7 | 7.8e-7 | 7.3e-7 |
| $N = 1500$ | 1.4e-7 | 4.4e-8 | 1.3e-7 | 8.3e-7 | 2.4e-6 | 3.2e-04 | 5.5e-7 | 8.8e-7 |

Table 5.3: NPIPM applied to the New Problems (5.15) with 8 random initial points

The table (5.3) show the convergence behaviour of our problems.

| Option N | Number of Initial points ( data in the box is iteration ) | | | | | | | |
|---|---|---|---|---|---|---|---|---|
| | 1 | 2 | 3 | 4 | 5 | 6 | 7 | 8 |
| $N = 2$ | 27 | 23 | 25 | 24 | 24 | 25 | 23 | 25 |
| $N = 3$ | 15 | 16 | 17 | 29 | 28 | 16 | 26 | 16 |
| $N = 15$ | 15 | 15 | 16 | 16 | 15 | 16 | 16 | 17 |
| $N = 61$ | 24 | 33 | 30 | 25 | 20 | 24 | 21 | 29 |
| $N = 500$ | 54 | 60 | 59 | 76 | 65 | 64 | 63 | 78 |
| $N = 1500$ | 93 | 85 | 89 | 87 | 100 | 100 | 86 | 99 |

Table 5.4: NPIPM applied to the New Problems (5.15) with 8 random initial points

The table (5.5) show the iteration of our problems.

***Conclude:*** as you can see in the table, when the size of the Pbs is small, nothing happens, but when the size of the problems is large we met the trouble, as in the case N = 500, N = 1500, Sometimes we don't get the convergence.

**Remark 5.5.** *Here we put the 1-norm, but we can't sure that* $u_i \neq 0$*, when* $u_i = 0$ *our function is not differentiable. Especially, When we have the condition*

$$min\,\{u_{3I}; (Au-\lambda u)_{3I}\} = 0 \quad \forall 1 \leqslant I \leqslant N \tag{5.22}$$

*In theory, sometimes the solution of the above equation is* $u_{3I} = 0$*, and* $(Au-\lambda u)_{3I} \geqslant 0$*, so In practice, sometimes we have* $\exists I, 1 \leqslant I \leqslant N$ *s.t* $u_{3I} = \tau$*, which* $\tau \approx 0$ *so we are in danger, So to avoid it, we put the 2-norm. It serve the motivation for the next subsection.*

### 5.2.3 adding one more equation in 2-norm

In this section we add the following equaion:

$$\sum_{i=1}^{m} u_i^2 - 1 = 0$$

And the keep the same form of the previous ones, just change l = $l = (2u_1, 2u_2, ..., 2u_m, 0, ..., 0) \in \mathbb{R}^{1,m+n}$ first m components is the vector 2u, components from m+1 to end = 0.

We have the systems after replacing 1-norm equation by 2-norm equation as following

$$\begin{cases} (Au-\lambda u)_1 = 0 \\ (Au-\lambda u)_2 = 0 \\ (Au-\lambda u)_3 - w_1 = 0 \\ (Au-\lambda u)_4 = 0 \\ (Au-\lambda u)_5 = 0 \\ (Au-\lambda u)_6 - w_2 = 0 \\ \dots \\ \dots \\ (Au-\lambda u)_{3N-2} = 0 \\ (Au-\lambda u)_{3N-1} = 0 \\ (Au-\lambda u)_{3N} - w_N = 0 \\ u_3 w_1 = \mu \\ u_6 w_2 = \mu \\ \dots \\ \dots \\ u_{3(N-1)} w_{N-1} = \mu \\ u_{3N} w_N = \mu \\ \sum_{i=1}^{m} u_i^2 - 1 = 0 \\ \frac{1}{2}\sum_{I=1}^{N} min\left\{u_{3I};0\right\}^2 + \frac{1}{2}\sum_{I=1}^{N} min\left\{w_I;0\right\}^2 + \mu^2 + \epsilon\mu = 0 \end{cases} \tag{5.23}$$

### 5.2.4 Numerical test of New Problems with adding 2-norm

Firstly, we compare the iteration of the systems with 2-norm equation with the systems with 1-norm equation.

| Option N | Number of Initial points ( data in the box is iteration ) | | | | | | | |
|---|---|---|---|---|---|---|---|---|
| | 1 | 2 | 3 | 4 | 5 | 6 | 7 | 8 |
| $N = 2$ | (27,26) | (23,26) | (25,27) | (24,26) | (24,28) | (25,23) | (23,27) | (25,25) |
| $N = 3$ | (15,19) | (16,29) | (17,18) | (29,18) | (28,18) | (16,18) | (26,30) | (16,19) |
| $N = 15$ | (15,17) | (15,19) | (16,18) | (16,18) | (15,19) | (16,33) | (16,19) | (17,19) |
| $N = 61$ | (24,29) | (33,24) | (30,27) | (25,22) | (20,27) | (24,24) | (21,26) | (29,23) |
| $N = 500$ | (54,28) | (60,28) | (59,27) | (76,26) | (65,28) | (64,27) | (63,26) | (78,28) |
| $N = 1500$ | (93,30) | (85,34) | (89,29) | (87,30) | (100,34) | (100,32) | (86,34) | (99,31) |

Table 5.5: Compare 2-way (5.15) with 8 random initial points

where (5,7), means the iteration of the way 1 and the way 2 respectively.

**Conclude: We can see that the way 2 is more better than the way 1 in the sense of iteration. So we choose the way 2 is our ways to solve the problems.**

In this part, we use N = 61, and 1000 initial points to have as most as possible eigenvalue and eigenvector. The way we implement it on Matlab is

i) detect in the vector er that components is divergence, i.e. components $er_i$ s.t

$$er_i > \sqrt{\epsilon_machine}$$

by

$$er(find(er > \sqrt{\epsilon_machine})) = [];$$

ii) Transpose $final_sol = final_sol'$

iii) num_row = size(final_sol,1); num_column = size(final_sol,2);

iv) final_sol = sortrows(final_sol,num_column);

Then we remove the same row as follow:

**Algorithm 6** Remove the same row

position = [];
**for** i=1:num_row-1 **do**
    $b = norm(final_sol(i,:) - final_sol(i+1,:),1);$
    if $b <= 1e-3$
    $position = [position; i];$
final_sol(position,:)=[];
final_sol = sortrows(final_sol,size(final_sol,2));

Test with the cases N = 15, i.e. $A \in \mathbb{R}^{45,45}$ with 10000 difference initial point we get:

$$\lambda_1 = 10544.0490$$
$$\lambda_2 = 10557.6192$$
$$\lambda_3 = 13175.0569$$
$$\lambda_4 = 13175.0569$$
$$\lambda_5 = 13335.8262$$

Special with the $\lambda_3$ we got 2 eigenvector. Sometimes, thanks to the special form of A and u we get some opposite eigenvector, let call it u and v.

| | | | | |
|---|---|---|---|---|
| $u_1 = 0$ | $u_{10} = 0$ | $u_{19} = 0$ | $u_{28} = 0$ | $u_{37} = 0$ |
| $u_2 = 0.0719$ | $u_{11} = 0.2289$ | $u_{20} = 0.3650$ | $u_{29} = 0.3267$ | $u_{38} = 0.1782$ |
| $u_3 = 0$ | $u_{12} = 0$ | $u_{21} = 0$ | $u_{30} = 0$ | $u_{39} = 0$ |
| $u_4 = 0$ | $u_{13} = 0$ | $u_{22} = 0$ | $u_{31} = 0$ | $u_{40} = 0$ |
| $u_5 = 0.1268$ | $u_{14} = 0.2786$ | $u_{23} = 0.3924$ | $u_{32} = 0.2786$ | $u_{41} = 0.1269$ |
| $u_6 = 0$ | $u_{15} = 0$ | $u_{24} = 0$ | $u_{33} = 0$ | $u_{42} = 0$ |
| $u_7 = 0$ | $u_{16} = 0$ | $u_{25} = 0$ | $u_{34} = 0$ | $u_{43} = 0$ |
| $u_8 = 0.1782$ | $u_{17} = 0.3267$ | $u_{26} = 0.3650$ | $u_{35} = 0.2289$ | $u_{44} = 0.0720$ |
| $u_9 = 0$ | $u_{18} = 0$ | $u_{27} = 0$ | $u_{36} = 0$ | $u_{45} = 0$ |
| $v_1 = 0$ | $v_{10} = 0$ | $v_{19} = 0$ | $v_{28} = 0$ | $v_{37} = 0$ |
| $v_2 = -0.0719$ | $v_{11} = -0.2289$ | $v_{20} = -0.3650$ | $v_{29} = -0.3267$ | $v_{38} = -0.1782$ |
| $v_3 = 0$ | $v_{12} = 0$ | $v_{21} = 0$ | $v_{30} = 0$ | $v_{39} = 0$ |
| $v_4 = 0$ | $v_{13} = 0$ | $v_{22} = 0$ | $v_{31} = 0$ | $v_{40} = 0$ |
| $v_5 = -0.1268$ | $v_{14} = -0.2786$ | $v_{23} = -0.3924$ | $v_{32} = -0.2786$ | $v_{41} = -0.1269$ |
| $v_6 = 0$ | $v_{15} = 0$ | $v_{24} = 0$ | $v_{33} = 0$ | $v_{42} = 0$ |
| $v_7 = 0$ | $v_{16} = 0$ | $v_{25} = 0$ | $v_{34} = 0$ | $v_{43} = 0$ |
| $v_8 = -0.1782$ | $v_{17} = -0.3267$ | $v_{26} = -0.3650$ | $v_{35} = -0.2289$ | $v_{44} = -0.0720$ |
| $v_9 = 0$ | $v_{18} = 0$ | $v_{27} = 0$ | $v_{36} = 0$ | $v_{45} = 0$ |

we have if u satisfies (5.23) then,

$$\begin{cases} (Au-\lambda u)_1 & = 0 \\ (Au-\lambda u)_2 & = 0 \\ min\,\{u_3;(Au-\lambda u)_3\} & = 0 \\ (Au-\lambda u)_4 & = 0 \\ (Au-\lambda u)_5 & = 0 \\ min\,\{u_6;(Au-\lambda u)_6\} & = 0 \\ \dots & \\ \dots & \\ \dots & \\ (Au-\lambda u)_{3N-2} & = 0 \\ (Au-\lambda u)_{3N-1} & = 0 \\ min\,\{u_{3N};(Au-\lambda u)_{3N}\} & = 0 \end{cases} \tag{5.24}$$

now we prove v is also satisfies (5.23). Easy to see that v ≈ -u, so

$$\begin{cases} (Av-\lambda v)_1 & = (A(-u)-\lambda(-u))_1 = (\lambda u - Au)_1 = 0 \\ (Av-\lambda v)_2 & = (A(-u)-\lambda(-u))_2 = (\lambda u - Au)_2 = 0 \\ (Av-\lambda v)_4 & = (A(-u)-\lambda(-u))_4 = (\lambda u - Au)_4 = 0 \\ (Av-\lambda v)_5 & = (A(-u)-\lambda(-u))_5 = (\lambda u - Au)_5 = 0 \\ \dots & \\ \dots & \\ \dots & \\ (Av-\lambda v)_{3N-2} & = (A(-u)-\lambda(-u))_{3N-2} = (\lambda u - Au)_{3N-2} = 0 \\ (Av-\lambda v)_{3N-1} & = (A(-u)-\lambda(-u))_{3N-1} = (\lambda u - Au)_{3N-1} = 0 \end{cases} \tag{5.25}$$

it remains to prove that

$$min\,\{v_{3I};(Av-\lambda v)_{3I}\}=0,\forall \quad 1\leqslant I\leqslant N \tag{5.26}$$

It's simpler when we consider a special ones, consider

$$min\,\{v_3;(Av-\lambda v)_3\}=0 \tag{5.27}$$

we have $v_3=u_3=0$, so it just need to prove that $(Av-\lambda v)_3\geqslant 0$
It can happens because we look on the matrix A with position
A(3,2),A(3,5),A(3,8),A(3,11),A(3,14),...,A(3,44) is very near to 0 so when it multiply by $v_2, v_5, v_8, ...$, it is near to 0 also.
so

$$(Av-\lambda v)_3=(Av)_3=\sum_{i=1}^{N}A_{3,j}v_j=\sum_{i=2,5,8,...}^{N}A_{3,j}v_j\geqslant 0 \tag{5.28}$$

which complete the proof of v is also satisfies (5.23).
Test with the cases N = 61, i.e. $A\in\mathbb{R}^{183,183}$ with 3000 difference initial point we get:

| | | |
|---|---|---|
| $\lambda_1=404.8448$ | $\lambda_7=495.0281$ | $\lambda_{13}=555.1931$ |
| $\lambda_2=405.7169$ | $\lambda_8=530.4680$ | $\lambda_{14}=562.0067$ |
| $\lambda_3=473.1207$ | $\lambda_9=532.3509$ | $\lambda_{15}=588.1418$ |
| $\lambda_4=480.9800$ | $\lambda_{10}=534.0137$ | $\lambda_{16}=611.9478$ |
| $\lambda_5=488.4186$ | $\lambda_{11}=536.5415$ | $\lambda_{17}=665.6148$ |
| $\lambda_6=493.7288$ | $\lambda_{12}=540.8070$ | $\lambda_{18}=672.7544$ |

**Special with the $\lambda_{17}$ we got 2 eigenvector and two eigenvector is the same form of the cases N = 15,**
Test with the cases N = 500, i.e. $A\in\mathbb{R}^{1500,1500}$ with 1500 difference initial point we found only ones eigenvalue.

| |
|---|
| $\lambda_1=8.168$ |

Test with the cases N = 1500, i.e. $A\in\mathbb{R}^{4500,4500}$ with 200 difference initial point we get:

| | | |
|---|---|---|
| $\lambda_1=3523.7924$ | $\lambda_7=6661.7263$ | $\lambda_{13}=7928.6206$ |
| $\lambda_2=3528.2314$ | $\lambda_8=6783.2016$ | $\lambda_{14}=8811.4789$ |
| $\lambda_3=5482.0470$ | $\lambda_9=6920.2344$ | $\lambda_{15}=8895.7633$ |
| $\lambda_4=5562.2832$ | $\lambda_{10}=7054.4159$ | |
| $\lambda_5=6152.4462$ | $\lambda_{11}=7503.9343$ | |
| $\lambda_6=6167.3012$ | $\lambda_{12}=7923.3079$ | |

## 5.3 NPIPM when applied to to linear programing problems

We want to solve LP problems of the form:

$$\text{(LP)} \begin{cases} \min c^T x \\ Ax = b \\ x \geqslant \boldsymbol{O} \end{cases} \tag{5.29}$$

where $c \in \mathbb{R}^n, x \in \mathbb{R}^n, A \in \mathbb{R}^{m,n}, b \in \mathbb{R}^m$.

Lagrangian function

$$\mathcal{L}(x, \lambda, s) = c^T x + \lambda(Ax - b) - s^T x \tag{5.30}$$

where $\lambda \in \mathbb{R}^m, s \in \mathbb{R}^n_+$

***KKT conditions:*** $x^*$ is a solution of (5.34) if and only if $\exists(\lambda^*, s^*) \in (\mathbb{R}^m, \mathbb{R}^m_+)/$

$$\begin{cases} c + A^T \lambda^* - s^* = \boldsymbol{O} \\ Ax = b \\ < s, x >= \boldsymbol{O}, \end{cases} \tag{5.31}$$

Using "NPIPM" scheme, we have to solve

$$(\tilde{\text{P}}) \begin{cases} c + A^T \lambda^* - s^* = \boldsymbol{O} \\ Ax = b \\ s \cdot x - \mu e = \boldsymbol{O} \\ f(x) + f(s) + (\mu^2 + \epsilon\mu) = \boldsymbol{O} \end{cases} \tag{5.32}$$

where

$$f(t) = \frac{1}{2} \sum_{i=1}^{n} min \left\{t_i; 0\right\}^2, \forall t \in \mathbb{R}^n$$

If $s \geqslant \boldsymbol{O}_n$, $x \geqslant \boldsymbol{O}_n$ At each step, we need to solve the following linear system.

$$\begin{bmatrix} A & \boldsymbol{O} & \boldsymbol{O} & \boldsymbol{O} \\ s & x & \boldsymbol{O} & -e \\ \boldsymbol{O} & -I & A^T & \boldsymbol{O} \\ \boldsymbol{O} & \boldsymbol{O} & \boldsymbol{O} & 2\mu + \epsilon \end{bmatrix} \begin{bmatrix} x \\ s \\ \lambda \\ \mu \end{bmatrix} = \begin{bmatrix} b \\ \boldsymbol{O} \\ \boldsymbol{O} \\ \boldsymbol{O} \end{bmatrix} \tag{5.33}$$

### 5.3.1 Numerical test of NPIPM to Linear Programming

We test with the following example:

$$\text{(LP)} \begin{cases} \min x_1 - x_2 \\ x_1 + x_2 = 3 \\ x \geqslant \boldsymbol{O}, \end{cases} \tag{5.34}$$

The optimal solution is $(x_1, x_2) = (0, 3)$, we run NPIPM we got the solution $(x_1, x_2) = (0, 3)$ after 11 iteration.

Another test

$$\text{(LP)} \begin{cases} \min x_1 - x_2 \\ x_1 + 2x_2 = 3 \\ x \geqslant \boldsymbol{O}, \end{cases} \tag{5.35}$$

The optimal solution is $(x_1, x_2) = (0, \frac{3}{2})$, we run NPIPM we got the solution $(x_1, x_2) = (0, \frac{3}{2})$ after 10 iteration.

Another test

$$(\text{LP})\begin{cases} \min \ +6x_1+8x_2+5x_3+9x_4 \\ x_1+x_2+x_3+x_4=1 \\ x\geqslant \boldsymbol{O}, \end{cases} \tag{5.36}$$

The optimal solution is $(x_1,x_2,x_3,x_4) = (0,1,0,0)$, we run NPIPM we got the solution $(x_1,x_2,x_3,x_4) = (0,1,0,0)$ after 11 iteration.

Another test

$$(\text{LP})\begin{cases} \min \ +0x_1+0x_2+3x_3-x_4 \\ x_1-3x_3+3x_4=6 \\ x_2-8x_3+4x_4=4 \\ x\geqslant \boldsymbol{O} \end{cases} \tag{5.37}$$

The optimal solution is $(x_1,x_2,x_3,x_4) = (3,0,0,1)$, we run NPIPM we got the solution $(x_1,x_2,x_3,x_4) = (3,0,0,1)$ after 14 iteration.

Another test

$$(\text{LP})\begin{cases} \min \ +6x_1+14x_2+13x_3+0x_4+0x_5 \\ \frac{1}{2}x_1+2x_2+x_3+x_4=24 \\ x_1+x_2+4x_3+0x_4+x_5=60 \\ x\geqslant \boldsymbol{O} \end{cases} \tag{5.38}$$

The optimal solution is $(x_1,x_2,x_3,x_4,x_5) = (0,0,0,24,12)$, we run NPIPM we got the solution $(x_1,x_2,x_3,x_4,x_5) = (0,0,0,24,12)$ after 17 iteration.

To run more complicated example, we need to reformulated the form

$$(\text{LP})\begin{cases} \min \ c^Tx \\ Ax=b \\ lbounds\leqslant x\leqslant ubounds, \end{cases} \rightarrow (\tilde{\text{LP}})\begin{cases} \min \ c^Ty \\ Ay=b \\ y\geqslant \boldsymbol{O}, \end{cases} \tag{5.39}$$

we put x = w+l, so Ax = A(w+lbounds) = b so Aw = b-A $\times$ lbounds, and add the following equation

$$w+z=ubounds-lbounds \tag{5.40}$$

this equation come form $w\leqslant ubound-lbound$ so w+z = ubounds-lbounds which some $z\geqslant 0$, So we have

$$\hat{A}=\begin{bmatrix} A & 0 \\ I & I \end{bmatrix}, \hat{b}=\begin{bmatrix} (b-A\times loubnds) \\ ubounds-lbounds \end{bmatrix}, \hat{c}=\begin{bmatrix} c \\ 0 \end{bmatrix}, \hat{w}=\begin{bmatrix} w \\ z \end{bmatrix} \tag{5.41}$$

Finally, we get the form

$$(\tilde{\text{LP}})\begin{cases} \min \ \hat{c}^T\hat{w} \\ \hat{A}\hat{w}=\hat{b} \\ \hat{w}\geqslant \mathbf{0}, \end{cases} \tag{5.42}$$

We test NPIPM with some famous example:

| Option | algorithm | Iteration | Objective | Error |
|---|---|---|---|---|
| 25FV47 | NPIPM | 800 | 5.4986e+03 | 9.3124 |
| | PCLP | 191 | 5.5018e+03 | 3.3576e-07 |
| ADLITTLE | NPIPM | 198 | 2.2549e+05 | 1.3917e-08 |
| | PCLP | 64 | 2.2549e+05 | 1.7494e-07 |
| AFIRO | NPIPM | 38 | -464.7531 | 9.6721e-09 |
| | PCLP | 41 | -464.7531 | 4.0236e-08 |
| AGG | NPIPM | 464 | -3.5992e+07 | 4.9373e-09 |
| | PCLP | 145 | -3.5992e+07 | 5.0222e-08 |
| AGG2 | NPIPM | 557 | -2.0239e+07 | 8.9232e-09 |
| | PCLP | 80 | -2.0239e+07 | 3.2182e-08 |
| AGG3 | NPIPM | 534 | 1.0312e+07 | 6.1719e-09 |
| | PCLP | 102 | 1.0312e+07 | 2.9821e-08 |
| BANDM | NPIPM | 167 | -158.6280 | 1.7530e-09 |
| | PCLP | 60 | -158.6280 | 4.5697e-09 |
| BEACONFD | NPIPM | 111 | 3.3592e+04 | 2.8388e-09 |
| | PCLP | 36 | 3.3592e+04 | 1.1053e-08 |
| BLEND | NPIPM | 37 | -30.8121 | 3.2570e-10 |
| | PCLP | 20 | -30.8121 | 1.6794e-09 |
| BNL1 | NPIPM | 800 | 1.9695e+03 | 7.9993 |
| | PCLP | 271 | 1.9776e+03 | 4.0141e-09 |
| BNL2 | NPIPM | 800 | 1.8111e+03 | 1.2436 |
| | PCLP | 256 | 1.8112e+03 | 6.5168e-09 |
| BOEING1 | NPIPM | 800 | 571.3892 | 1519752 |
| | PCLP | 300 | 880.1620 | 95.27 |
| BOEING2 | NPIPM | 800 | -257.1791 | 3.15e+40 |
| | PCLP | 300 | -277.5541 | 0.4714 |
| BORE3D | NPIPM | 2 | 6.0093e-08 | 0 |
| | PCLP | 41 | 6.0093e-08 | 2.4225e-08 |
| BRANDY | NPIPM | 800 | 1.5185e+03 | 3.0e-09 |
| | PCLP | 52 | 1.5185e+03 | 3.6506e-10 |
| CAPRI | NPIPM | 240 | 1.9126e+03 | 2.1e-09 |
| | PCLP | 114 | 1.9126e+03 | 1.065e-09 |
| CYCLE | NPIPM | 800 | -3.8299e+12 | 90.37 |
| | PCLP | 189 | NaN | NaN |
| CZPROB | NPIPM | 800 | 1.9666e+06 | 31924 |
| | PCLP | 192 | 2.1825e+06 | 2.7435e-07 |
| D2Q06C | NPIPM | 800 | 1.1710e+05 | 8436.3 |
| | PCLP | 300 | 1.2253e+05 | 1.2108 |
| D6CUBE | NPIPM | 62 | 314.9167 | 1.1e-09 |
| | PCLP | 29 | 314.9167 | 4.0994e-14 |
| DEGEN2 | NPIPM | 800 | -1.4352e+03 | 2.41 |
| | PCLP | 300 | -1.4352e+03 | 4.1898e-08 |
| DEGEN3 | NPIPM | 800 | -987.2940 | 7.4363e-05 |
| | PCLP | 300 | -9.8729e+02 | 5.4086e-08 |
| FFFFF800 | NPIPM | 749 | 5.5568e+05 | 4.0866e-09 |
| | PCLP | 247 | 5.5568e+05 | 6.3554e-06 |
| ISRAEL | NPIPM | 800 | -4.8569e+05 | 3.5972e+05 |
| | PCLP | 35 | -6.3871e+05 | 1.4205e+03 |
| LOTFI | NPIPM | 800 | -10.7248 | 59.5059 |
| | PCLP | 300 | NaN | NaN |

Table 5.6: Comparison of NPIPM.m, PCLP.m for LP problems in Netlib

**Remark 5.6.** *PCLP here means Predictor Corrector algorithms for Linear Programming. Note that in all of famous cases above, we don't have pre-process and post-process of our matrix $A$ and vector $b$. So that the results is not exact to compare with any methods. Just to have an idea about our algorithm.*

***Conclude:*** We can see that in the cases that PCLP & NPIPM work, PCLP with less iteration than NPIPM such as ADLITTLE, AFIRO, AGG, AGG2, ... Moreover, in some cases NPIPM is diverge but PCLP is converge such as 25FV47, BLN1, BLN2, BRANDY. So that in this situation, we can conclude that PCLP is more effective than NPIPM for solving Linear Program.